\documentclass[a4paper,10pt]{article}
\usepackage{amsmath}
\usepackage{amsfonts}
\usepackage{amssymb}
\usepackage{array,longtable}
\usepackage{amscd}
\usepackage[cp1251]{inputenc}
\usepackage[english]{babel}
\usepackage[usenames]{color}
\usepackage[matrix,arrow,curve]{xy}
\usepackage{tikz}

\begin{document}
\newtheorem{definition}{Definition}
\newtheorem{theorem}{Theorem}
\newtheorem{lemma}{Lemma}
\newtheorem{proposition}{Proposition}
\newtheorem{criterion}{Criterion}
\newtheorem{observation}{Observation}
\newtheorem{corollary}{Corollary}
\newcommand{\diag}{{\rm diag}}
\newcommand{\tr}{{\rm Tr}}

\begin{center}
{\large {\bf UNIFYING MATRIX STABILITY CONCEPTS WITH A VIEW TO APPLICATIONS}} \\[0.5cm]
{\bf Olga Y. Kushel} \\[0.5cm]
Shanghai University, \\ Department of Mathematics, \\ Shangda Road 99, \\ 200444 Shanghai, China \\
kushel@mail.ru
\end{center}

\begin{abstract}
Multiplicative and additive $D$-stability, diagonal stability, Schur $D$-stability, $H$-stability are classical concepts which arise in studying linear dynamical systems. We unify these types of stability, as well as many others, in one concept of $({\mathfrak D}, {\mathcal G}, \circ)$-stability, which depends on a stability region ${\mathfrak D} \subset {\mathbb C}$, a matrix class ${\mathcal G}$ and a binary matrix operation $\circ$. This approach allows us to unite several well-known matrix problems and to consider common methods of their analysis. In order to collect these methods, we make a historical review, concentrating on diagonal and $D$-stability. We prove some elementary properties of $({\mathfrak D}, {\mathcal G}, \circ)$-stable matrices, uniting the facts that are common for many partial cases. Basing on the properties of a stability region $\mathfrak D$ which may be chosen to be a concrete subset of $\mathbb C$ (e.g. the unit disk) or to belong to a specified type of regions (e.g. LMI regions) we briefly describe the methods of further development of the theory of $({\mathfrak D}, {\mathcal G}, \circ)$-stability. We mention some applications of the theory of $({\mathfrak D}, {\mathcal G}, \circ)$-stability to the dynamical systems of different types.

Hurwitz stability, $D$-stability, diagonal stability, eigenvalue clustering, LMI regions, Lyapunov equation.

MSC[2010] 15A48, 15A18, 15A75
\end{abstract}
\tableofcontents
\part{Motivation and history}
\section{Introduction} Theory of dynamical systems gave rise to many classical matrix problems basing on Lyapunov's idea that stability of a system of ODE can be established through the spectral properties of the matrix of the system. In this connection, a number of matrix classes related to different stability types were introduced and studied. However, the further development of systems analysis including robustness and control led to new matrix problems. Rapid development in this area in the last decades caused some separation between classical matrix stability theory and its systems applications. Nowadays the main approaches to stability problems are the following ones.
\begin{enumerate}
\item[-] {\bf Matrix approach.} Here, we refer to the review papers by Hershkowitz \cite{HER1} and Datta \cite{DATTA1} as well as to a number of papers by Berman (see \cite{BERH}-\cite{BERVW}), Carlson (\cite{CARL1}-\cite{CAS}), C. Johnson (\cite{JOHN6}-\cite{JOHN2}), Hershkowitz and Schneider (\cite{HER1} - \cite{HERSS2}, \cite{SCHNE}) and many others. This line mainly focuses on long-standing open problems (for the examples, see \cite{HERK2}, \cite{HOLS}), connected to the classes of structured matrices, introduced in 1950s--1970s, classical methods of matrix analysis and graph theory (see \cite{BERH}).
\item[-] {\bf Systems and control approach.} Here, we refer to the special monograph by Kaszkurevicz and Bhaya \cite{KAB} which collects the results on different types of diagonal stability and their applications. This research line also includes modern books in robust control theory (see \cite{BH}, \cite{BAR}, \cite{BHK2}, \cite{RUGH}), a number of books and papers describing different approaches to stability and related problems (see \cite{BGF}, \cite{GUT2}, \cite{HENG}, \cite{JONCK}, \cite{JU2}, \cite{KOG}, \cite{KUSHN}, \cite{MART}), including modern and classical polynomial methods (see \cite{JEM}, \cite{MAR}, \cite{RAS}). This research states problems connected to the study of different types of dynamical systems and outlines some ways of their treatment.
\end{enumerate}

In this paper, we refresh the link between these two approaches. Besides surveying the existing results, we provide the unified proofs of some properties common to many different types of stability.

\subsection{Overview of Paper} The definitions of matrix classes, mentioned in this paper, are given in the Appendix, made in the form of a matrix dictionary. The reader not familiar with these definitions may consult it when necessary. The definitions, connected to matrix minors and compound matrices, can also be found there.

Part I of the paper consists of Sections 1-5, where we give motivation of unifying several concepts of matrix stability under a general one and present a historical overview of the most studied partial cases. Section 1 gives introduction, basic definitions and examples of new results obtained using the concept introduced in the paper. We collect the examples of partial cases of $({\mathfrak D}, {\mathcal G}, \circ)$-stability which appeared in the literature, describing the corresponding stability regions $\mathfrak D$, matrix classes $\mathcal G$ and binary operations $\circ$. In Section 2, we list several results on matrix (polynomial) stability and generalized stability, also known as eigenvalue clustering in a given region of ${\mathbb C}$. We group this results according to the methods we will use later. In Section 3, we present the detailed survey of the historical development of the main important partial cases of $({\mathfrak D}, {\mathcal G}, \circ)$-stability, namely, multiplicative and additive $D$-stability and diagonal stability. In Section 4, we consider binary operations and their properties, stating some open problems. The properties of binary operations connected to the matrix spectra form a basis for the further results. Section 5 deals with matrix classes and their properties, mainly focusing on the properties of symmetric positive definite matrices. Positive (negative) definiteness plays a crucial role in studying different stability types.

Part II of the paper consists of Sections 6-10. The main trust of Part II is the unified proof of basic properties most cases do share. In Section 6, we provide the proof of inclusion relations and topological properties, using the methods of abstract algebra and group theory. Section 7 deals with the unified proof of basic properties (transposition, inversion, multiplication by a scalar) of $({\mathfrak D},{\mathcal G},\circ)$-stable matrices. In Section 8, we study so-called "Lyapunov regions" described by generalizations of the Lyapunov theorem. We define a generalization of the widely used concept of diagonal stability and study its relations to $({\mathfrak D},{\mathcal G},\circ)$-stability. Section 9 deals with the qualitative approach and its possible generalizations. In this section, we also list other methods of studying $({\mathfrak D},{\mathcal G},\circ)$-stability. Section 10 provides the ways of the further development of $({\mathfrak D},{\mathcal G},\circ)$-stability theory and outlines some open problems.

Part III consists of Sections 11-14 dealing with the existing and potential applications. In Section 11, we apply the above theory to the perturbed families of dynamical systems of different kinds. Section 12 deals with the global asymptotic stability of nonlinear systems and applications of diagonal stability and its generalizations. In Section 13, we consider the recent application of diagonal stability to passivity and network stability analysis. Section 14 describes classical dynamical models which give raise to the partial cases of $({\mathfrak D},{\mathcal G},\circ)$-stability, mentioned in the previous sections.

\subsection{Unifying concept} Let ${\mathcal M}^{n \times n}$ denote the set of all real $n \times n$ matrices and $\mathbf A \in {\mathcal M}^{n \times n}$. Let $\sigma({\mathbf A})$ denote the spectrum of $\mathbf A$ (i.e. the set of all eigenvalues of $\mathbf A$ defined as zeroes of its characteristic polynomial $f_{\mathbf A}(\lambda):= \det(\lambda{\mathbf I - {\mathbf A}})$).

{\bf ${\mathfrak D}$-stability.} Let ${\mathfrak D} \subset {\mathbb C}$ be a set with the property: $\lambda \in {\mathfrak D}$ implies its complex conjugate $\overline{\lambda} \in {\mathfrak D}$ (this property is needed since we study matrices with real entries). An $n \times n$ matrix ${\mathbf A}$ is called {\it stable with respect to $\mathfrak D$} or simply {\it $\mathfrak D$-stable} if $\sigma({\mathbf A}) \subset {\mathfrak D}$. In this case, ${\mathfrak D}$ is called a {\it stability region}. Note, that we do not impose any restrictions (e.g. connectivity or convexity) on ${\mathfrak D}$. Consider the following most well-known examples.
\begin{enumerate}
\item[\rm 1.] {\bf Stability.} An $n \times n$ real matrix $\mathbf A$ is called {\it Hurwitz stable (semistable)} or just {\it stable (semistable)} if all its eigenvalues have negative (respectively, nonpositive) real parts (see, for example, \cite{BELL}, \cite{KAB}, \cite{MAM}). In a number of books and papers in matrix theory, {\it positive stability} is used: $\mathbf A$ is called {\it positive stable} if all its eigenvalues have positive real parts.
\item[\rm 2.] {\bf Schur stability.} An $n \times n$ real matrix $\mathbf A$ is called {\it Schur stable} if all its eigenvalues lie inside the unit circle, i.e. the spectral radius $\rho(\mathbf A) < 1$ (see \cite{BHK}, \cite{KAB}). This property is mostly referred as {\it convergence of matrices} (see \cite{HOJ}, p. 137).
\item[\rm 3.] {\bf Aperiodicity.} An $n \times n$ real matrix $\mathbf A$ is called {\it aperiodic} if all its eigenvalues are real (see, for example, \cite{GUJU}, p. 860 and \cite{GUT2}, p. 92).
\end{enumerate}

{\bf $({\mathfrak D},{\mathcal G},\circ)$-stability}. Given a stability region ${\mathfrak D} \subset {\mathbb C}$, a matrix class ${\mathcal G} \subset {\mathcal M}^{n \times n}$ and a binary operation $\circ$ on ${\mathcal M}^{n \times n}$, we call an $n \times n$ matrix $\mathbf A$ {\it left (right) $({\mathfrak D},{\mathcal G},\circ)$-stable} if $\sigma({\mathbf G}\circ{\mathbf A}) \subset {\mathfrak D}$ (respectively, $\sigma({\mathbf A}\circ{\mathbf G}) \subset {\mathfrak D}$) for any matrix $\mathbf G$ from the class ${\mathcal G}$.  Later we will show that, for the most used binary operations $\circ$, the notion of left $({\mathfrak D},{\mathcal G},\circ)$-stability coincides with the notion of right $({\mathfrak D},{\mathcal G},\circ)$-stability. Thus we use the term "$({\mathfrak D},{\mathcal G},\circ)$-stable" if a matrix $\mathbf A$ is both left and right $({\mathfrak D},{\mathcal G},\circ)$-stable.

Consider the following $({\mathfrak D},{\mathcal G},\circ)$-stable matrix classes. They are grouped with a view to applications, with respect to:
\begin{enumerate}
\item[\rm 1.] {\bf Binary operation}. The choice of a binary operation $\circ$ represents the type of the perturbations of a dynamical system.
\item[\rm 2.] {\bf Stability region}. The choice of a region $\mathfrak D$ is defined by the type of a dynamical system (e.g. continuous-time, discrete-time, fractional order, etc) or by studied system properties (e.g. the transient response of a system, oscillation and Hopf bifurcation phenomena).  Later we mostly consider the concepts, introduced in the literature for the following stability regions: the open left-hand side of the complex plane
$${\mathbb C}^-:=\{\lambda \in {\mathbb C}: {\rm Re}(\lambda) < 0\}$$
and the interior of the unit disk.
$$D(0,1):= \{\lambda \in {\mathbb C}: |\lambda| < 1\}.$$
These regions correspond to continuous- and discrete-time linear systems, respectively.
\item[\rm 3.] {\bf Matrix class}. The choice of a matrix class $\mathcal G$ is determined by the properties of the analyzed dynamical model. Besides that, some matrix classes may be introduced "artificially" for analyzing the stability properties of some other types of system perturbations.
\end{enumerate}

\subsection{Multiplicative $({\mathfrak D},{\mathcal G})$-stability} Here, the operation $\circ$ is fixed as matrix multiplication.

First, we consider the {\bf stability region ${\mathfrak D} = {\mathbb C}^{-}$}.

{\bf Multiplicative $D$-stability (Arrow and McManus, 1958).} An $n \times n$ real matrix $\mathbf A$ is called (multiplicative) {\it $D$-stable} if ${\mathbf D}{\mathbf A}$ is stable for every positive diagonal matrix $\mathbf D$ (i.e., an $n \times n$ matrix $\mathbf D$ with positive entries on its principal diagonal, while the rest are zero). Later on, using the term "$D$-stability", we mostly refer to multiplicative $D$-stability. Here, ${\mathcal G}$ is the class of positive diagonal matrices. This concept was introduced in \cite{AM} in connection with some problems of mathematical economics. The literature on multiplicative $D$-stability is particularly rich due to a lot of applications (see Part III for the selected details, also see \cite{LOG}, \cite{QR}, \cite{KAB} and references therein).

{\bf $H$-stability (Arrow and McManus, 1958).} An $n \times n$ real matrix $\mathbf A$ is called (multiplicative) {\it $H$-stable} if ${\mathbf H}{\mathbf A}$ is stable for every symmetric positive definite matrix $\mathbf H$. Here, ${\mathcal G}$ is the class of symmetric positive definite matrices. This matrix class also arises in \cite{AM} under the name of $S$-stability and later studied (see \cite{CARL3}, \cite{CAS}, \cite{OSS}) under the name of $H$-stability.

The following concept "interpolates" stability and $D$-stability.

{\bf Multiplicative $D(\alpha)$-stability (Khalil and Kokotovic, 1979).} Let, as usual, $[n]$ denotes the set of indices $\{1, \ \ldots, \ n\}$. Given a positive integer $p$, $1 \leq p \leq n$, let $\alpha = (\alpha_1, \ \ldots \ \alpha_p)$ be a partition of $[n]$. A diagonal matrix $\mathbf D$ is called an {\it $\alpha$-scalar matrix} if ${\mathbf D}[\alpha_k]$ is a scalar matrix for every $k = 1, \ \ldots, \ p$, i.e.
    $${\mathbf D} = \diag\{d_{11} {\mathbf I}[\alpha_1], \ \ldots, \ d_{pp} {\mathbf I}[\alpha_p]\}.$$
(Here, as before, ${\mathbf D}[\alpha_k]$ denotes a principal submatrix spanned by rows and columns with indices from $\alpha_k$).
    ${\mathbf D}$ is called a {\it positive $\alpha$-scalar matrix} if, in addition, $d_{ii} > 0$, $i = 1, \ \ldots, \ p$.
    Khalil and Kokotovic introduce the following definition based on the given above matrix class (see \cite{KHAK2}, \cite{KHAK1} ). An $n \times n$ matrix $\mathbf A$ is called {\it $D(\alpha)$-stable} (relative to the partition $\alpha = (\alpha_1, \ \ldots \ \alpha_p)$) if ${\mathbf D}{\mathbf A}$ is stable for every positive $\alpha$-scalar matrix ${\mathbf D}$. (Originally, this property was called {\it block $D$-stability}). The matrix class ${\mathcal G}$ for this case is the class of positive $\alpha$-scalar matrices.

This recently introduced concept "interpolates" $D$-stability and $H$-stability.

 {\bf $H(\alpha)$-stability (Hershkowitz and Mashal, 1998).} Given a positive integer $p$, $1 \leq p \leq n$, the set $\alpha = (\alpha_1, \ \ldots \ \alpha_p)$, where each $\alpha_i$ is a nonempty subset of $[n]$, $\alpha_i\bigcap\alpha_j =\emptyset$ and $[n] = \bigcup_i \alpha_i$, is called a {\it partition} of $[n]$. Without loss of generality, we may assume that each $\alpha_i$, $i = 1, \ \ldots, \ p$, consists of contagious indices. A block diagonal matrix $\mathbf H$ of the form
    $${\mathbf H} = \diag\{H[\alpha_1], \ \ldots, \ H[\alpha_p]\},$$ where each $H[\alpha_i]$ is a principal submatrix of ${\mathbf H}$ formed by rows and columns with indices from $\alpha_i$, $i = 1, \ \ldots, \ p$, is called an {\it $\alpha$-diagonal matrix} (see \cite{HERM}).
The following concept was provided in \cite{HERM}: given a partition $\alpha = (\alpha_1, \ \ldots \ \alpha_p)$, an $n \times n$ real matrix $\mathbf A$ is called (multiplicative) {\it $H(\alpha)$-stable} if ${\mathbf H}{\mathbf A}$ is stable for every symmetric positive definite $\alpha$-diagonal matrix $\mathbf H$. Here, ${\mathcal G}$ is the class of symmetric positive definite $\alpha$-diagonal matrices denoted by $H(\alpha)$.

{\bf Interval $D$-stability (Johnson, 1975 and Kosov, 2010).} Consider a matrix parallelepiped of the form:
 $$ \Theta = \diag\{d_{ii}, \ \ 0 < d_{ii}^{min} < d_{ii} < d_{ii}^{max} < + \infty, \ \  i  = 1, \ \ldots, \ n\}.$$
An $n \times n$ matrix $\mathbf A$ is called {\it $D$-stable with respect to $\Theta \subset {\mathcal M}^{n \times n}$} if ${\mathbf D}{\mathbf A}$ is stable for every matrix ${\mathbf D} \in \Theta$. This concept fist appeared in a remark in \cite{JOHN5} by Johnson under the name of "partial $D$-stability". However, the results based on this concept appeared only recently (see \cite{KOS}, where this concept was independently introduced).

{\bf Ordered $D$-stability (Kushel, 2016).} The following definition was provided in \cite{KU1}: given a positive diagonal matrix ${\mathbf D} = \diag\{d_{11}, \ \ldots, \ d_{nn}\}$ and a permutation $\tau = (\tau(1), \ \ldots, \ \tau(n))$ of the set of indices $[n]:= \{1, \ \ldots, \ n\}$, we call the matrix ${\mathbf D}$ {\it ordered with respect to $\tau$}, or {\it $\tau$-ordered}, if it satisfies the inequalities $$d_{\tau(i)\tau(i)} \geq d_{\tau(i+1)\tau(i+1)}, \qquad i = 1, \ \ldots, \ n-1. $$ A matrix $\mathbf A$ is called {\it $D$-stable with respect to the order $\tau$}, or {\it $D_\tau$-stable}, if the matrix ${\mathbf D}{\mathbf A}$ is positive stable for every $\tau$-ordered positive diagonal matrix $\mathbf D$.

Note, that if we consider the closed stability region $\overline{{\mathfrak D}} = \{\lambda \in {\mathbb C}: {\rm Re}(\lambda)\leq 0\}$, we shall obtain the corresponding concepts of {\it $D$-semistability}, {\it $H$-semistability} and others.

Now consider the stability region ${\mathfrak D} = D(0,1)$, i.e. the {\bf interior of the unit disk}.

{\bf Schur $D$-stability (Bhaya and Kaszkurewicz, 1993).} An $n \times n$ real matrix $\mathbf A$ is called {\it Schur $D$-stable} if ${\mathbf D}{\mathbf A}$ is Schur stable for every diagonal matrix $\mathbf D$ with $\|{\mathbf D}\|< 1$ (i.e. an $n \times n$ diagonal matrix $\mathbf D$ with $|d_{ii}|< 1$, $i = 1, \ \ldots, \ n$). Here ${\mathfrak D} = \{z \in {\mathbb C}: |z|< 1\}$, ${\mathcal G}$ is the class of diagonal matrices with $\|{\mathbf D}\|< 1$. This matrix class was defined in \cite{BHK} in connection with the study of discrete-time systems, and studied in \cite{CA2} as {\it convergent multiples}.

{\bf Vertex stability (Bhaya and Kaszkurewicz, 1993).} A matrix $\mathbf A$ is called {\it vertex stable} if $\rho({\mathbf D}{\mathbf A}) < 1$ for any real diagonal matrix $\mathbf D$ with $|{\mathbf D}|=1 $, i.e. with $d_{ii} = \pm1, \ i = 1, \ \ldots, \ n$ (note that this matrix class contains only a finite number of matrices). This matrix class was also defined in \cite{BHK}. The concept of vertex stability was introduced for characterizing Schur $D$-stability. The paper \cite{CA2} provides different proofs of the basic results on Schur $D$-stability and vertex stability.

{\bf Other stability regions.}
This two concepts may be considered as spectral localization inside/outside the boundary of a stability region.

{\bf $D$-hyperbolicity (Abed, 1986).} An $n \times n$ real matrix $\mathbf A$ is called {\it $D$-hyperbolic} if the eigenvalues of ${\mathbf D}{\mathbf A}$ have nonzero real parts for every real nonsingular $n \times n$ diagonal matrix $\mathbf D$. This definition was provided in \cite{AB1}, see also \cite{AB2}, in connection with studying Hopf bifurcation phenomena of linearized system of differential equations. In this case, ${\mathfrak D}$ is the complex plane without imaginary axes, ${\mathcal G}$ is the class of nonsingular diagonal matrices.

{\bf $D$-positivity and $D$-aperiodicity (Barkovsky, Ogorodnikova, 1987).} An $n \times n$ real matrix $\mathbf A$ is called (multiplicative) {\it $D$-positive ($D$-negative)} if all the eigenvalues of ${\mathbf D}{\mathbf A}$ are positive (respectively, negative) for every positive diagonal matrix $\mathbf D$. This definition was given in \cite{BAO} in connection with oscillation properties of mechanical systems. The stability region ${\mathfrak D}$ in this case is the positive direction of the real axes, ${\mathcal G}$ is the class of positive diagonal matrices. It is natural to define also the following matrix class: an $n \times n$ real matrix $\mathbf A$ is called (multiplicative) {\it $D$-aperiodic} if all the eigenvalues of ${\mathbf D}{\mathbf A}$ are real for every diagonal matrix $\mathbf D$. Here, we extend the stability region ${\mathfrak D}$ to the whole real axes and also extend the class ${\mathcal G}$ to the whole class of diagonal matrices from ${\mathcal M}^{n \times n}$.

{\bf More general concepts.}

{\bf ${\mathcal G}$-stability (Cain, DeAlba, Hogben, Johnson, 1998).} The following generalization of $D$-stability was defined in \cite{CADHJ} (see \cite{CADHJ}, p. 152). By varying the class ${\mathcal G} \subset {\mathcal M}^{n \times n}$, the following definitions were obtained: an $n \times n$ real matrix $\mathbf A$ is called {\it ${\mathcal G}$-stable} ({\it ${\mathcal G}$-convergent}) if ${\mathbf G}{\mathbf A}$ is positive stable (respectively, convergent) for every matrix ${\mathbf G}$ from the selected matrix class ${\mathcal G}$. The concept of "set product", when a matrix class ${\mathcal G}$ is multiplied by a given matrix $\mathbf A$ is further studied in \cite{BHK3} and \cite{CALN}.

{\bf Polyhedron stability (Geromel, de Oliveira, Hsu, 1998).}  Some generalization of Schur and vertex stability is studied in \cite{GEOH}: given a matrix ${\mathbf A} \in {\mathcal M}^{n \times n}$, and a convex polyhedron ${\mathcal B} \subset {\mathcal M}^{n \times n}$, defined as the convex hull of the finite number of its {\it extreme} matrices ${\mathbf B}_1, \ \ldots, \ {\mathbf B}_N$. The stability (with respect to a given region) of matrix set $${\mathcal A}: = \{{\mathbf A}{\mathbf B}: {\mathbf B} \in {\mathcal B}\} $$ is considered. In this case, the matrix class ${\mathcal G}$ we may consider as all the points of matrix polyhedron. The corresponding notion of {\it polyhedron vertex stability}, where the matrix class ${\mathcal G}$ consists of a finite number of the extreme matrices, may be used for characterizing polyhedron stability. Note, that originally this concept was introduced to cover the two most important stability regions: ${\mathbb C}^-$ and $D(0,1)$.

\subsection{Hadamard $({\mathfrak D},{\mathcal G})$-stability} Now we set the binary operation $\circ$ to be Hadamard (entry-wise) matrix multiplication (for the definitions and properties see \cite{JOHN3}). First, consider the stability region ${\mathfrak D} = {\mathbb C}^{-}$.

{\bf Hadamard $H$-stability (Johnson, 1974).} A real matrix ${\mathbf A}$ is called {\it Hadamard $H$-stable} if ${\mathbf H}\circ{\mathbf A}$ is stable for every symmetric positive definite matrix ${\mathbf H} \in {\mathcal M}^{n \times n}$. Here ${\mathcal G}$ is the class of symmetric positive definite matrices. This definition was introduced by Johnson (see \cite{JOHN3}, p. 304) and was called by the author {\it Schur stability}. But, since the term "Schur stable" is already reserved for matrices which spectral radius is less than 1, we refer to this property as to {\it Hadamard $H$-stability}.

{\bf Finite-rank Hadamard stability (Johnson and van den Driessche, 1988).} The following matrix classes were introduced in \cite{JOHD} in order to "interpolate" matrix properties and properties of sign pattern classes (for the definitions of sign-pattern classes and sign-stability, see the Appendix), in particular, the classes of $D$-stable and sign-stable matrices. An $n \times n$ real matrix $\mathbf A$ is called {\it $B_k$-stable} (to belong to the class $B_k$) if the Hadamard product ${\mathbf B}\circ{\mathbf A}$ is positive stable for every entry-wise positive matrix ${\mathbf B} \in {\mathcal M}^{n \times n}$ with ${\rm rank}({\mathbf B}) \leq k$. Here, ${\mathcal G}$ is the class of entry-wise positive matrices of finite rank $k$ and $\circ$ is the Hadamard matrix multiplication. By varying the class ${\mathcal G}$, the authors also introduce the class of $B_k^+$-stable matrices: an $n \times n$ real matrix $\mathbf A$ is called {\it $B_k^+$-stable} (to belong to the class $B_k^+$) if the Hadamard product ${\mathbf B}\circ{\mathbf A}$ is positive stable for every matrix ${\mathbf B} \in {\mathcal M}^{n \times n}$ such that $${\mathbf B} = {\mathbf B}_1 + {\mathbf B}_2 + \ldots + {\mathbf B}_k,$$ where each ${\mathbf B}_i$ is entry-wise positive, ${\rm rank}({\mathbf B}_i) = 1$, $i = 1, \ 2, \ \ldots, \ k$.

Now we consider ${\mathbb C} \setminus \{0\}$ as a stability region. Classes of nonsingular matrices that preserve nonsingularity under certain perturbations arise in many problems connected to stability. As an example of $({\mathfrak D}, {\mathcal G}, \circ)$ -stability with ${\mathfrak D} = {\mathbb C} \setminus \{0\}$, we consider the following matrix class, which "interpolates" nonsingular and sign-nonsingular matrices.

{\bf Hadamard nonsingularity (Johnson and van den Driessche, 1988).} This concept was introduced in \cite{JOHD} (see \cite{JOHD}, p. 368). An $n \times n$ matrix $\mathbf A$ is called {\it $B_k$-nonsingular} if the Hadamard product ${\mathbf B}\circ{\mathbf A}$ is nonsingular for every entry-wise positive matrix ${\mathbf B} \in {\mathcal M}^{n \times n}$ with ${\rm rank}({\mathbf B}) \leq k$ (in \cite{JOHD}, this matrix class is denoted by $L_{n,k}$). For strong forms of nonsingularity, see also \cite{DJD}.

\subsection{Additive $({\mathfrak D},{\mathcal G})$-stability}
Now let the binary operation $\circ$ be matrix addition. In spite of being naturally connected to multiplicative forms and a huge variety of applications, the concept of additive stability seems to have less attention in literature. Here and later on we consider the stability region ${\mathfrak D} = {\mathbb C}^-$.

{\bf Additive $D$-stability (Cross, 1978).} An $n \times n$ real matrix $\mathbf A$ is called {\it additive $D$-stable} if $-{\mathbf D}+{\mathbf A}$ is stable for every positive diagonal matrix $\mathbf D$. According to this definition, ${\mathcal G}$ is the class of {\bf negative} diagonal matrices (but if to consider positive stability, ${\mathcal G}$ will be changed to the class of positive diagonal matrices). This class was first defined in \cite{CROSS} (referring to the study of diffusion models of biological systems \cite{HAD}) under the name of {\it strong stability}.

{\bf Additive interval $D$-stability (Romanishin and Sinitskii, 2002)}
For studying additive $D$-stability, the following subset of ${\mathcal M}^{n \times n}$ is considered (see \cite{KOS}, \cite{RS1}):
$$\Theta_0 = \prod(0, \ d_{ii}^{max})= $$ $$ \diag (d_{ii}, \ 0 < d_{ii} < d_{ii}^{max} < +\infty, \ \ i = 1, \ 2, \ \ldots \ n). $$
An $n \times n$ matrix ${\mathbf A}$ is called {\it interval additive $D$-stable} with respect to $\Theta_0$ if $-{\mathbf D} +{\mathbf A}$ is stable for every ${\mathbf D} \in \Theta_0$.

{\bf Additive $H(\alpha)$-stability (Gumus and Xu, 2017).} The concept of additive $H(\alpha)$-stability was introduced in \cite{GUH1}: given a partition $\alpha = (\alpha_1, \ \ldots \ \alpha_p)$, an $n \times n$ real matrix $\mathbf A$ is called {\it additive $H(\alpha)$-stable} if $-{\mathbf H}+{\mathbf A}$ is stable for every symmetric positive definite $\alpha$-diagonal matrix $\mathbf H$. The matrix class ${\mathcal G}$ here is the class of symmetric negative definite $\alpha$-diagonal matrices.

The above three concepts were often considered together with the corresponding concepts of multiplicative $({\mathfrak D}, {\mathcal G})$-stability. However, the following case of the same nature, which is of great importance for the systems theory, is considered separately by quite different methods.

{\bf Finite rank perturbations (since 1960s).}  Given an $n \times n$ matrix $\mathbf A$, consider its finite-rank perturbation of the form
\begin{equation}\label{pert} \widetilde{{\mathbf A}} = {\mathbf A} + {\mathbf B}, \end{equation}
where ${\mathbf B} \in {\mathcal M}^{n \times n}$ with ${\rm rank}({\mathbf B}) \leq k$, $k = 1, \ \ldots, \ n$. If ${\rm rank}({\mathbf B}) = 1$, Equality \eqref{pert} may be written as
$$ \widetilde{{\mathbf A}} = {\mathbf A} + x \otimes y,$$
where $x, \ y \in {\mathbb R}^n$. The general problem is as follows: given a stability region $\mathfrak D$, and a class of vectors $V \subset {\mathbb R}^n$, when $\sigma(\widetilde{{\mathbf A}}) \subset {\mathfrak D}$ for all $x, \ y, \in V$?

The partial case of the above problem was studied by Barkovsky (see \cite{BARK}, \cite{BAY}): given an $n \times n$ matrix $\mathbf A$ and two vectors $x, \ y \in {\mathbb R}^n$, when all the matrices $\widetilde{{\mathbf A}}_{\tau} = {\mathbf A} +\tau (x \otimes y)$ have real spectra? This problem can be considered as establishing $({\mathfrak D},{\mathcal G},\circ)$-stability, where the stability region $\mathfrak D$ is the real axes, ${\mathcal G}$ is a parametric rank-one matrix family of the form $\{\tau (x\otimes y)\}_{\tau \in {\mathbb R}}$, the operation $\circ$ is matrix addition. For the problems of this type, see also \cite{FGR}, \cite{BIR}.
\subsection{Example of application to dynamical system stability} Now consider an example of applications of several types of $({\mathfrak D}, {\mathcal G}, \circ)$-stability to the stability of dynamical systems. Consider the system of second-order differential equations
\begin{equation}\label{syst2}
\ddot{x} = {\mathbf A}\dot{x} + {\mathbf B}x, \qquad x \in {\mathbb R^n}.
\end{equation}

The dynamics of System \eqref{syst2} is determined by a $2n \times 2n$ matrix of the form
\begin{equation}\label{C}{\mathbf C}:= \begin{pmatrix}{\mathbf A} & {\mathbf B} \\
 {\mathbf I} & {\mathbf O} \\ \end{pmatrix}, \end{equation}
 where ${\mathbf A}, \ {\mathbf B} \in {\mathcal M}^{n \times n}$ and ${\mathbf I}$ is $n \times n$ identity matrix.

It is well-known that System \eqref{syst2} is asymptotically stable if and only if the matrix $\mathbf C$ is stable, i.e. all its eigenvalues have negative real parts (see \cite{NIS}).
The following sufficient for stability conditions were established in \cite{NIS} (see \cite{NIS}, Theorems 2 and 3, Corollary 1. For the definition of {\it negative diagonally dominant (NDD)} matrices, see Appendix).
\begin{criterion}[\cite{NIS}]\label{Crit1}
Let ${\mathbf A}=\{a_{ij}\}_{i,j=1}^n$ and ${\mathbf B}=\{b_{ij}\}_{i,j = 1}^n$ be real $n \times n$ matrices, and the $2n \times 2n$ matrix ${\mathbf C}$ be defined by \eqref{C}. Let $a_{ii} < 0$ and $b_{ii} < 0$ for all $i = 1, \ \ldots, \ n$. If, in addition, $b_{ij} = 0$ for all $i \neq j$ (i.e. ${\mathbf B}$ be negative diagonal) and ${\mathbf A}$ is NDD, then ${\mathbf C}$ is stable.
\end{criterion}

\begin{criterion}[\cite{NIS}]\label{Crit2}
Let ${\mathbf A}=\{a_{ij}\}_{i,j=1}^n$ be a real $n \times n$ matrix with  $a_{ii} < 0$ for all $i = 1, \ \ldots, \ n$ and let the $2n \times 2n$ matrix ${\mathbf C}$ be defined as follows: $${\mathbf C}= \begin{pmatrix}{\mathbf A} & b{\mathbf I} \\
 {\mathbf I} & {\mathbf O} \\ \end{pmatrix}, $$ where $ b < 0$. Then ${\mathbf C}$ is stable if and only if $\mathbf A$ is stable.
\end{criterion}

\begin{criterion}[\cite{NIS}]\label{Crit3}
Let ${\mathbf B}=\{b_{ij}\}_{i,j = 1}^n$ be a real $n \times n$ matrix with $b_{ii} < 0$ for all $i = 1, \ \ldots, \ n$, and let the $2n \times 2n$ matrix ${\mathbf C}$ be defined by $${\mathbf C}= \begin{pmatrix} a{\mathbf I} & {\mathbf B} \\
 {\mathbf I} & {\mathbf O} \\ \end{pmatrix},$$ where $a < 0$. If all the eigenvalues of ${\mathbf B}$ are real and negative, then ${\mathbf C}$ is stable.
\end{criterion}

Basing on Criterion \ref{Crit1}, we establish the following criterion of (multiplicative) $D$-stability, which may be used for establishing stability of some perturbation of System \eqref{syst2}.

\begin{theorem}
 Let ${\mathbf A}=\{a_{ij}\}_{i,j=1}^n$ and ${\mathbf B}=\{b_{ij}\}_{i,j = 1}^n$ be real $n \times n$ matrices, and let the $2n \times 2n$ matrix ${\mathbf C}$ be defined by \eqref{C}. Let $a_{ii} < 0$ and $b_{ii} < 0$ for all $i$. If, in addition, $b_{ij} = 0$ for all $i \neq j$ and ${\mathbf A}$ is NDD, then ${\mathbf C}$ is multiplicative $D$-stable.
\end{theorem}

{\bf Proof.} Given a $2n \times 2n$ positive diagonal matrix $\mathbf D$, write it in the block-diagonal form:
$${\mathbf D} = {\rm diag}\{{\mathbf D}_{11}, \ {\mathbf D}_{22}\},$$
where ${\mathbf D}_{11}$, ${\mathbf D}_{22}$ are $n \times n$ positive diagonal matrices. Then
$${\mathbf D}{\mathbf C} = \begin{pmatrix} {\mathbf D}_{11} & {\mathbf O} \\
  {\mathbf O} & {\mathbf D}_{22} \\ \end{pmatrix}\begin{pmatrix} {\mathbf A} & {\mathbf B} \\
 {\mathbf I} & {\mathbf O} \\ \end{pmatrix} = \begin{pmatrix} {\mathbf D}_{11}{\mathbf A} & {\mathbf D}_{11}{\mathbf B} \\
 {\mathbf D}_{22} & {\mathbf O} \\ \end{pmatrix}.$$
To study the spectrum of ${\mathbf D}{\mathbf C}$, consider the similarity transformation $\widetilde{{\mathbf C}}:=\widetilde{{\mathbf D}}^{-1}({\mathbf D}{\mathbf C})\widetilde{{\mathbf D}}$, where $\widetilde{{\mathbf D}} = {\rm diag}\{{\mathbf I}, \ {\mathbf D}^{-1}_{22}\}$.
Then $$\widetilde{{\mathbf C}} = \begin{pmatrix} {\mathbf I} & {\mathbf O} \\
  {\mathbf O} & {\mathbf D}^{-1}_{22} \\ \end{pmatrix}\begin{pmatrix} {\mathbf D}_{11}{\mathbf A} & {\mathbf D}_{11}{\mathbf B} \\
 {\mathbf D}_{22} & {\mathbf O} \\ \end{pmatrix}\begin{pmatrix} {\mathbf I} & {\mathbf O} \\
  {\mathbf O} & {\mathbf D}_{22} \\ \end{pmatrix} = \begin{pmatrix} {\mathbf D}_{11}{\mathbf A} & {\mathbf D}_{11}{\mathbf B}{\mathbf D}_{22} \\
 {\mathbf I} & {\mathbf O} \\ \end{pmatrix}.$$
 Clearly, $\sigma(\widetilde{{\mathbf C}}) = \sigma({\mathbf D}{\mathbf C})$ and $\widetilde{{\mathbf C}}$ is of Form \eqref{C}. Moreover, ${\mathbf D}_{11}{\mathbf A}$ is NDD for all positive diagonal ${\mathbf D}_{11}$ and ${\mathbf D}_{11}{\mathbf B}{\mathbf D}_{22}$ is negative diagonal for all positive diagonal ${\mathbf D}_{11}$ and ${\mathbf D}_{22}$. Thus applying Criterion \ref{Crit1}, we obtain the stability of ${\mathbf D}{\mathbf C}$.
 $\square$

Let us consider the perturbations of System \eqref{syst2} of the following form:
\begin{equation}\label{syst2per1}
\ddot{x} = {\mathbf D}{\mathbf A}\dot{x} + {\mathbf B}x, \qquad x \in {\mathbb R^n},
\end{equation}
where ${\mathbf D}$ is $n \times n$ positive diagonal matrix.

This type of system perturbations corresponds to the following perturbation of the matrix ${\mathbf C}$:
\begin{equation}\label{Cper}\widetilde{{\mathbf C}}:= \begin{pmatrix}{\mathbf D}{\mathbf A} & {\mathbf B} \\
 {\mathbf I} & {\mathbf O} \\ \end{pmatrix}. \end{equation}

Matrix perturbation \eqref{Cper} can be described with the help of {\it block Hadamard products} (for the definition and studies see \cite{CHO}, \cite{HMN}).

Given ${\mathbf H}, \ {\mathbf G} \in {\mathcal M}^{2n \times 2n}$, partitioned into $n \times n$ blocks as follows:
$${\mathbf H} =  \begin{pmatrix}{\mathbf H}_{11} & {\mathbf H}_{12} \\
{\mathbf H}_{21} & {\mathbf H}_{22} \\ \end{pmatrix}, \qquad {\mathbf G} =  \begin{pmatrix}{\mathbf G}_{11} & {\mathbf G}_{12} \\
{\mathbf G}_{21} & {\mathbf G}_{22} \\ \end{pmatrix}.$$
Then their block Hadamard product is defined by
$${\mathbf H}\diamond {\mathbf G}:= \begin{pmatrix}{\mathbf H}_{11}{\mathbf G}_{11} & {\mathbf H}_{12}{\mathbf G}_{12} \\
{\mathbf H}_{21}{\mathbf G}_{21} & {\mathbf H}_{22}{\mathbf G}_{22} \\ \end{pmatrix},$$
where ${\mathbf H}_{ij}{\mathbf G}_{ij}$, $i,j = 1, \ 2$, denotes the "usual" matrix product of $n \times n$ blocks ${\mathbf H}_{ij}$ and ${\mathbf G}_{ij}$.

Now let us define a matrix class ${\mathcal G}_1$ by
$${\mathcal G}_1:= \{ {\mathbf G} \in {\mathcal M}^{2n \times 2n}: {\mathbf G} = \begin{pmatrix}{\mathbf D} & {\mathbf I} \\
{\mathbf I} & {\mathbf I} \\ \end{pmatrix} \},$$
where ${\mathbf D}$ is $n \times n$ positive diagonal matrix, ${\mathbf I}$ is $n \times n$ identity. Then, for a block matrix ${\mathbf C}$ of Form \eqref{C} and an arbitrary matrix ${\mathbf G} \in {\mathcal G}_1$, we obtain:

$${\mathbf G}\diamond{\mathbf C}= \begin{pmatrix}{\mathbf D}{\mathbf A} & {\mathbf B} \\
 {\mathbf I} & {\mathbf O} \\ \end{pmatrix}  = \widetilde{{\mathbf C}}.$$

Using Criterion \ref{Crit2}, we obtain the following result.
\begin{theorem} Let ${\mathbf A}=\{a_{ij}\}_{i,j=1}^n$ with $a_{ii} < 0$ for all $i = 1, \ \ldots, \ n$
and let the $2n \times 2n$ matrix ${\mathbf C}$ be defined as follows: $${\mathbf C}= \begin{pmatrix}{\mathbf A} & b{\mathbf I} \\
 {\mathbf I} & {\mathbf O} \\ \end{pmatrix}, $$ where $ b < 0$.
Then ${\mathbf C}$ is $({\mathcal G}_1, \diamond)$-stable if and only if $\mathbf A$ is $D$-stable.
\end{theorem}

Now consider the perturbations of System \eqref{syst2} of the following form:
\begin{equation}\label{syst2per1}
\ddot{x} = {\mathbf A}\dot{x} + {\mathbf D}{\mathbf B}x, \qquad x \in {\mathbb R^n},
\end{equation}
where ${\mathbf D}$ is $n \times n$ positive diagonal matrix.

For this, we define a matrix class ${\mathcal G}_2$ by
$${\mathcal G}_2:= \{ {\mathbf G} \in {\mathcal M}^{2n \times 2n}: {\mathbf G} = \begin{pmatrix}{\mathbf I} & {\mathbf D} \\
{\mathbf I} & {\mathbf I} \\ \end{pmatrix} \},$$
where ${\mathbf D}$ is $n \times n$ positive diagonal matrix, ${\mathbf I}$ is $n \times n$ identity. Then, for a block matrix ${\mathbf C}$ of Form \eqref{C} and an arbitrary matrix ${\mathbf G} \in {\mathcal G}_2$, we obtain:

$${\mathbf G}\diamond{\mathbf C}= \begin{pmatrix}{\mathbf A} & {\mathbf D}{\mathbf B} \\
 {\mathbf I} & {\mathbf O} \\ \end{pmatrix}.$$

 Criterion \ref{Crit3} immediately implies:

\begin{theorem} Let ${\mathbf B}=\{b_{ij}\}_{i,j = 1}^n$ be a real $n \times n$ matrix with $b_{ii} < 0$ for all $i = 1, \ \ldots, \ n$, and let the $2n \times 2n$ matrix ${\mathbf C}$ be defined by $${\mathbf C}= \begin{pmatrix} a{\mathbf I} & {\mathbf B} \\
 {\mathbf I} & {\mathbf O} \\ \end{pmatrix},$$ where $a < 0$. If the matrix ${\mathbf B}$ is $D$-negative, then ${\mathbf C}$ is $({\mathcal G}_2, \diamond)$-stable.
\end{theorem}

Here, we may impose on $\mathbf B$ any condition which guarantee $D$-negativity, for example, take $-{\mathbf B}$ strictly totally positive.

\section{The historical development of $\mathfrak D$-stability (with a view to robust problems)}

In the last decades, matrix and polynomial $\mathfrak D$-stability, also known as matrix (respectively, polynomial) {\it root clustering} has become an attractive area for researchers. Giving a brief overview, the theory has been developed from the simplest and the most used partial cases to more and more sophisticated and general stability regions (see, for example, \cite{GUT2}, \cite{JU2}, \cite{JU3}, \cite{JUA}). However, the regions described by polynomial conditions are not easy to study. Thus, different kinds of linearizations and other representations of complicated regions in a simpler form, are introduced (\cite{CHIGA}, \cite{CGA}, \cite{PABB}, \cite{BAM}, \cite{BHPM}). Here, we separate the following two approaches:
\begin{enumerate}
\item[\rm -] polynomial approach, by the transition from the system matrix $\mathbf A$ to its characteristic polynomial $f({\mathbf A})$ and then applying to $f({\mathbf A})$ the results on polynomial root clustering;
\item[\rm -] matrix approach, by considering some quadratic forms, constructed from $\mathbf A$ (e.g. Lyapunov theorem and different its modifications) and by the localization of the eigenvalues in some regions of the complex plane defined by the matrix entries (e.g. Gershgorin theorem).
\end{enumerate}
 We give a brief historical overview of both of the approaches, focusing on the crucial (for the development of "super"-stability theory) results. As it will be shown later for the most studied partial case ${\mathfrak D} = {\mathbb C}^-$, different stability criteria lead to different approaches to the study of multiplicative and additive $D$-stability and imply sufficient for $D$-stability conditions of different nature. For studying $({\mathfrak D}, {\mathcal G}, \circ)$-stability, it would be natural to consider the generalizations of stability criteria to the case of other stability regions $\mathfrak D$. For the convenience of further analysis, we collect here the results connected to stability and ${\mathfrak D}$-stability.

 Given ${\mathbf A}, \ {\mathbf B} \in {\mathcal M}^{n \times n}$, we use the notation ${\mathbf A} \prec {\mathbf B}$ $({\mathbf A} \succ {\mathbf B})$ if the matrix ${\mathbf A} - {\mathbf B}$ is negative definite (respectively, positive definite).

\subsection{Polynomial approach} Given a polynomial $p(z)$, we call it {\it ${\mathfrak D}$-stable} if $p(z) = 0$ implies $z \in {\mathfrak D}$ for any $z \in {\mathbb C}$. In the case when ${\mathfrak D} = {\mathbb C}^-$, ${\mathfrak D}$-stable polynomials are called just {\it stable}. Let us consider a perturbation $\widetilde{{\mathbf A}}$ of a matrix $\mathbf A$. Obviously, the transition from the study of a perturbed matrix $\widetilde{{\mathbf A}}$ to some perturbation of the characteristic polynomial $f_{\mathbf A}$ of the initial matrix $\mathbf A$ and back, cause certain difficulties. Thus we do not mention here a lot of results on polynomial root clustering, but only those, which allow us to make this transition easily. These results are applicable to the analysis of $({\mathfrak D}, {\mathcal G}, \circ)$-stability, where we come from the study of the initial characteristic polynomial $f_{\mathbf A}$ to the study of a perturbed polynomial family $f_{{\mathbf G}\circ{\mathbf A}}$, where ${\mathbf G}$ varies along the class $\mathcal G$.

{\bf Classical stability criteria.} For the classical examples of the stability regions ${\mathfrak D}$, this concept goes back to Descartes \cite{DEC} and is developed in the papers by Cauchy \cite{CAU}, Sturm \cite{STU}, Hermite \cite{HERMI}, Routh \cite{ROU1}-\cite{ROU4}, Hurwitz \cite{HUR} (${\mathfrak D}$ is the left-hand side of the complex plane) and Schur \cite{SHU2}, Cohn \cite{COH} (${\mathfrak D}$ is the open unit disk). Studying polynomials all whose zeroes are real or positive (for the beginning, see \cite{WAR}) can be also considered as a partial case of $\mathfrak D$-stability (${\mathfrak D}$ is the real line or its positive direction).

{\bf Kharitonov stability criterion for interval polynomials.} One of the most prominent results in robust stability of polynomials is Kharitonov stability criterion obtained in 1978. In the case, when the coefficients of a polynomial are not exactly defined, we consider the following family of polynomials:
    \begin{equation}\label{pf}F^n(z) : = \{f(z): f(z) = \sum_{i=0}^n a_ix^i ; \ a_i^- \leq a_i \leq a_i^+; \ a_n \leq 0 \}.\end{equation}
The family \eqref{pf} is called an {\it interval polynomial} and denoted
$$ F^n(z) = \sum_{i=0}^n[a_i^-,a_i^+]z^{n-i}, \qquad a_0^- \neq 0.$$
Kharitonov stability criterion surprisingly gives a necessary and sufficient condition for the stability of an infinite number of polynomials by testing only four polynomials of a special form (the so-called Kharitonov polynomials) (see \cite{KHAR1}, p. 2086-2087, Theorems 1 and 2):
$$k_1(z) = z^n + a^1_1z^{n-1} + \ldots + a^1_n, $$
where $$ a^1_{n-2k} =\left\{\begin{array}{cc} a^+_{n-2k} & \mbox{if $k$ is even}; \\
a^-_{n-2k} & \mbox{if $k$ is odd} \end{array}\right. \qquad a^1_{n-2k-1} =\left\{\begin{array}{cc} a^+_{n-2k-1} & \mbox{if $k$ is even}; \\
a^-_{n-2k-1} & \mbox{if $k$ is odd} \end{array}\right.$$

$$k_2(z) = z^n + a^2_1z^{n-1} + \ldots + a^2_n, $$
where $$ a^2_{n-2k} =\left\{\begin{array}{cc} a^-_{n-2k} & \mbox{if $k$ is even}; \\
a^+_{n-2k} & \mbox{if $k$ is odd} \end{array}\right. \qquad a^2_{n-2k-1} =\left\{\begin{array}{cc} a^-_{n-2k-1} & \mbox{if $k$ is even}; \\
a^+_{n-2k-1} & \mbox{if $k$ is odd} \end{array}\right.$$

$$k_3(z) = z^n + a^3_1z^{n-1} + \ldots + a^3_n, $$
where $$ a^3_{n-2k} =\left\{\begin{array}{cc} a^-_{n-2k} & \mbox{if $k$ is even}; \\
a^+_{n-2k} & \mbox{if $k$ is odd} \end{array}\right. \qquad a^3_{n-2k-1} =\left\{\begin{array}{cc} a^+_{n-2k-1} & \mbox{if $k$ is even}; \\
a^-_{n-2k-1} & \mbox{if $k$ is odd} \end{array}\right.$$

$$k_4(z) = z^n + a^4_1z^{n-1} + \ldots + a^4_n, $$
where $$ a^4_{n-2k} =\left\{\begin{array}{cc} a^+_{n-2k} & \mbox{if $k$ is even}; \\
a^-_{n-2k} & \mbox{if $k$ is odd} \end{array}\right. \qquad a^4_{n-2k-1} =\left\{\begin{array}{cc} a^-_{n-2k-1} & \mbox{if $k$ is even}; \\
a^+_{n-2k-1} & \mbox{if $k$ is odd} \end{array}\right.$$
\begin{theorem}[Kharitonov] \label{KH} An interval polynomial (i.e. all the members of the family \eqref{pf}) is stable if and only if the four Kharitonov polynomials $k_1(z)$, $k_2(z)$, $k_3(z)$, $k_4(z)$ are stable.
\end{theorem}

The following ways of generalizing the Kharitonov theorem are considered in literature: with respect to different structures of polynomial uncertainties  and with respect to different kinds of stability regions (the so-called Kharitonov regions). Both of these ways may be used for studying different types of $({\mathfrak D}, {\mathcal G}, \circ)$-stability.

\subsection{Lyapunov theorem approach} Here, we place the most important results connected to Lyapunov theorem and its generalizations, in chronological order. As the Lyapunov theorem and Lyapunov equation analysis play crucial role in the study of multiplicative and additive $D$-stability, the generalizations of the Lyapunov theorem provide a natural tool for studying $D$-hyperbolicity, Schur $D$-stability and other concepts. In Section 8, we discuss the relations between solvability of generalized Lyapunov equations for different regions $\mathfrak D$ and $({\mathfrak D}, {\mathcal G}, \circ)$-stability in more details.

{\bf 1892 --- Lyapunov.}
Remind that an $n \times n$ real matrix $\mathbf A$ is called {\it Hurwitz stable} or just {\it stable} if all its eigenvalues have negative real parts. The approach we analyze in this subsection is based on the necessary and sufficient condition of matrix stability, proved by Lyapunov (see, for example, \cite{BELL}, \cite{GANT2}, for exact formulation see \cite{GANT}, \cite{HER1}, p. 164, Theorem 2.4).
    \begin{theorem}[Lyapunov] An $n \times n$ matrix $\mathbf A$ is stable if and only if there exists a symmetric positive definite matrix $\mathbf H$ such that the matrix
    $${\mathbf W}:={\mathbf H}{\mathbf A} + {\mathbf A}^{T}{\mathbf H}$$
    is negative definite.
    \end{theorem}
Equivalently, we analyze the solvability of the Lyapunov equation
\begin{equation}\label{lyap}
{\mathbf H}{\mathbf A} + {\mathbf A}^{T}{\mathbf H} = {\mathbf W},
\end{equation}
where ${\mathbf W}$ is a symmetric negative definite matrix, in the class of symmetric positive definite matrices.

The matrix ${\mathbf A}$ is stable if and only if, for any given negative definite matrix ${\mathbf W}$, the Lyapunov equation \eqref{lyap} has a unique symmetric solution ${\mathbf H}$, and this solution ${\mathbf H}$ is positive definite (see \cite{CHEN2}, p. 132). Partial cases, when $\mathbf H$ belongs to a specified subclass of positive definite matrices, are of great interest. The following concept has enormous number of applications: a matrix ${\mathbf A} \in {\mathcal M}^{n \times n}$ is called {\it diagonally stable} if the Lyapunov equation \eqref{lyap} has a positive diagonal solution ${\mathbf D}$, in other words, if the matrix $${\mathbf W}:={\mathbf D}{\mathbf A} + {\mathbf A}^{T}{\mathbf D}$$
is negative definite for some positive diagonal matrix ${\mathbf D}$. In this case, $\mathbf D$ is called a {\it Lyapunov scaling factor}.

{\bf 1952 --- Stein.} Here, we mention an analogous statement for Schur stability (see \cite{STE}, also \cite{TAU1}, \cite{TAU2}, \cite{WIM}). Remind, that an $n \times n$ real matrix $\mathbf A$ is called {\it Schur stable} if all its eigenvalues lie inside the unit circle, i.e. the spectral radius $\rho(\mathbf A) < 1$.

    \begin{theorem}[Stein] An $n \times n$ matrix $\mathbf A$ is Schur stable if and only if there exists a symmetric positive definite matrix $\mathbf H$ such that the matrix
    \begin{equation}\label{ste}{\mathbf W}:= {\mathbf A}^{T}{\mathbf H}{\mathbf A} - {\mathbf H}\end{equation}
    is negative definite.
    \end{theorem}

{\bf 1962 --- Ostrowski and Schneider.}
    Here, instead of Lyapunov theorem, we deal with the following theorem proved in \cite[p. 76]{OSS} (see \cite{OSS}, p. 76, Theorem 1).
    \begin{theorem}[Ostrowski, Schneider]\label{OS} An $n \times n$ matrix $\mathbf A$ has no pure imaginary eigenvalues (i.e. with zero real parts) if and only if there exists a symmetric matrix $\mathbf H$ such that the matrix
    $${\mathbf W}:={\mathbf H}{\mathbf A} + {\mathbf A}^{T}{\mathbf H}$$
    is positive definite. Then we have ${\rm In}({\mathbf H}) = {\rm In}({\mathbf A})$.
    \end{theorem}

{\bf 1969 --- Hill.} For a class of more general stability regions $\mathfrak D$, we need the following generalization of Lyapunov theorem, obtained by Hill (see \cite{HIL}, also \cite{WIM}, p. 140, Theorem 1).

\begin{theorem}[Generalized Lyapunov]\label{TGenLyap} If an $n \times n$ matrix $\mathbf A$ satisfies the matrix equation
\begin{equation}\label{GenLyap} \sum_{i,j = 0}^{n-1}c_{ij}({\mathbf A}^T)^i{\mathbf H}{\mathbf A}^j = {\mathbf W}, \qquad c_{ij} = c_{ji}
\end{equation}
with a symmetric positive definite matrix $\mathbf H$, then
\begin{enumerate}
\item[\rm 1.] ${\mathbf W}$ is a symmetric positive definite matrix implies $$f(\lambda) := \sum_{i,j = 0}^{n-1}c_{ij}\overline{\lambda}^i\lambda^j > 0;$$
\item[\rm 2.] ${\mathbf W}$ is a symmetric positive semidefinite matrix implies $$f(\lambda) := \sum_{i,j = 0}^{n-1}c_{ij}\overline{\lambda}^i\lambda^j \geq 0;$$
\item[\rm 3.] ${\mathbf W} = {\mathbf 0}$ implies $$f(\lambda) := \sum_{i,j = 0}^{n-1}c_{ij}\overline{\lambda}^i\lambda^j = 0.$$
\end{enumerate}
\end{theorem}

{\bf 1981 -- Gutman and Jury.} The theory of stability in some generalized regions known as {\it root clustering} was mainly developed in 1980s by Gutman and Jury (see \cite{GUJU}) and continued in \cite{GUT}, \cite{GUT2}, where a lot of special classes of stability regions were analyzed. A generalization of Lyapunov theorem was introduced for a special type of regions called GLE (Generalized Lyapunov Equation) regions. As the examples of GLE regions study, we refer to \cite{HHP}, \cite{LLK}, \cite{MEME}, \cite{MAO} (disk regions), \cite{ANB}, \cite{BICJ1}, \cite{ANB}, \cite{BICJ2}, \cite{BICJ3}, \cite{MOMOK}, \cite{SGH} (sector regions).

{\bf 1996 --- Chilali and Gahinet.} Since it is hard to analyze GLE regions, defined by matrix inequalities of Form \eqref{GenLyap}, due to the polynomial nature of the conditions, their linearization catched enormous attention. A subset ${\mathfrak D} \subset {\mathbb C}$ that can be defined as
\begin{equation}\label{LMI} {\mathfrak D} = \{z \in {\mathbb C}: \ {\mathbf L} + {\mathbf M}z+{\mathbf M}^T\overline{z} \prec 0\},\end{equation}
where ${\mathbf L}, {\mathbf M} \in {\mathcal M}^{n \times n}$, ${\mathbf L}^T = {\mathbf L}$, is called an {\it LMI region} with the {\it characteristic function} $f_{\mathfrak D} = {\mathbf L} + z{\mathbf M}+\overline{z}{\mathbf M}^T$ (see \cite{CHIGA}, \cite{CGA}). LMI regions are dense in the set of convex regions that are symmetric with respect to the real axis. Thus they include a lot of regions of great importance (including the left-hand side of the complex plane and the unit disk) and a huge variety of other stability regions can be approximated by LMI regions.
Let us recall the following result (see \cite{CHIGA}, p. 360, Theorem 2.2).
\begin{theorem}[Lyapunov theorem for LMI regions]
Given an LMI region $\mathfrak D$, defined by \eqref{LMI}, a  matrix $\mathbf A$ is $\mathfrak D$-stable if and only if there is a symmetric positive definite matrix $\mathbf H$ such that the matrix
    \begin{equation}\label{LMIeq}{\mathbf W} :=\ {\mathbf L}\otimes {\mathbf H} + {\mathbf M}\otimes({\mathbf H}{\mathbf A})+{\mathbf M}^T\otimes({\mathbf A}^T{\mathbf H})\end{equation}
    is negative definite.
\end{theorem}
Note, that the dimension of matrices $\mathbf A$ and $\mathbf H$ is supposed to be $n$, while the dimension of matrices $\mathbf L$ and $\mathbf M$ is supposed to be $m$ which does not depend on $n$.

{\bf 2000 --- Peaucelle, Arzelier, Bachelier, Bernussou.} Due to their convexity, LMI regions does not include some important regions (e.g. the exterior of the unit disk), that are nonconvex. Thus the following kind of possibly nonconvex regions was introduced in \cite{PABB}. Let ${\mathbf R} \in {\mathcal M}^{2d \times 2d}$ be a symmetric matrix partitioned as
$${\mathbf R} = \begin{pmatrix}{\mathbf R}_{11} & {\mathbf R}_{12} \\ {\mathbf R}_{12}^T & {\mathbf R}_{22} \end{pmatrix},$$
where ${\mathbf R}_{11} = {\mathbf R}_{11}^T \in {\mathcal M}^{d \times d}$ and ${\mathbf R}_{11} = {\mathbf R}_{11}^T \in {\mathcal M}^{d \times d}$. A subset ${\mathfrak D} \subset {\mathbb C}$ that can be defined as
\begin{equation}\label{EMI} {\mathfrak D} = \{z \in {\mathbb C}: \ {\mathbf R}_{11} + {\mathbf R}_{12}z+{\mathbf R}_{12}^T\overline{z} + {\mathbf R}_{22}z\overline{z} \prec 0\},\end{equation}
 is called an {\it EMI region} (Ellipsoidal Matrix Inequality). If we do not impose any restrictions (e.g. positive definiteness) on ${\mathbf R}_{22}$, then EMI regions are not necessarily convex. The Lyapunov characterization of EMI regions can be easily deduced from generalized Lyapunov Theorem.
\begin{theorem}[Lyapunov theorem for EMI regions]
Given an EMI region $\mathfrak D$, defined by \eqref{EMI}, a  matrix $\mathbf A$ is $\mathfrak D$-stable if and only if there is a symmetric positive definite matrix $\mathbf H$ such that the matrix
    \begin{equation}\label{EMIeq} {\mathbf W} :=\ {\mathbf R}_{11}\otimes {\mathbf H} + {\mathbf R}_{12}\otimes({\mathbf H}{\mathbf A})+{\mathbf R}_{12}^T\otimes({\mathbf A}^T{\mathbf H})+{\mathbf R}_{22}\otimes({\mathbf A}^T{\mathbf H}{\mathbf A})\end{equation}
    is negative definite.
\end{theorem}

All the above versions of the Lyapunov theorem will be used later for introducing generalizations of the concept of diagonal stability.

\subsection{Other necessary and sufficient stability criteria} Here, we collect results that provide an alternative to the well-known Routh--Hurwitz conditions. These results will be used later for studying multiplicative (Duan and Patton) and additive (Li and Wang) $D$-stability. Generalizations of these results for other regions $\mathfrak D$ would be also of interest.

{\bf 1998 --- Li and Wang.} The following stability criterion is provided in \cite{LIW} (for the definition of compound and additive compound matrices, see  Appendix).
\begin{theorem} Let ${\mathbf A} \in {\mathcal M}^{n \times n}$, consider its second additive compound matrix ${\mathbf A}^{[2]}$. For ${\mathbf A}$ to be Hurwitz stable, it is necessary and sufficient that ${\mathbf A}^{[2]}$ is Hurwitz stable and $(-1)^n\det({\mathbf A}) > 0$.
\end{theorem}

{\bf 1998 --- Duan and Patton.} The necessary and sufficient characterization of stable matrices by matrix factorization was deduced from Lyapunov equation (see \cite{DUP}).
\begin{theorem}\label{DU} A matrix ${\mathbf A} \in {\mathcal M}^{n \times n}$ is stable if and only if there is a (not necessarily symmetric) negative definite matrix $\mathbf G$ and a symmetric positive definite matrix $\mathbf L$ such that
$${\mathbf A} = {\mathbf G}{\mathbf L}. $$
\end{theorem}

\subsection{Stability of certain matrix classes}
Here, we collect stability criteria which are suitable for special forms of matrix. These criteria will be used later in connection with multiplicative and additive $D$-stability. It would be of interest to study their generalizations in connection with a more general question: "When matrix $\mathfrak D$-stability would imply $({\mathfrak D}, \mathcal G, \circ)$-stability?"

{\bf 1931 --- Gershgorin.} The following prominent result describes easily computable domain that contains all the eigenvalues of a matrix (see, for example, \cite{HOJ}, p. 344).
\begin{theorem}[Gershgorin]\label{GER} Let ${\mathbf A} = \{a_{i,j}\}_{i,j = 1}^n \in {\mathcal M}^{n \times n}$, define
$$R_i := \sum_{j =1; \ i \neq j}^n |a_{ij}|, \qquad 1 \leq i \leq n.$$ Let $D(a_{ii},R_i) \subset {\mathbb C}$ be a closed disk centered at $a_{ii}$ with the radius $R_i$. Then all the eigenvalues of ${\mathbf A}$ are located in the union of $n$ discs $$G({\mathbf A}):= \bigcup_{i=1}^n D(a_{ii},R_i). $$
\end{theorem}
\begin{corollary} Strictly diagonally dominant matrices with negative principal diagonal entries are stable.
\end{corollary}
The proof of the corollary just shows that the union of Gershgorin disks $G({\mathbf A})$ of a strictly diagonally dominant matrix $\mathbf A$ with negative principal diagonal is located in the left-hand side of the complex plane.

{\bf 1978 --- Tyson, Othmer.} The following stability result was first considered in \cite{TYO} (see also \cite{THR1} (Appendix A) for the exact proof and \cite{THR2} for the further study) with the application to sequences of biochemical reactions.
\begin{theorem}[Secant criterion]\label{Sect} Let ${\mathbf A} \in {\mathcal M}^{n \times n}$ be of the form
\begin{equation}\label{STR2}{\mathbf A} = \begin{pmatrix} -\alpha_1 & 0 & \ldots & 0 & - \beta_n \\
\beta_1 & -\alpha_2 & \ddots & \ddots & 0 \\
0 & \beta_2 & -\alpha_3 & \ddots & \vdots \\
\vdots & \ddots & \ddots & \ddots & 0 \\
0 & \ldots & 0 & \beta_{n-1} & - \alpha_n \\ \end{pmatrix},\end{equation} where $\alpha_i > 0, \ \beta_i > 0, \ i = 1, \ \ldots, \ n.$
Then $\mathbf A$ is Hurwitz stable if
$$ \frac{\beta_1 \ldots \beta_n}{\alpha_1 \ldots \alpha_n} < \sec(\frac{\pi}{n})^n.$$
\end{theorem}

\section{The historical development of diagonal and $D$-stability}
 We provide a detailed survey of the most well-studied partial cases of $({\mathfrak D}, {\mathcal G}, \circ)$-stability, also collecting the methods of their study. Together with multiplicative and additive $D$-stability, we consider the concept of diagonal stability, i.e. the existence of a positive diagonal solution of the Lyapunov equation. Also note that in matrix literature, due to the convenience reasons, positive stability is often referred: a matrix ${\mathbf A} \in {\mathcal M}^{n \times n}$ is called {\it positive stable} if all its eigenvalues have positive real parts (see \cite{HER1}, \cite{JOHN1}). The concepts of diagonal stability, multiplicative and additive $D$-stability all over this section may be based on stability (${\mathfrak D} = {\mathbb C}^-$) as well as on positive stability (${\mathfrak D} = {\mathbb C}^+$). The reader can easily understand from the content, what kind of stability is supposed.
  
  Sketches of the proofs in this section are given for the convenience of future generalizations to the other cases of $({\mathfrak D}, {\mathcal G}, \circ)$-stability.

\subsection{1950-1960s. Basic definitions and statements} In this period, the general problem of characterizing multiplicative $D$-stable matrices was raised. Elementary properties of $D$-stable matrices were studied and necessary for $D$-stability conditions were analyzed. Some important matrix classes, such as $M$-matrices and negative definite matrices are shown to be $D$-stable. The following two general methods were developed for studying $D$-stability.
\begin{enumerate}
\item[\rm 1.] {\bf Lyapunov equation analysis.} The conditions sufficient for $D$-stability may be derived from the solvability of Lyapunov equation by putting the right-hand side ${\mathbf W}$ to be of a special form and deriving conditions for $\mathbf A$, or by imposing additional properties on the solution ${\mathbf P}$ (e.g. to be positive diagonal).
\item[\rm 2.] {\bf Qualitative approach.} This approach uses the study of the spectral properties of matrix families, which members have entries of prescribed signs (positive, negative or zero). Such a structure is obviously preserved under multiplication by a positive diagonal matrix.
\end{enumerate}
The details are given below.

{\bf 1956 --- Enthoven and Arrow.} The problem of $D$-stability was raised when studying the "expected price" model. $D$-stability conditions for a Metzler matrix were established (see \cite{ENT}).

{\bf 1956-1958 --- Arrow and McManus.} In \cite{AM1}, the study of Metzler matrices was continued. In \cite{AM}, the general problem of $D$-stability characterization was raised: "If $\mathbf A$ is stable, in what circumstances is ${\mathbf D}{\mathbf A}$ stable, where $\mathbf D$ is diagonal?" The following sufficient condition was proved.
\begin{theorem}\cite{AM}
Negative definite matrices are $D$-stable.
\end{theorem}

The elementary properties of $D$-stable matrices and transformations which preserves $D$-stability were studied. Here, we mention the following statement.
\begin{theorem}\cite{AM} Let $\mathbf S$ be a nonsingular diagonal matrix, $\mathbf A$ is $D$-stable if and only if ${\mathbf S}{\mathbf A}{\mathbf S}^{-1}$ is $D$-stable.
\end{theorem}

{\bf 1958 --- Fisher and Fuller.} The following powerful result was first proved by Fisher and Fuller (see \cite{FIF}), with a number of simpler proofs appeared later (see \cite{FISH}, \cite{BAL}). This result was used by many authors for establishing conditions sufficient for stability. The proofs, based on the analysis of the characteristic polynomial, are on independent interest as examples of applying polynomial criteria to obtain results on matrix eigenvalue localization.

 \begin{theorem}[Fisher, Fuller]\cite{FIF}\label{FIF} Let $\mathbf A$ be an $n \times n$ real matrix all of whose leading principal minors are positive. Then there is an $n \times n$ positive diagonal matrix $\mathbf D$ such that all the roots of ${\mathbf D}{\mathbf A}$ are positive and simple.
\end{theorem}

{\bf 1961 --- Taussky.} Recall the following fact: if a matrix $\mathbf A$ is stable, the Lyapunov equation \eqref{lyap} is solvable for any negative definite right-hand side ${\mathbf W}$. In \cite{TAU1}, certain sufficient for stability conditions were derived from the Lyapunov equation, putting ${\mathbf W}:=-{\mathbf I}$. These conditions lead to the following criterion of $D$-stability.
\begin{theorem}\cite{TAU1} Let ${\mathbf A} = {\mathbf B}-\alpha{\mathbf I}$, where $\alpha>0$ and $\mathbf B$ be a skew-symmetric matrix. Then ${\mathbf D}{\mathbf A}$ is stable for every positive diagonal matrix $\mathbf D$.
\end{theorem}

{\bf 1965 --- Quirk and Ruppert.} The {\it qualitative stability} approach was developed in \cite{QR}. The core ideas are as follows. Denote ${\mathcal A}$ the set of all matrices, sign-similar to a given matrix ${\mathbf A}$ (for the definition of sign-similarity, see Appendix). Then ${\mathbf A}$ is called {\it sign-stable} or {\it qualitative stable}, if any matrix from ${\mathcal A}$ is stable. Note, that we may consider ${\mathcal A}$ as an interval matrix, with the entries belong to one of the sets $(0, \ +\infty)$, $(- \infty, \ 0)$ or $\{0\}$. The following inclusion was established.
 \begin{theorem}\cite{QR} Sign-stable matrices are $D$-stable.
 \end{theorem}

It was stated referring \cite{AM} that diagonal stability is a sufficient condition for $D$-stability.
\begin{theorem}\cite{QR} If there exists a positive diagonal matrix $\mathbf D$ such that ${\mathbf W}: = {\mathbf D}{\mathbf A} + {\mathbf A}^T{\mathbf D}$ is negative definite, then $\mathbf A$ is $D$-stable.
\end{theorem}

The following important necessary condition was proved.
\begin{theorem}\cite{QR} If $\mathbf A$ is $D$-stable then $\mathbf A$ is almost Hicksian.
\end{theorem}

The concept of {\it total stability} is introduced: a matrix ${\mathbf A}$ is called {\it totally stable} if every principal submatrix of $\mathbf A$ is $D$-stable (note, that this matrix property is sometimes referred as {\it total $D$-stability} (see, for example, \cite{KAB})). Such matrices, which are known to be Hicksian, are connected to the results of Metzler \cite{ME}.

\subsection{1970s. Crucial results} In this decade, the most important results which formed the basis of the later study were obtained. Certain new classes of multiplicative $D$-stable matrices were described. The analysis of known classes of $D$-stable matrices, such as $M$-matrices, was continued, some generalizations of known classes appeared. Biological applications gave independent interest to diagonal stability, including different characterizations of diagonally stable matrices. Implications between diagonal stability, multiplicative and additive $D$-stability were established and intensively studied.

The following new methods for studying $D$-stability were developed.
\begin{enumerate}
\item[\rm 1.] {\bf Forbidden boundary approach.} The core is as follows. Let ${\mathbf A}$ be stable, i.e. all its eigenvalues be located in the left-hand side of the complex plane. Consider $\{f({\mathbf A})\}$, a family of continuous perturbations of ${\mathbf A}$. Assume that each member $\widetilde{{\mathbf A}} \in \{f({\mathbf A})\}$ does not have eigenvalues with zero real parts, i.e. located on the boundary of the left half-plane. Then all the family $\{f({\mathbf A})\}$ is stable. This reasoning connects the results on zero localization outside imaginary axes (so-called hyperbolicity) with the results on matrix stability.
\item[\rm 2.] {\bf Studying of small-dimensional cases $(n = 1,2,3)$.} The property of $D$-stability is not easy to verify even in the finite-dimensional case.
\item[\rm 3.] {\bf Applications of Gershgorin theorem.} This method is based on the analysis of the perturbations of Gershgorin disks of a matrix $\mathbf A$ under multiplication of $\mathbf A$ by a positive diagonal matrix $\mathbf D$.
\item[\rm 4.] {\bf Polynomial methods,} characterized by applying classical polynomial results to establish relations between spectral properties of matrices and submatrices.
\end{enumerate}
Now see the details.

{\bf 1974 --- Carlson.} The following stability result is established using "forbidden boundary" approach.

\begin{theorem}\cite{CARL2} \label{CARL} Sign-symmetric $P$-matrices are stable.
\end{theorem}

The sketch of the proof consists of the following implications.
\begin{enumerate}
\item[\rm Step 1.] ${\mathbf A}$ is a $P$-matrix $\Rightarrow$ Fisher--Fuller stabilization process leads to a stable matrix ${\mathbf D}_0{\mathbf A}$, where ${\mathbf D}_0$ is a positive diagonal matrix.

\item[\rm Step 2.] Consider the family of continuous perturbations of the form $\{{\mathbf D}_t{\mathbf A}\}$, ${\mathbf D}_t = t{\mathbf I} + (1-t){\mathbf D}_0$, $t \in [0,1]$. Each member $\widetilde{{\mathbf A}} \in \{{\mathbf D}_t{\mathbf A}\}$ is a sign-symmetric $P$-matrix
 $\Rightarrow$ $\widetilde{{\mathbf A}}$ does not have any eigenvalues on the imaginary axes $\Rightarrow$ All matrices from the family $\{{\mathbf D}_t{\mathbf A}\}$ are stable including the initial matrix ${\mathbf A}$.
\end{enumerate}
As a simple corollary, it was established by Johnson that {\it sign-symmetric $P$-matrices are $D$-stable.}

{\bf 1974 ---Johnson.} The problem of $D$-stability was intensively studied by Johnson in several papers \cite{JOHNN}--\cite{JOHN5} which played crucial role in the development of $D$-stability theory.

Probably the most important among this papers is \cite{JOHN1} published in 1974. The paper starts with the following general observation, which outlines the area of searching for the sufficient conditions for $D$-stability.
\begin{theorem}\cite{JOHN1} Any condition on matrices which implies stability and which is preserved under positive diagonal multiplication, is sufficient for $D$-stability.
\end{theorem}

In this paper, Johnson also collected the elementary properties of $D$-stable matrices, that were in fact established earlier in 1950s.
\begin{theorem}[Elementary properties of $D$-stable matrices]\cite{JOHN1} If $\mathbf A$ is $D$-stable then $\mathbf A$ is nonsingular and each of the following matrices are also $D$-stable:
\begin{enumerate}
\item[\rm 1.] ${\mathbf A}^T$
\item[\rm 2.] ${\mathbf A}^{-1}$
\item[\rm 3.] ${\mathbf P}^T{\mathbf A}{\mathbf P},$ where $\mathbf P$ is any permutation matrix.
\item[\rm 4.] ${\mathbf D}{\mathbf A}{\mathbf E}$, where $\mathbf D$, $\mathbf E$ are positive diagonal matrices.
\end{enumerate}
\end{theorem}

The study of the gap between necessary and sufficient conditions for $D$-stability leads to the following natural question:
{\it supposing ${\mathbf A} \in P_0^+$, what additional conditions on $\mathbf A$ would imply $D$-stability?}

Johnson provided a list of sufficient for $D$-stability conditions, collecting already known and establishing new ones. He also showed that {\it none of these conditions are necessary for $D$-stability}.
\begin{theorem}\cite{JOHN1}\label{CLJ} The following matrix classes are $D$-stable.
\begin{enumerate}
\item[\rm 1.] Diagonally stable matrices.
\item[\rm 2.] $M$-matrices
\item[\rm 3.] Strictly diagonally dominant matrices with positive principal diagonal entries.
\item[\rm 4.] Triangular matrices with positive principal diagonal entries.
\item[\rm 5.] Sign-stable matrices.
\item[\rm 6.] Tridiagonal $P$-matrices.
\item[\rm 7.] Oscillatory matrices.
\item[\rm 8.] Hadamard $H$-stable matrices.
\item[\rm 9.] Sign-symmetric $P$-matrices.
\end{enumerate}
\end{theorem}
A number of certain conditions for small-dimensional cases $(n = 2, \ 3, \ 4)$ is analyzed in \cite{JOHN1}. The case $n = 4$ is also considered in \cite{JOHN4}.

In his paper \cite{JOHN5}, Johnson re-stated the problem of $D$-stability in terms of multivariate polynomials. The following statement was proved.
\begin{theorem}[Johnson]\label{JO} Let ${\mathbf A} \in {\mathcal M}^{n \times n}$ be stable. Then $\mathbf A$ is $D$-stable if and only if ${\mathbf A} \pm i{\mathbf D}$ is nonsingular for any positive diagonal matrix $\mathbf D$.
\end{theorem}
 The author considered the real and imaginary parts of $\det({\mathbf A} + i{\mathbf D})$ as multivariate polynomials, and showed that  $D$-stability is equivalent to the property that the system of two multivariate polynomial equations has no positive solution. Thus Johnson made an important conclusion that {\it the $D$-stability of $\mathbf A$ depends entirely on the sequence of principal minors of $\mathbf A$.}
These ideas lead to many characterizations of $D$-stability, e.g. in terms of structured singular values.

{\bf 1975 --- Araki.} $M$-matrices are studied with a view to applications to dynamical systems. The following theorem characterizing diagonal stability was obtained.
\begin{theorem}\cite{AR}\label{AR} A $Z$-matrix is diagonally stable if and only if it is an $M$-matrix.
\end{theorem}

{\bf 1976 --- Goh.} Sufficient condition for the global stability of Lotka--Volterra model of a two species interactions were obtained in \cite{GOH1}. This stability problem leads to the matrix problem of establishing diagonal stability of a $2 \times 2$ matrix. Easy-to-verify sufficient conditions for diagonal stability of a matrix ${\mathbf A} = \{a_{ij}\}_{i,j = 1}^2$ (namely, $a_{11}, a_{22} < 0$, $\det{\mathbf A} > 0$) were established and explained in biological terms. This results gave rise to a number of results connected to the stability of Lotka-Volterra model.

{\bf 1976-1977 --- Cain.} The study of small dimensional cases was continued by Cain. The complete description of $3 \times 3$ real $D$-stable matrices was given in \cite{CA}. In \cite{BC}, the more general problem of characterization the matrices ${\mathbf A} \in {\mathcal M}^{3 \times 3}$ for which ${\rm In}({\mathbf D}{\mathbf A}) = {\rm In}({\mathbf A})$ for any positive diagonal matrix $\mathbf D$, is studied. Relations among matrix minors were obtained by analyzing the characteristic polynomial.

{\bf 1978 --- Cross.} The paper \cite{CROSS} studies both additive and multiplicative $D$-stability, using common approach through Lyapunov equation analysis and diagonal stability concept:
\begin{theorem} \cite{CROSS}
If $\mathbf A$ is diagonally stable then $\mathbf A$ is both multiplicative and additive $D$-stable.
\end{theorem}
The necessary for $D$-stability conditions are analyzed by studying the characteristic polynomial and applying the classical polynomial results.
\begin{theorem} \cite{CROSS} If $\mathbf A$ is additive (multiplicative) $D$-stable then all principal submatrices of $\mathbf A$ are additive (multiplicative) $D$-semistable. If $\mathbf A$ is diagonally stable then all principal submatrices of $\mathbf A$ are diagonally stable.
\end{theorem}
The necessary conditions for $D$-stability and diagonal stability are compared: {\it if $\mathbf A$ is multiplicative or additive $D$-stable then it is a $P_0^+$ matrix, while if $\mathbf A$ is diagonally stable then it is a $P$-matrix}. Necessary and sufficient conditions for all the three stability types for matrices of order $2$ and $3$ were established. For normal and $Z$-matrices, it was shown that all these three stability types are equivalent to stability.

{\bf 1978 --- Barker, Berman and Plemmons.} In \cite{BARBERPL}, the authors studied diagonal stability, starting with the elementary properties of diagonally stable matrices. The following equivalent characterization which gave rise to several criteria of diagonal stability, was obtained.
\begin{theorem}\cite{BARBERPL}\label{BBP} A matrix $\mathbf A$ is diagonally stable if and only if for every nonzero positive semidefinite matrix $\mathbf B$, ${\mathbf B}{\mathbf A}$ has a positive principal diagonal entry.
\end{theorem}
As simple corollaries, the authors derived the criterion of diagonal stability for triangular matrices and for $2 \times 2$ matrices. Some simple necessary for diagonal stability conditions connected to being a $P$-matrix and principal submatrices properties, were obtained. In remarks, the authors listed {\it the matrix classes, for which the conditions of diagonal stability, stability of all principal submatrices (later we refer to this property as to {\it partial stability}) and being a $P$-matrix are equivalent.} These classes include:
\begin{enumerate}
\item[\rm 1.] Triangular matrices;
\item[\rm 2.] $2 \times 2$ matrices;
\item[\rm 3.] $Z$-matrices.
\item[\rm 4.] Symmetric matrices.
\end{enumerate}
For the last two classes, the equivalence of stability, diagonal stability and being a $P$-matrix was pointed out earlier by Cross \cite{CROSS}.

Another important result of this paper is the inductive construction of the positive diagonal solution ${\mathbf D}$ of the Lyapunov equation.

\begin{theorem}\cite{BARBERPL} Let an $n \times n$ matrix ${\mathbf A}$ be partitioned as
$${\mathbf A} = \begin{pmatrix}{\mathbf A}|_{n-1} & \overline{a}_{n} \\ (\underline{a}_{n})^T & a_{nn} \end{pmatrix}.$$
Suppose ${\mathbf A}|_{n-1}$ is diagonally stable, i.e. there is a positive diagonal matrix ${\mathbf D}_{11}$ such that
$${\mathbf D}_{11}{\mathbf A}|_{n-1} + {{\mathbf A}|_{n-1}}^T{\mathbf D}_{11} = {\mathbf W}_{11}, $$
where ${\mathbf W}_{11}$ is $(n-1)\times(n-1)$ symmetric positive definite matrix. Then there exists a positive value $d_{nn}$ such that
$${\mathbf W}:={\mathbf D}{\mathbf A} + {\mathbf A}^T{\mathbf D}$$
is symmetric positive definite for $${\mathbf D} = \begin{pmatrix}{\mathbf D}_{11} & 0 \\ 0 & d_{nn} \end{pmatrix}$$
if and only if
$$a_{nn} > (\overline{a}_{n})^T{{\mathbf A}|_{n-1}}^{-1}{\mathbf D}_{11}\underline{a}_{n} + ((\overline{a}_{n})^T{{\mathbf A}|_{n-1}}^{-1}\overline{a}_{n})(({\mathbf D}_{11}\underline{a}_{n})^T{{\mathbf A}|_{n-1}}^{-1}{\mathbf D}_{11}\underline{a}_{n}). $$
\end{theorem}

{\bf 1977-1978 --- Moylan and Hill.} Positive diagonally dominant and $M$-matrices were studied in \cite{MOY}, their diagonal stability was established.
\begin{theorem}\cite{MOY}
Positive diagonally dominant matrices are diagonally stable.
\end{theorem}
In \cite{MOYH}, diagonal stability was studied from the point of view of applications to the stability of large-scale systems. Sufficient conditions of diagonal stability were analyzed.

{\bf 1978 --- Datta.} The following general question was raised by Datta: "When does stability imply $D$-stability?". Forbidden boundary approach lead to inertia methods for proving $D$-stability, developed in in \cite{DATTA}.
Basing on inertia results and the analysis of Lyapunov equation with the right-hand side of a special form, the following statement was proved.
\begin{theorem}Let ${\mathbf A} \in {\mathcal M}^{n \times n}$ be an upper (lower) Hessenberg matrix with non-zero subdiagonal (superdiagonal). Let there exists a symmetric matrix $\mathbf H$ such that the matrix
    $$-{\mathbf W}:={\mathbf H}{\mathbf A} + {\mathbf A}^{T}{\mathbf H}$$
    is positive semidefinite having a column of the form $w = (\alpha, \ 0, \ \ldots, \ 0)^T, \alpha \neq0$ ($w = ( 0, \ \ldots, \ 0, \ \alpha)^T, \alpha \neq0$, respectively). Then
    \begin{enumerate}
    \item[(i)] $\mathbf A$ does not have any pure imaginary eigenvalues.
    \item[(ii)] ${\rm In}({\mathbf H}) = {\rm In}({\mathbf A})$.
    \item[(iii)] If ${\mathbf H}$ is a positive diagonal matrix, then $\mathbf A$ is $D$-stable.
     \end{enumerate}
    \end{theorem}

He derived from this theorem that {\it a nonderogatory stable matrix $\mathbf A$ in its Schwarz canonical form is $D$-stable} and {\it a nonderogatory stable matrix $\mathbf A$ in its Routh canonical form is totally stable}.

{\bf 1978-1979 --- Berman.} The questions of $D$-stability and diagonal stability, as well as a number of related matrix classes (e.g. $Z$-matrices, $M$-matrices, symmetric, triangular and normal matrices) were studied by Berman and his co-authors (see \cite{BERW}, \cite{BERVW} and the book \cite{BERPL}).

\subsection{1980s. Robust problems}

This decade is characterized by the rapid growth of interest to robust stability problems. New applications of $D$-stability appeared and the question of robustness of $D$-stability was raised. Robustness of $D$-stability and diagonal stability was studied from the point of view of topology in ${\mathcal M}^{n \times n}$. New applications of diagonal stability led to the intensive research in this field, new equivalent characterizations and algorithms for checking diagonal stability. New results, based on the analysis of the Lyapunov equation appeared. The following methods were developed for studying $D$-stability.

\begin{enumerate}
\item[\rm 1.] {\bf Graph-theoretical methods}. There is a variety of applications of graph theory to the study of stability and "super"-stability properties.
\item[\rm 2.] {\bf Generalizations of diagonal dominance conditions}. An example of such a generalization is representing a block partition of a matrix and studying the relations between principal diagonal blocks and the rest of the blocks.
\end{enumerate}

{\bf 1979-1982 --- Khalil and Kokotovic.} In \cite{KHAK2}, \cite{KHAK1}, new applications of $D$-stability were considered (see Section 15.4 for the details). In \cite{KHAL3}, the equivalent characterization of diagonal stability through the convex function minimization was given. Namely, a convex subclass $$ {\mathcal V}:=\{{\mathbf D} \in {\mathcal D}: \ 0 \leq d_{ii} \leq 1\}$$ of the class of positive diagonal matrices was considered and a continuous convex function $$g({\mathbf D}) := \lambda_{max}({\mathbf D}{\mathbf A} + {\mathbf A}^T{\mathbf D})$$ was defined.
\begin{theorem} A matrix $\mathbf A$ is diagonally stable if and only if $$\min_{{\mathbf D} \in {\mathcal V}}g({\mathbf D}) < 0.$$
\end{theorem}

Basing on the above theorem, the algorithm for checking diagonal stability was presented in \cite{KHAL4}.
{\it  For $x \in {\mathbb R}^n$, define ${\mathbf D}:={\rm diag}(x_1, \ \ldots, \ x_n)$. Then $$g(x) := g({\mathbf D})= \lambda_{max}({\mathbf D}{\mathbf A} + {\mathbf A}^T{\mathbf D}) = \max_{v \in V}v^T({\mathbf D}{\mathbf A} + {\mathbf A}^T{\mathbf D})v:=\max_{v \in V}f(x,v),$$
 where $V=\{v \in {\mathbb R}^n:\|v\|=1\}$. Considering $X = \{x \in {\mathbb R}^n: 0 \leq x_i \leq 1\}$, the diagonal stability of $\mathbf A$ is equivalent to the existence of $x\in X$ such that $g(x)<0$.} In fact, some min-max problem is solved on every step.

{\bf 1980 --- Hartfiel.} The study of small perturbations of price stability models led to the study of perturbations of $D$-stable matrices in \cite{HART}.  Given a topology induced by a norm on ${\mathcal M}^{n \times n}$, the necessary and sufficient condition for a $D$-stable matrix to be in the topological interior of the set of $D$-stable matrices were studied. Topological results, which, in modern language, describe robust properties of $D$-stability, were obtained. Hartfiel showed by a counterexample, that the set of $D$-stable matrices is not open (i.e. the property of $D$-stability is not robust) and raised the question how to describe its interior. The following result concerning diagonal stability was obtained.
\begin{theorem} The set of diagonally stable matrices is open in ${\mathcal M}^{n \times n}$
\end{theorem}
As it follows, the property of diagonal stability is robust.
It was shown by a counterexample, that the set of diagonally stable matrices does not coincide with the interior of the set of $D$-stable matrices.

\begin{theorem}\label{HAR} If a matrix $\mathbf A$ lies in the interior of the set of $D$-stable matrices then each principal submatrix of $\mathbf A$ and $({\mathbf A}[\alpha])^{-1}$ is $D$-stable for every $\alpha \subset [n]$.
\end{theorem}
In fact, the above necessary for robust $D$-stability condition was claimed by the author to be necessary and sufficient. However, the sufficiency fails. Later, the corrected version obtained by adding some conditions appeared in the paper by Cain \cite{CA1}. Since the condition of Theorem \ref{HAR} describes the class of totally stable matrices, we have the following inclusions:
$$\mbox{robustly $D$-stable matrices} \subset \mbox{totally stable matrices} \subset \mbox{$D$-stable matrices},$$
where each inclusion is proper.

{\bf 1980 --- Togawa.} The study of the set of $D$-stable matrices from the topological point of view was continued in \cite{TOG}, focusing on its boundary points. Since the set of $D$-stable matrices is neither closed nor open, the boundary points were shown to include those $D$-stable matrices which are not robustly $D$-stable as well as matrices which are $D$-semistable but not $D$-stable. The second case was characterized for $n \leq 4$, while for the first case the following necessary conditions were obtained.

\begin{theorem} If $\mathbf A$ is a $D$-stable matrix which is not robustly $D$-stable then at least one of $k \times k$ principal submatrices of $\mathbf A$ or ${\mathbf A}^{-1}$, is a boundary point of $D$-semistable matrices, for some $k < n$.
\end{theorem}

The set of $D$-semistable matrices was shown to be closed in ${\mathcal M}^{n \times n}$. For $n=4$, the perturbations of the form ${\mathbf A} + t{\mathbf I}$, where $t \geq 0$ (i.e. small diagonal shifts), were considered. The example of a $4 \times 4$ matrix, which is $D$-stable, but loses this property under diagonal shifts, was given (the so-called {\it Togawa matrix}). This fact again shows that $D$-stability is not a robust property, even for the perturbations of such a specific structure.

{\bf 1981 --- Kimura.} On the basis of Metzlerian and Hickian matrices, generalizations of positive diagonally dominant and $M$-matrices (namely, matrices with dominant diagonal blocks) appeared. In (\cite{KIM}), the sufficient conditions of $D$-stability was generalized, using the following construction.

\begin{enumerate}
\item[\rm Step 1.] Given an {\it absolute vector norm} (i.e. a vector norm, satisfying the condition $\|x\| = \||x|\|$, where  $|x| = (|x^1|, \ \ldots, \ |x^n|)$, for any $x \in {\mathbb R}^n$), define the corresponding {\it Minkowski norm} on ${\mathcal M}^{n \times n}$ as follows:
$$\|{\mathbf A}\| = \sup_{x \neq 0}\frac{\|Ax\|_1}{\|x\|_2}.$$ The set of $D$-semistable matrices is shown to be closed in the corresponding topology of ${\mathcal M}^{n \times n}$.
\item[\rm Step 2.] Consider the partition $\alpha$ of the set of indices $[n]$ into $k$ nonempty and pairvise non-intersecting sets $(\alpha_1, \ \ldots, \ \alpha_k)$, union of which covers the whole $[n]$. This partition defines the corresponding partition of the matrix ${\mathbf A}$:
$${\mathbf A} = \{{\mathbf A}_{\alpha_i, \alpha_j}\}, \qquad i,j = 1, \ \ldots, \ k,$$
where ${\mathbf A}_{\alpha_i, \alpha_j}$ is the submatrix of ${\mathbf A}$ lying on the intersection of the rows with indices from $\alpha_i$ and the columns with the indices from $\alpha_j$.
\item[\rm Step 3.] A condition which guarantee $D$-stability is imposed on principal diagonal blocks ${\mathbf A}_{\alpha_i \alpha_i}$. In particular, each $-{\mathbf A}_{\alpha_i \alpha_i}$ is assumed to be an $M$-matrix.
\item[\rm Step 4.] A condition which guarantee $D$-stability is imposed on block matrix $\{{\mathbf A}_{\alpha_i, \alpha_j}\}$, the Minkowski norms of the blocks are used instead of absolute values of the entries. In particular, the condition of {\it generalized diagonal dominance with respect to the given Minkowski norm} is assumed for ${\mathbf A}{\mathbf D}_{\alpha}$, where ${\mathbf D}_{\alpha}$ is a block diagonal matrix (with respect to the partition $\alpha$).
\end{enumerate}

\begin{theorem}\cite{KIM} Given a Minkowski norm induced by an absolute vector norm, an $n \times n$ matrix ${\mathbf A}$ and the partition $\alpha$ of $[n]$. Suppose that every diagonal block ${\mathbf A}_{\alpha_i, \alpha_i}$ $(i,j = 1, \ \ldots, \ k)$ is an $M$-matrix. Then the existence of positive values $d_1, \ \ldots, \ d_k$, satisfying either
$$\sum_{j \neq i} \|A_{\alpha_i\alpha_i}^{-1}\|\|A_{\alpha_i\alpha_j}\|d_j < d_i \qquad i = 1, \ \ldots, \ k.$$
or
$$\sum_{j \neq i} \|A_{\alpha_i\alpha_i}^{-1}|A_{\alpha_i\alpha_j}|\|d_j < d_i \qquad i = 1, \ \ldots, \ k,$$ where $|{\mathbf A}|$ denotes the matrix which consists of the absolute values of the entries of ${\mathbf A}$,
implies that $\mathbf A$ is $D$-stable.
\end{theorem}
\begin{corollary} Let $\mathbf A$ satisfy the conditions of the above theorem and assume further that either $|A_{\alpha_i\alpha_j}| = A_{\alpha_i\alpha_j}$ or $|A_{\alpha_i\alpha_j}| = -A_{\alpha_i\alpha_j}$ for all $j \neq i$, $i = 1, \ \ldots, \ k$. Then $\mathbf A$ is $D$-stable if there exist positive values $d_1, \ \ldots, \ d_k$ such that
$$\sum_{j \neq i} \|A_{\alpha_i\alpha_i}^{-1}A_{\alpha_i\alpha_j}\|d_j < d_i \qquad i = 1, \ \ldots, \ k.$$
\end{corollary}

{\bf 1982 --- Carlson, Datta, Johnson.} The results on diagonal and $D$-stability of tridiagonal matrices, based on the Lyapunov equation analysis were obtained in \cite{CARDJ}.

\begin{theorem} \cite{CARDJ} For a tridiagonal matrix $\mathbf A$, the following conditions are equivalent:
\begin{enumerate}
\item[\rm (i)] $\mathbf A$ is diagonally stable;
\item[\rm (ii)] $\mathbf A$ is totally stable
\item[\rm (iii)] $\mathbf A$ is a $P$-matrix.
\end{enumerate}
\end{theorem}

Note, that $D$-stability of tridiagonal $P$-matrices was pointed out earlier by Johnson, however, in general case $D$-stability does not imply diagonal stability.

A necessary and sufficient criterion of $D$-stability was proved for irreducible tridiagonal $P_0^+$-matrices (see \cite{CARDJ}, p. 301, Theorem 3).

{\bf 1983-1985 --- Berman and Hershkowitz.} In \cite{BERH}, the authors continued the research started in \cite{BERW} and \cite{BERVW}, studying the inclusion relations between diagonally stable, positive stable and $P$-matrices and uniting the results previously obtained in \cite{BERW}, \cite{BERVW}. Some special matrix classes were described through their graph properties. The equivalence between diagonal stability, positive stability and positivity of principal minors was established for special matrix classes ($Z$-matrices, symmetric, triangular). For normal matrices, it was shown by a counterexample, that stability is equivalent to diagonal stability but not equivalent to positivity of principal minors.

The following theorem illustrates the graph-theoretical approach to diagonal-stability analysis.
\begin{theorem}\cite{BERH} If $\mathbf A$ is a $P$-matrix and if the nondirected graph of $\mathbf A$ is a forest then $\mathbf A$ is diagonally stable.
\end{theorem}
Note, that the matrix class mentioned in the above theorem includes Jacobi matrices. The following still open questions on well-known matrix classes were raised: "Are oscillatory or strictly totally positive matrices diagonally stable?"
The study of $D$-stability by graph-theoretic methods was continued in \cite{BERH2}.

{\bf 1984 --- Cain.} "Characterizing $D$-stable matrices is one of the prominent unsolved problems of matrix theory". Topological study of the set of $D$-stable matrices with real and complex entries was continued in \cite{CA1}. This set was shown to be bounded with some complicated algebraic surfaces. However, it was shown that the interiors of the set of $D$-stable and $D$-semistable matrices coincide. The main question considered in \cite{CA1} is which of the $D$-stable matrices are robustly $D$-stable.
\begin{theorem}\cite{CA1} Let $\mathbf A$ be robustly $D$-stable. Then
\begin{enumerate}
\item[\rm 1.] ${\mathbf A}^{-1}$ is robustly $D$-stable;
\item[\rm 2.] Any principal submatrix ${\mathbf A}[\alpha]$, $\alpha \subseteq [n]$ is robustly $D$-stable.
\end{enumerate}
\end{theorem}
Necessary and sufficient conditions for robust $D$-stability in terms of principal submatrices were established.
\begin{theorem}\cite{CA1} Let ${\mathbf A} \in {\mathcal M}^{n \times n}$, $n > 1$ be $D$-stable. Then the following conditions are equivalent.
\begin{enumerate}
\item[\rm 1.] ${\mathbf A}$ is robustly $D$-stable;
\item[\rm 2.] All the matrices of the form $({\mathbf A}[\alpha])^{-1}[\beta]$ are robustly $D$-stable for every $\alpha \subset [n]$, $\beta \subseteq [\alpha] $, $|\beta|< n$.
\item[\rm 3.] All the principal submatrices ${\mathbf A}[\alpha]$ and ${\mathbf A}^{-1}[\alpha]$ are robustly $D$-stable for every $\alpha \subset [n]$, $|\alpha|= n-1$;
\item[\rm 4.] All the principal submatrices ${\mathbf A}[\alpha]$ are robustly $D$-stable and all the principal submatrices ${\mathbf A}^{-1}[\alpha]$ are $D$-stable for every $\alpha \subset [n]$, $|\alpha|< n$;
\item[\rm 5.] All the matrices of the form $({\mathbf A}[\alpha])^{-1}[\beta]$ are $D$-stable for every $\alpha \subset [n]$, $\beta \subseteq [\alpha] $, $|\beta|< n$.
\end{enumerate}
\end{theorem}
This is the correct version of the characterization of robust $D$-stability, claimed earlier by Hartfiel. Note that without assuming $D$-stability of $\mathbf A$, the above statement does not hold.

The equivalent characterization of robust $D$-stability was given for small dimensional cases.
\begin{theorem}\cite{CA1} For $n < 3$, ${\mathbf A}$ is robustly $D$-stable if and only if all the principal submatrices of $\mathbf A$ and ${\mathbf A}^{-1}$ are $D$-stable.
\end{theorem}

{\bf 1985 --- Geromel.} In \cite{GER}, the algorithm proposed by Khalil for finding a positive diagonal solution of the Lyapunov equation, was improved.

{\bf 1985 --- Redheffer.} The problems of finding equivalent characterizations and easy-to-verify sufficient conditions of diagonal stability were considered in \cite{REDH1}, \cite{REDH2}.

In \cite{REDH2}, the following criteria of diagonal stability was proved.

\begin{theorem}\cite{REDH2}\label{RED} Let ${\mathbf A} \in {\mathcal M}^{n \times n}$ be nonsingular, and let ${\mathbf A}|_{n-1}$ and ${\mathbf A}^{-1}|_{n-1}$ denote the leading principal $(n-1)\times(n-1)$ submatrices of $\mathbf A$ and ${\mathbf A}^{-1}$, respectively, obtained by deleting the last row and column. Let $\mathbf D$ be a positive diagonal matrix, ${\mathbf D}|_{n-1}$ be its $(n-1)\times (n-1)$ leading principal matrix. Then
\begin{enumerate}\item[\rm (i)] If ${\mathbf W}:={\mathbf D}{\mathbf A} + {\mathbf A}^T{\mathbf D} > 0$, then $a_{nn} > 0$, ${\mathbf W}|_{n-1}={\mathbf D}|_{n-1}{\mathbf A}|_{n-1} + {\mathbf A}^T|_{n-1}{\mathbf D}|_{n-1} > 0$ and $\widetilde{{\mathbf W}}={\mathbf D}|_{n-1}{\mathbf A}^{-1}|_{n-1} + ({\mathbf A}^{-1})^T|_{n-1}{\mathbf D}|_{n-1} > 0$.
\item[\rm (ii)] If $a_{nn} > 0$, ${\mathbf W}|_{n-1}={\mathbf D}|_{n-1}{\mathbf A}|_{n-1} + {\mathbf A}^T|_{n-1}{\mathbf D}|_{n-1} > 0$ and $\widetilde{{\mathbf W}}={\mathbf D}|_{n-1}{\mathbf A}^{-1}|_{n-1} + ({\mathbf A}^{-1})^T|_{n-1}{\mathbf D}|_{n-1} > 0$, there is $d_{nn} > 0$ such that the extended positive diagonal matrix ${\mathbf D} = {\rm diag}\{{\mathbf D}|_{n-1}, d_{nn}\}$ will satisfy the inequality ${\mathbf D}{\mathbf A} + {\mathbf A}^T{\mathbf D} > 0$.
\end{enumerate}
\end{theorem}
This result, which reduce the study of diagonal stability of an $n \times n$ matrix $\mathbf A$ to the study of (simultaneous) diagonal stability of two matrices of a smaller size, plays an important role in the further development of diagonal stability analysis.

{\bf 1985-1986 --- Abed.} With a view to the new applications of $D$-stability to multi-parameter singular perturbations introduced in \cite{KHAK1}, \cite{KHAK2} (see Section 15), Abed defined the concept of "strong $D$-stability": a matrix ${\mathbf A}$ is called {\it strongly $D$-stable} ({\it robustly $D$-stable}) if $\mathbf A$ is $D$-stable and there is a positive constant $\epsilon > 0$ such that ${\mathbf A} + {\mathbf G}$ is also $D$-stable for each ${\mathbf G} \in {\mathcal M}^{n \times n}$ with $\|{\mathbf G}\| < \epsilon$ (see \cite{AB1}). (Here we use the term "robust $D$-stability", since the term "strong stability" is often used in literature (see, for example, \cite{CROSS}) for additive $D$-stability.) This definition obviously describes the set of interior points of $D$-stable matrices. Similarly, the author introduced robust $D(\alpha)$-stability. Seemingly unaware of the results by Hartfiel, Togawa and Cain, the author showed that $D$-stability is not a robust property (see \cite{AB1} for the corresponding counterexample), but diagonally stable matrices are robustly $D$-stable, and proved the corresponding condition for robust $D(\alpha)$-stability (see \cite{AB1}, Proposition 1). He emphasizes the importance of identifying those classes of $D$-stable matrices for which $D$-stability is a robust property. This notion leads to a variety of problems, where different classes of $D$-stable matrices together with different structures of small perturbations are considered. The applications of robust $D$-stability and $D(\alpha)$-stability were analyzed in \cite{ABT}.

{\bf 1985-1988 --- Hershkowitz and Schneider.} In a number of papers, matrix methods applicable to the study of $D$-stability were developed.
In \cite{HERSS1}, the results on diagonal stability of $H$-matrices (for the definition, see the Appendix) were obtained.  In more details, the following criterion of diagonal stability for $H_+$-matrices was obtained (see \cite{HERSS1}, p. 132, Theorem 4.2).

\begin{theorem} Let $\mathbf A$ be an $H_+$-matrix. Then $\mathbf A$ is diagonally stable if and only if $\mathbf A$ is nonsingular.
\end{theorem}
This result generalizes results of Araki (see Theorem \ref{AR} above). Note, that the class of $H_+$-matrices was considered earlier in \cite{BARBERPL}, where some sufficient for its diagonal stability conditions were pointed out.

Recall, that a matrix $\mathbf A$ is called {\it diagonally semistable} if there is a positive diagonal matrix $\mathbf D$ such that ${\mathbf W}:= {\mathbf D}{\mathbf A} + {\mathbf A}^T{\mathbf D}$ is positive semidefinite. Together with the study of diagonally stable and semistable matrices, Hershkowitz and Schneider introduced the class of {\it diagonally near stable} matrices: a matrix $\mathbf A$ is called {\it diagonally near stable} if $\mathbf A$ is diagonally semistable but not diagonally stable. The authors described diagonally semistable $H_+$-matrices and and determined which $H_+$-matrices are diagonally near stable.

The problem of the uniqueness of Lyapunov scaling factor of diagonally stable and semistable matrices was studied in \cite{HERSS3}, \cite{HERSS2} by using graph-theoretic methods.

{\bf 1987 --- Hu.} Another algorithm for checking diagonal stability was proposed in \cite{HU}. The algorithm is based on solving an infinite system of linear equations.
\begin{theorem} Given an (entry-wise) positive matrix ${\mathbf A} \in {\mathcal M}^{n \times n}$. Then $\mathbf A$ is diagonally stable if and only if the infinite system of linear inequalities
\begin{equation}\label{SYS}
({\mathbf D}(x){\mathbf A}x)^Ty \geq 1 \qquad \mbox{for all} \ y\in S^n,
\end{equation}
where ${\mathbf D}(x) = {\rm diag}\{x_1, \ \ldots, \ x_n\}$, $S^n = \{y \in {\mathbb R}^n: \ y^Ty = 1\}$, has a solution $x \in {\mathbb R}^n$.
Moreover, for any solution $x_0$ of system \eqref{SYS}, ${\mathbf D}(x_0)$ gives a positive diagonal solution to the Lyapunov equation (i.e. a Lyapunov diagonal scaling).
\end{theorem}
If ${\mathbf A}$ is diagonally stable, the proposed algorithm finds the diagonal scaling ${\mathbf D}$ in a finite number of steps.

\subsection{1990s. LMI methods and other new approaches}

In this decade, new equivalent characterizations and algorithms for checking diagonal stability were developed. The algorithms coded into  MATLAB provided a convenient and fast opportunity for checking diagonal stability. Additive $D$-stability was studied from the topological point of view, new criteria of $D$-stability for matrices of a special form were obtained through the Lyapunov equation analysis. Results on $D$-stability and inertia were collected in survey papers. The following new approaches were developed.

\begin{enumerate}
\item[\rm 1.] {\bf Structured singular value approach}. This approach uses powerful tools of the control theory developed for testing robust stability and studying structured uncertainties.
\item[\rm 2.] {\bf LMI-based approach}. This very promising approach gives sufficient for $D$-stability conditions in terms of solving (checking feasibility) of some system of LMIs.
\end{enumerate}

{\bf 1991 --- Kraaijevanger.} Theorem \eqref{BBP} proved in \cite{BARBERPL} leads to the following equivalent characterization of diagonally stable matrices through Hadamard products.

\begin{theorem}\cite{KR} Given a matrix ${\mathbf A} \in {\mathcal M}^{n \times n}$, the following statements are equivalent.
\begin{enumerate}
\item[\rm 1.] ${\mathbf A}$ is diagonally stable;
\item[\rm 2.] ${\mathbf A}\circ{\mathbf S}$ is a $P$-matrix for all symmetric ${\mathbf S}$ satisfying ${\mathbf S} = \{s_{ij}\}$, $s_{ij} \geq 0$, $s_{ii} \neq 0$ (or equivalently, all $s_{ii} = 1$).
\end{enumerate}
\end{theorem}

The proof was based on the property of Hadamard products to preserve positive definiteness. This property together with the analysis of Lyapunov equation gave the following result.
\begin{theorem}\cite{KR} If $\mathbf A$ is diagonally stable then so is ${\mathbf A}\circ{\mathbf S}$ for all symmetric ${\mathbf S}$ satisfying ${\mathbf S} = \{s_{ij}\}$, $s_{ij} \geq 0$, $s_{ii} \neq 0$.
\end{theorem}
Simple necessary and sufficient characterizations of $3 \times 3$ diagonally stable matrices were obtained.

{\bf 1992 --- Hershkowitz.} The basic results on diagonal and $D$-stability were collected in the review paper \cite{HER1}. The relative question was asked: "How far is a $P$-matrix $\mathbf A$ from being stable?" Characterization of multiplicative and additive $D$-stability together with diagonal stability was mentioned among the most important problems of matrix stability (see \cite{HER1}, pp. 162-163). It was demonstrated by an example, that "none of this three types of matrix stability can be characterized by the spectrum of a matrix".
 Problems, related to $D$-stability (e.g. the problem of matrix stabilization) were discussed: "Given a square matrix $\mathbf A$, can we find a diagonal matrix $\mathbf D$ such that the matrix ${\mathbf D}{\mathbf A}$ is stable?" Actually, Fisher--Fuller theorem \ref{FIF} mentioned there gives even more: the positivity of the spectrum of ${\mathbf D}{\mathbf A}$. Basing on results of \cite{BERW} and \cite{CROSS}, the diagramm showing the relations between the classes of stable matrices, was provided. "The problem of characterizing the various types of matrix stability is, in general, a hard open problem and it has been solved only for matrices of order less than or equal to 4". The equivalence results for the class of $Z$-matrices were mentioned.

{\bf 1992 --- Sun.} In \cite{SUN}, perturbations of the class of additive $D$-stable matrices is studied. The notion of {\it total additive $D$-stability} is considered (for the beginning, see \cite{CROSS}): a matrix $\mathbf A$ is called {\it totally additive $D$-stable} if all its principal submatrices are additive $D$-stable. By analogue with multiplicative $D$-stability, it was pointed out, that the set of additive $D$-stable matrices is neither open no closed with respect to the usual topology of ${\mathcal M}^{n \times n}$. It was shown that the closure of the set of additive $D$-stable matrices coincides with the set of totally additive $D$-semistable matrices, while the interior coinsides with the set of totally additive $D$-stable matrices. As in the case of multiplicative $D$-stability, it was shown, that the set of diagonally stable matrices belongs to the interior of the set of additive $D$-stable matrices, but does not coincide with this interior.

{\bf 1994 --- Boyd et al.} Basing on the methods provided in \cite{BGF}, LMI-solvers available to test the feasibility of the LMI ${\mathbf P}{\mathbf A}+{\mathbf A}^T{\mathbf P}\prec0$ with the constraints on the matrix $\mathbf P$ (e.g. ${\mathbf P}$ is positive diagonal) were developed.

 {\bf 1995 --- Chen, Fan and Yu}. An equivalent characterization of $D$-stability in terms of structured singular values was proved in \cite{CFY}. Structured singular values were originally introduced as a tool for studying linear control systems. Here, denote ${\mathcal D}$ the set of all diagonal matrices and ${\mathcal D}^+$ the set of all positive diagonal matrices. Denote $\overline{\sigma}({\mathbf D})$ the largest singular value of a matrix ${\mathbf D} \in {\mathcal D}$.

The {\it real structured singular value $\mu_{{\mathcal D}}({\mathbf A})$} of the matrix ${\mathbf A}$ is defined as $$\mu_{{\mathcal D}}({\mathbf A}) := \frac{1}{\min_{{\mathbf D}\in {\mathcal D}}\{\overline{\sigma}({\mathbf D}): \det({\mathbf I} - {\mathbf A}{\mathbf D})=0\}}$$ and is set equal to $0$, if there is no diagonal matrix ${\mathbf D}$ such that  $\det({\mathbf I} - {\mathbf A}{\mathbf D})=0$.

\begin{theorem}\cite{CFY}
An $n \times n$ matrix $\mathbf A$ is $D$-stable if and only if $\mathbf A$ is stable and
$$\mu_{\mathcal{D}}((j{\mathbf I} + {\mathbf A})^{-1}(j{\mathbf I} - {\mathbf A})) \leq 1. $$
\end{theorem}

This result is based on the equivalent characterization \eqref{JO} by Johnson. Note, that the original statement of \cite{CFY}, in fact, describes totally stable matrices and was later improved by Lee and Edgar \cite{LEE}.

The methods of structured singular values allowed the authors to obtain the necessary condition for $D$-stability

\begin{theorem}\cite{CFY}
An $n \times n$ matrix $\mathbf A$ is $D$-stable only if $\mathbf A$ is stable and
$$\mu_{\mathcal{D}}(({\mathbf I} + {\mathbf A})^{-1}({\mathbf I} - {\mathbf A})) < 1. $$
\end{theorem}
This method was also applied for characterizations of diagonal stability. Denote
$$\hat{\mu}({\mathbf A}) := \inf_{{\mathbf D} \in {\mathcal D}^+}\overline{\sigma}({\mathbf D}{\mathbf A}{\mathbf D}^{-1}). $$

\begin{theorem}\cite{CFY}
An $n \times n$ matrix $\mathbf A$ is Lyapunov diagonally stable if and only if $\mathbf A$ is stable and
$$\hat{\mu}(({\mathbf I} + {\mathbf A})^{-1}({\mathbf I} - {\mathbf A})) < 1. $$
\end{theorem}

The characterization of robust $D$-stability in terms of structured singular values was also obtained.

 {\bf 1997-2001 --- Kanovei and Logofet.}  The following special form of Jacobi matrices was introduced in \cite{LOG2} (see also \cite{ZAL}).
 A matrix $\mathbf A$ is called a {\it Svicobian}, if it can be represented in the following form:
$${\mathbf A} = ({\mathbf S}-{\mathbf Q}){\mathbf D}, $$
where the matrix ${\mathbf S} = {\mathbf F}^T - {\mathbf F}$ is skew-symmetric, $\mathbf Q$ is nonnegative diagonal, ${\mathbf D}$ is positive diagonal. The intersection of multiplicative and additive $D$-stable matrices (so-called $D_aD$-stability) was considered. Sufficient condition of $D_aD$-stability of Svicobians was introduced in terms of associated digraph properties (the Black-White Test).

  For $n = 4$, a verifiable criterion of $D$-stability was proved by Kanovei and Logofet (see \cite{KL}) on the basis of the Routh--Hurwitz criterion. In \cite{KL2}, the elementary properties of $D$-stable matrices were studied. In \cite{KANO}, the following inequality was derived from necessary condition of $D$-stability (i.e. being a $P_0^+$-matrix) and the Routh--Hurwitz criterion.
\begin{theorem} Given a matrix ${\mathbf A} \in {\mathcal M}^{n \times n}$, $n \geq 4$, have no less than two zero elements on the principal diagonal. Then if $-{\mathbf A}$ is $D$-stable we have ${\mathbf A} \in P_0^+$ and for any $i,j,k$ such that $1 \leq i,j,k \leq n$, $i \neq j$, $a_{ii}=a_{jj} = 0$, $a_{kk} \neq 0$, the following minor inequalities hold
$$A[k]A[i,j] \geq A[i,j,k]. $$
\end{theorem}

{\bf 1998 --- Cain et al.} Stable and convergent bounded linear operators in complex Hilbert space were analyzed in \cite{CADHJ}. Diagramms showing relations between the classes of positive definite, stable, $H$-stable, $D$-stable and diagonally stable matrices were provided. For matrices with complex entries, congruence classes with invertible and unitary matrices were studied. It was shown that the classes of positive definite and $H$-stable matrices are invariant under congruent transformation with an invertible matrix.

{\bf 1998 --- Geromel, de Oliveira, Hsu.} A method of checking $D$-stability by solving a number of LMI was proposed in \cite{GEOH}. The problem of $D$-stability was re-stated in the following way.
Consider the set of matrices
$${\mathcal A} = \{{\mathbf A}{\mathbf D}: \ {\mathbf D} \in {\mathcal D}\}, $$
where ${\mathbf A} \in {\mathcal M}^{n \times n},$ ${\mathcal D}$ is a convex polyhedron, spanned by
\begin{equation}\label{Polytope}{\mathbf D}_i = {\rm diag}\{\frac{\epsilon}{n-1}, \ \ldots, \ \frac{\epsilon}{n-1}, \ 1 - \epsilon, \ \frac{\epsilon}{n-1}, \ \ldots, \ \frac{\epsilon}{n-1} \}, \qquad i =1, \ \ldots, \ n, \ \epsilon > 0.\end{equation}
As $\epsilon\rightarrow 0$, we get the set of all positive diagonal matrices.

This approach is based on the concept of {\it transfer functions}.
\begin{enumerate}
\item[\rm Step 1.] By analyzing Lyapunov equation, we get, that if $\mathbf A$ is diagonally stable and ${\mathbf P}$ is the corresponding Lyapunov scaling factor, then ${\mathbf P}{\mathbf D}$ is the Lyapunov scaling factor for ${\mathbf A}{\mathbf D}$.
\item[\rm Step 2.] Consider the rational function of the form
$$T(s):={\mathbf C}(s{\mathbf I}- {\mathbf A})^{-1}{\mathbf B}+{\mathbf D}. $$
It is said to be {\it extended strictly positive real} if it is analytic in ${\mathbb C}^+_0$ and $T(jw)+T(-jw)^T> 0$ for all $w \in [0, + \infty]$. For checking ESPR property, LMI conditions are used.
\item[\rm Step 3.] Given two matrices $\mathbf A$ and $\mathbf B$, the condition of stability of ${\mathbf A}{\mathbf B}$ is equivalent to the ESPR property of the corresponding transfer function $$T(s) = ({\mathbf G}^T{\mathbf A} +{\mathbf H}s)(s{\mathbf I}-{\mathbf B}{\mathbf A})^{-1}$$ for some choice of matrices $\mathbf G$ and $\mathbf H$.
\item[\rm Step 4.] Through the LMI conditions verifying the ESPR property of the corresponding transfer functions, the results on simultaneous stability and Hurwitz stability of $\mathcal A$ for arbitrary convex polytope are obtained.
\end{enumerate}
\begin{theorem} For ${\mathcal B} = {\rm conv}_{N}({\mathbf B})_i$, the following properties are equivalent.
\begin{enumerate}
\item[\rm (i)] The set of matrices $\mathcal A$ is simultaneously stable, that is, there exist a constant positive definite matrix $\mathbf P$ such that ${\mathbf P}{\mathbf A}{\mathbf B} + {\mathbf B}^T{\mathbf A}^T{\mathbf P} \prec 0$ for any ${\mathbf B} \in {\mathcal B}$.
\item[\rm (ii)]There exist a positive definite matrix $\mathbf P$ and arbitrary matrices $\mathbf G$ and $\mathbf H$ satisfying the LMI
$$\begin{pmatrix}{\mathbf G}{\mathbf B}_i+{\mathbf B}_i^T{\mathbf G}^T & {\mathbf P}{\mathbf A} - {\mathbf G} + {\mathbf B}_i^T{\mathbf H}^T \\
{\mathbf A}^T{\mathbf P}-{\mathbf G}^T + {\mathbf H}{\mathbf B}_i & -{\mathbf H}-{\mathbf H}^T \\ \end{pmatrix}\prec 0, \qquad i = 1, \ 2, \ \ldots, \ N.$$
\end{enumerate}
\end{theorem}

Sufficient condition for multiplicative $D$-stability were derived from the following result.

\begin{theorem}For ${\mathcal B} = {\rm conv}_{N}({\mathbf B})_i$, the set of matrices $\mathcal A$ is stable if there exist positive definite matrices ${\mathbf P}_i$, $i = 1, \ 2, \ \ldots, \ n$, and arbitrary matrices $\mathbf G$ and $\mathbf H$ satisfying the LMI
$$\begin{pmatrix}{\mathbf G}{\mathbf B}_i+{\mathbf B}_i^T{\mathbf G}^T & {\mathbf P}_i{\mathbf A} - {\mathbf G} + {\mathbf B}_i^T{\mathbf H}^T \\
{\mathbf A}^T{\mathbf P}_i-{\mathbf G}^T + {\mathbf H}{\mathbf B}_i & -{\mathbf H}-{\mathbf H}^T \\ \end{pmatrix}\prec  0, \qquad i = 1, \ 2, \ \ldots, \ N.$$
\end{theorem}

The above theorem implies the result that diagonally stable matrices are $D$-stable.

{\bf 1999 --- Datta.} A survey paper \cite{DATTA1} included the inertia theorems and their generalizations to arbitrary regions $\mathfrak D$. Applications of inertia results to $D$-stability were mentioned. Main results of \cite{CARDJ} were presented.

{\bf 1999 --- R. Johnson, Tesi.} Basing on geometric characterization of $P_0$ matrices, the authors obtained the following equivalent characterization of $D$-stability closely connected to the results of Johnson \eqref{JO} (see \cite{JOHN5}).
\begin{theorem}\label{RJO} A matrix ${\mathbf A} \in {\mathcal M}^{n \times n}$ is $D$-stable if and only if it is stable and
$$\begin{vmatrix}{\mathbf A} & {\mathbf D} \\ -{\mathbf D} & {\mathbf A} \end{vmatrix}\neq 0$$
for all positive diagonal matrices $\mathbf D$.
\end{theorem}

The stronger condition is shown to be sufficient for $D$-stability.

\begin{theorem} A matrix ${\mathbf A} \in {\mathcal M}^{n \times n}$ is $D$-stable if it is stable and
$$\begin{vmatrix}{\mathbf A} & {\mathbf D}_1 \\ -{\mathbf D}_2 & {\mathbf A} \end{vmatrix}\neq 0$$
for each pair $({\mathbf D_1}, {\mathbf D}_2)$ of positive diagonal matrices.
\end{theorem}

By analyzing conditions, which would guarantee the non-negativity of the coefficients of the corresponding characteristic polynomial, the authors deduced Carlson's theorem \ref{CARL} from the above statement.

 The above results lead to one more sufficient condition.

\begin{theorem}A matrix ${\mathbf A} \in {\mathcal M}^{n \times n}$ is $D$-stable if it is stable and all the coefficients of the polynomial
$$F(d_{11}, \ \ldots, \ d_{nn}) = \begin{vmatrix}{\mathbf A} & {\mathbf D} \\ -{\mathbf D} & {\mathbf A} \end{vmatrix}$$
are nonnegative.
\end{theorem}

The authors also provided certain conditions for robust $D$-stability and for small-dimensional case $n=3$.

\subsection{2000s. New applications and new methods}
In this decade, new applications stimulated research interest to diagonal and additive $D$-stability. The results on different stability types and their applications were collected in the monograph \cite{KAB}. Applications to reaction-diffusion systems led to the new results on additive $D$-stability. More LMI-based conditions for $D$-stability were developed. New applications of diagonal stability led to the study of cyclic matrices. The following new approach to the study of $D$-stability appeared.
\begin{enumerate}
\item[\rm 1.] An approach based on {\bf the Kharitonov criterion} \eqref{KH}. The perturbations of a matrix $\mathbf A$ caused by multiplication by (or addition of) a positive diagonal matrix lead to the specific perturbations of its characteristic polynomial $f_{{\mathbf A}}$. Then we take a bounded subclass of the class of positive diagonal matrices and deduce some estimates on the coefficients of the perturbed polynomial. Using these estimates, we define an interval polynomial. Since all the polynomials of the form $f_{{\mathbf D}{\mathbf A}}$ ($f_{{\mathbf A} + {\mathbf D}}$), where $\mathbf D$ is an arbitrary positive diagonal matrix from the fixed bounded subclass, belong to the obtained interval polynomial, we get sufficient for interval $D$-stability conditions by applying the Kharitonov criterion to it.
\item[\rm 2.] An approach based on {\bf the Kalman--Yakubovich--Popov mathematical apparatus} and other results from dynamical systems theory. Here, we use the transition from the study of a matrix $\mathbf A$ to the study of the corresponding dynamical system $\dot{x} = {\mathbf A}x$.
\item[\rm 3.] An approach based on {\bf the factorization of a matrix class}. This approach is based on the stability criterion by Duan and Patton (see Subsection 2.3, Theorem \ref{DU}) and similar ones. The core idea is the representation of a matrix $\mathbf A$ as a product of two (or more) factors, each of which belongs to a specified matrix class. Then the problem of multiplicative $D$-stability is reduced to the study of invariance of the corresponding matrix classes (and their subclasses) under multiplication by a positive diagonal matrix. This approach allows us to obtain a variety of conditions and new subclasses.
\end{enumerate}

{\bf 2000 --- Kaszkurewicz, Bhaya.} A book \cite{KAB} collects a lot of results on diagonal and $D$-stability, including the classical cases of multiplicative and additive $D$-stability (where the stability region is the left-hand side ${\mathbb C}^-$ of the complex plane) and Schur diagonal and $D$-stability (where the stability region is the unit disk $D(0,1)$), special classes of diagonally stable matrices, a lot of mathematical models, etc.

{\bf 2001 --- Lee, Edgar.} Generalized singular value criterion for robust $D$-stability was proved in \cite{LEE}, improving the corresponding results from \cite{CFY}.

{\bf 2001 --- Wang, Li.} The problem of finding necessary and sufficient conditions of additive $D$-stability was raised again in \cite{WANGL} due to its applications to the stability of reaction-diffusion systems. Searching for these conditions, Wang and Li in fact re-introduced the class of Hicksian matrices under the name of "strict minor conditions" (they defined a matrix $\mathbf A$ to satisfy the minor conditions, if $-{\mathbf A}$ is a $P_0$-matrix and to satisfy strict minor conditions if $-{\mathbf A}$ is a $P$-matrix). The following criteria of additive $D$-stability was proved for $n \leq 3$.
\begin{theorem}\cite{WANGL} Let $n \leq 3$ and an $n \times n$ matrix $\mathbf A$ be stable. Then $\mathbf A$ is additive $D$-stable if and only if $-{\mathbf A}$ is a $P_0$-matrix.
\end{theorem}

Wang and Li also conjectured, that even for the case of an arbitrary $n$, {\it additive $D$-stability is equivalent to being a $P_0$-matrix.}
The necessity part of the above conjecture was pointed out earlier by Cross (see \cite{CROSS}), while the sufficiency part was answered negatively by a counterexample in \cite{SADR}.

However, for the arbitrary dimension $n$, the authors provided criterion of additive $D$-stability imposing some additional conditions on the second additive compound matrix ${\mathbf A}^{[2]}$, based on Lozinsk\~{ii} measures.

{\bf 2002 --- Romanishin, Sinitskii.} An attempt to obtain sufficient conditions for additive $D$-stability on the basis of Kharitonov criterion was made in \cite{RS1}. The authors introduced the class of {\it partially stable (semistable)} matrices (i.e. matrices whose all principal submatrices are stable (semistable)). Partial semistability was again mentioned as a necessary condition for additive $D$-stability.

 The authors made an attempt to reduce the problem of additive $D$-stability of a matrix to the problem of stability of a certain interval polynomial, which can be solved by applying Kharitonov criterion (Theorem \ref{KH}).
For a given $n \times n$ matrix $\mathbf A$ and an arbitrary positive diagonal matrix $\mathbf D$, they wrote the characteristic polynomial of ${\mathbf A} + {\mathbf D}$ as a parameter-dependent function $f_{{\mathbf A}+{\mathbf D}}(\lambda):= f(d_{11}, \ \ldots, \ d_{nn}, \ \lambda)$ and considered its coefficients as increasing functions of positive parameters $d_{ii}$ $(i = 1, \ \ldots, \ n)$. Thus, for bounded $d_{ii}$, we may consider $f_{{\mathbf A}+{\mathbf D}}(\lambda)$ to belong to some interval polynomial.

 However, the main result of \cite{RS1} was disproved by a counterexample in \cite{KOS1}. Note that the counterexample, provided by Kosov for the principal results of \cite{RS1}, uses the Togawa matrix, constructed in \cite{TOG}.

{\bf 2002 --- Kafri.} In the paper \cite{KAF}, it was proved that all 13 conditions sufficient for $D$-stability mentioned by Johnson in \cite{JOHN1} (see Theorem \ref{CLJ} and more) are also sufficient for robust $D$-stability.

{\bf 2003 --- Bhaya, Kaszkurewicz, Santos.} Relations between different stability classes (namely, positive definite matrices (symmetric or not), $D$-stable, simultaneously $D$-stable, diagonally stable, positive diagonally dominant, etc) were studied in \cite{BHK3}, based on the results of \cite{CADHJ}. The core idea was the analysis of the Lyapunov equation which can be regarded as
checking ${\mathbf P}{\mathbf A}+{\mathbf A}^T{\mathbf P}$ for positive definitness.

A new concept generalizing diagonal stability was introduced as follows: a matrix ${\mathbf A}$ is said {\it to belong to the class $\mathcal M$} if the solution $\mathbf P$ of the Lyapunov equation is row diagonally dominant. Positive diagonally dominant matrices were shown to be diagonally stable. The matrices from $\mathcal M$ were shown to be not necessarily $D$-stable, their relation to $D$-stability was posed as an open problem.

Besides a number of known results on implications between classes of stable matrices, the following new classes were considered: namely, the class ${\mathcal W}_{dom}$ of matrices ${\mathbf A} \in {\mathcal M}^{n \times n}$ for which there exists a symmetric positive definite matrix $\mathbf W$ such that ${\mathbf W}{\mathbf A}$ is positive diagonally dominant. Matrices from ${\mathcal W}_{dom}$ are shown to be stabilizable by a diagonal matrix. This gives an interesting link to the Fisher--Fuller result (see Theorem \ref{FIF}), which derives stabilizability from the existence of a sequence of nested positive minors.

The following theorem shows the connection between stabilization by a diagonal matrix and Lyapunov diagonal stability.
\begin{theorem} $\mathbf A$ is stabilizable by a diagonal matrix if and only if there is a symmetric positive definite matrix $\mathbf W$
such that ${\mathbf W}{\mathbf A}$ is diagonally stable. \end{theorem}
Inclusion relations and some description of the classes were written in terms of {\it set products}. Low dimensional characterization for the cases $n = 2$ and $n = 3$ were given. The results of \cite{BHK3} were improved in \cite{CALN}.

{\bf 2005 --- Oliveira and Peres.} Developing methods of \cite{GEOH}, the authors of \cite{OLP} considered the same convex polytope $\mathcal D$, defined by \eqref{Polytope} and represented any $n \times n$ positive diagonal matrix $\mathbf D$ by its coordinates: $${\mathbf D}={\mathbf D}(x) = \sum_{i=1}^nx_i{\mathbf D}_i,$$ where $x = (x_1, \ \ldots, \ x_n)$ is a parameter vector such that $x_i \geq 0$, $\sum_i x_i = 1$. Then $D$-stability of $\mathbf A$ was shown to be equivalent to the existence of a parameter-dependent positive definite matrix ${\mathbf P}(x)$ such that
$${\mathbf P}(x)({\mathbf A}{\mathbf D}(x)) + ({\mathbf A}{\mathbf D}(x))^T{\mathbf P}(x)\prec0 $$
holds for all parameter vectors $x$.

The most obvious, but the most conservative way here was to put ${\mathbf P}(x):={\mathbf P}$ to be constant for all vectors $x$. In \cite{GEOH}, $n$ different positive definite matrices ${\mathbf P}_i$ were found through LMI conditions for the vertices ${\mathbf D}_i$ and two fixed matrices ${\mathbf G}$ and ${\mathbf H}$. Then the matrix ${\mathbf P}(x)$ was defined as their convex combination: $${\mathbf P}(x):= \sum_{i=1}^nx_i{\mathbf P}_i.$$

The key idea of \cite{OLP} was, that matrices $\mathbf G$ and $\mathbf H$ that guarantee the ESPR property of the corresponding transfer function could be chosen different for each $i$, satisfying some additional LMI. Then it was obtained by estimates, that their convex combinations ${\mathbf G}(x)$ and ${\mathbf H}(x)$ satisfy conditions obtained in \cite{GEOH}. This makes the results of \cite{GEOH} less conservative.

\begin{theorem}\cite{OLP} Given a matrix ${\mathbf A} \in {\mathcal M}^{n \times n}$, the matrix ${\mathbf A}{\mathbf D}$ is stable for each ${\mathbf D} \in {\mathcal D} = {\rm conv}_{n}({\mathbf D})_i$,  if there exist positive definite matrices ${\mathbf P}_i$ and arbitrary matrices ${\mathbf G}_i$, ${\mathbf H}_i$, $i = 1, \ 2, \ \ldots, \ n$ satisfying the LMIs
$$\begin{pmatrix}{\mathbf G}_i{\mathbf D}_i+{\mathbf D}_i^T{\mathbf G}_i^T & {\mathbf P}_i{\mathbf A}_i - {\mathbf G}_i + {\mathbf D}_i^T{\mathbf H}_i^T \\
{\mathbf A}^T{\mathbf P}_i-{\mathbf G}_i^T + {\mathbf H}_i{\mathbf D}_i & -{\mathbf H}_i-{\mathbf H}_i^T \\ \end{pmatrix}\prec -{\mathbf I}, \qquad i = 1, \ 2, \ \ldots, \ n.$$
$$\begin{pmatrix} {\mathbf X}_{11}(i,j) & {\mathbf X}_{12}(i,j) \\
{\mathbf X}_{21}(i,j) & {\mathbf X}_{22}(i,j) \end{pmatrix}
\prec \frac{2}{n-1}{\mathbf I}, \qquad i = 1, \ 2, \ \ldots, \ n-1, \ j =i+1, \ \ldots, \ n,$$
where $${\mathbf X}_{11}(i,j):={\mathbf G}_i{\mathbf D}_j+{\mathbf D}_j^T{\mathbf G}_i^T + {\mathbf G}_j{\mathbf D}_i+{\mathbf D}_i^T{\mathbf G}_j^T;$$
$$ {\mathbf X}_{12}(i,j):= ({\mathbf P}_i+{\mathbf P}_j){\mathbf A}_i - {\mathbf G}_i-{\mathbf G}_j + {\mathbf D}_i^T{\mathbf H}_j^T + {\mathbf D}_j^T{\mathbf H}_i^T; $$
$${\mathbf X}_{21}(i,j):={\mathbf A}^T({\mathbf P}_i+{\mathbf P}_j)-{\mathbf G}_i^T-{\mathbf G}_j^T + {\mathbf H}_i{\mathbf D}_j + {\mathbf H}_j{\mathbf D}_i; $$
$$ {\mathbf X}_{22}(i,j):=-{\mathbf H}_i-{\mathbf H}_i^T -{\mathbf H}_j-{\mathbf H}_j^T.$$
Moreover, the Lyapunov factor here is defined by ${\mathbf P}(x)$ for each ${\mathbf A}{\mathbf D}(x)$.
\end{theorem}

{\bf 2005 --- Satnoianu, van den Driessche.} In \cite{SADR}, the application of additive $D$-stability to the reaction-diffusion equations was analyzed, basing on the results of \cite{CROSS}. Using the counterexample by Togawa \cite{TOG}, the authors concluded that {\it not all stable $P$-matrices are additive $D$-stable}, which disproves the conjecture of Wang and Li \cite{WANGL}.

 {\bf 2005 --- Logofet.} A particularly useful review paper \cite{LOG} (see also an earlier book \cite{LOG1}) collected several notions which imply stability (namely, diagonal stability, multiplicative and additive $D$-stability, total stability, sign-stability). This paper was mainly motivated by mathematical ecology problems. Several matrix classes, which are known to belong to some of the stability classes (positive diagonally dominant matrices, $M$-matrices, normal matrices) were studied and a space diagramms (so-called "matrix flower") was provided to show the inclusion relations between the classes of stable matrices. Some applications of $D$-stability to Lotka--Volterra model were presented.

{\bf 2006 --- Arcak, Sontag.} On the basis of the stability result of \cite{TYO} \label{Sect}, the following criterion of diagonal stability for matrices in Form \eqref{STR2} was established in \cite{ARCS1} (see also \cite{ARCS2}).

\begin{theorem} Let ${\mathbf A} \in {\mathcal M}^{n \times n}$ be of the form
$${\mathbf A} = \begin{pmatrix} -\alpha_1 & 0 & \ldots & 0 & - \beta_n \\
\beta_1 & -\alpha_2 & \ddots & \ddots & 0 \\
0 & \beta_2 & -\alpha_3 & \ddots & \vdots \\
\vdots & \ddots & \ddots & \ddots & 0 \\
0 & \ldots & 0 & \beta_{n-1} & - \alpha_n \\ \end{pmatrix}, \qquad \alpha_i > 0, \ \beta_i > 0, \ i = 1, \ \ldots, \ n.$$
Then $\mathbf A$ is diagonally stable if and only if
$$ \frac{\beta_1 \ldots \beta_n}{\alpha_1 \ldots \alpha_n} < \sec(\frac{\pi}{n})^n.$$
\end{theorem}

The proof was based on considering matrices of circulant and skew-circulant structure, obtained from Form \eqref{STR2} and elementary properties of diagonally stable matrices.

{\bf 2007 --- Chu.} In \cite{CHU}, the following problem was considered with the application to Lotka--Volterra model: for which classes of matrices, the equivalence between diagonal stability, total stability and $D$-stability holds?
Note that by definitions, we have the following implications for arbitrary ${\mathbf A} \in {\mathcal M}^{n \times n}$:

$$\mathbf A \  \mbox{is diagonally stable} \ \Rightarrow \ \mathbf A \ \mbox{is totally stable} \Rightarrow$$
  $$\mathbf A \ \mbox{is $D$-stable} \Rightarrow \ \mathbf A \ \mbox{is stable}.$$
For special classes of matrices (e.g. $M$-matrices), the reverse implications hold. In \cite{CHU}, certain sufficient conditions for the reverse implications were established in terms of solvable Lie algebras.

{\bf 2008 --- Cain et al.} In \cite{CALN}, the generalization of the results of \cite{CADHJ} to the case of complex matrices was provided.
The authors described two major approaches to stability:
\begin{enumerate}
\item[\rm 1.] Eigenvalue localization by Gersgorin theorem \eqref{GER};
\item[\rm 2.] Analysis of the solvability of Lyapunov equation \eqref{lyap}.
\end{enumerate}
The authors analyzed the solvability of the Lyapunov equation in a given class $\mathcal P$ and described the class of matrices $\mathcal A$ defined by the following set product: a real matrix $\mathbf A$ is said to belong to the class $\mathcal A$ if and only if there is a positive diagonal matrix $\mathbf D$ such that ${\mathbf A}{\mathbf D}$ is symmetric positive definite. Matrices from ${\mathcal A}$ were shown to be diagonally stable. This results referred and corrected the results from \cite{BHK3}. Sufficient conditions for real sign-symmetric matrices to belong to the class ${\mathcal A}$ were given.

{\bf 2006-2010 --- Shorten, Narendra.} The approach used in \cite{NAS} is based on the existence of the common solution of Lyapunov equation for some matrix family (recall, that given a matrix family $\{{\mathbf A}_i\}_{i=1}^m$, a positive definite matrix $\mathbf P$ is called its {\it common Lyapunov solution} if ${\mathbf P}{\mathbf A}_i + {\mathbf A}_i^T{\mathbf P} = {\mathbf W}_i \prec 0$ for all $i = 1, \ \ldots, \ m$). Shorten and Narendra re-stated the result of Redheffer (see Theorem \ref{RED}) in terms of common Lyapunov solutions.

\begin{theorem} Let an $n \times n$ matrix $\mathbf A$ be stable with negative principal diagonal entries. Let ${\mathbf A}|_{n-1}$ and ${\mathbf A}^{-1}|_{n-1}$ denote the leading principal $(n-1)\times(n-1)$ submatrices of $\mathbf A$ and ${\mathbf A}^{-1}$, respectively, obtained by deleting the last row and column. Then, the matrix $\mathbf A$ is diagonally stable if and only if there is a common diagonal Lyapunov solution for ${\mathbf A}|_{n-1}$ and ${\mathbf A}^{-1}|_{n-1}$.
\end{theorem}

As we see from the initial statement of Theorem \ref{RED}, ${\mathbf D}|_{n-1}$ plays the role of this common diagonal Lyapunov solution.

The main result of Shorten and Narendra reduces the problem of finding a diagonal solution of the Lyapunov equation for an $n \times n$ matrix to the problem of finding a common Lyapunov solution of two $(n-1) \times (n-1)$ matrices, similarly as Redheffer's result does. They use the transition from studying the matrix $\mathbf A$ to studying the corresponding dynamical system $\dot{x} = {\mathbf A}x$. This allows to use the methods, developed in systems theory, namely, the Kalman--Yakubovich--Popov lemma, which gives necessary and sufficient conditions for the system $\dot{x} = {\mathbf A}x$ and its perturbation to share the same Lyapunov function. The LMI form of the Kalman--Yakubovich--Popov result is as follows.

\begin{theorem}[Kalman--Yakubovich--Popov] Let an $n \times n$ matrix $\mathbf A$ allows the following partition:
$${\mathbf A} = \begin{pmatrix}{\mathbf A}_{11} & {\mathbf A}_{12} \\
{\mathbf A}_{21} & {\mathbf A}_{22} \\ \end{pmatrix},$$
where a $k \times k$ submatrix ${\mathbf A}_{11}$ is stable, $(n-k) \times (n-k)$ submatrix ${\mathbf A}_{22}$ is negative definite and ${\mathbf A}_{12}$, ${\mathbf A}_{21}$ are $k \times (n-k)$ and $(n-k) \times k$ submatrices, respectively. Then, there is a $k \times k$ positive definite matrix $\mathbf P$ such that
$$\begin{pmatrix}{\mathbf P} & {\mathbf 0} \\
{\mathbf 0} & {\mathbf I} \\ \end{pmatrix}\begin{pmatrix}{\mathbf A}_{11} & {\mathbf A}_{12} \\
{\mathbf A}_{21} & {\mathbf A}_{22} \\ \end{pmatrix} + \begin{pmatrix}{\mathbf A}_{11} & {\mathbf A}_{12} \\
{\mathbf A}_{21} & {\mathbf A}_{22} \\ \end{pmatrix}^T\begin{pmatrix}{\mathbf P} & {\mathbf 0} \\
{\mathbf 0} & {\mathbf I} \\ \end{pmatrix} \prec 0$$
if and only if ${\rm Re}(H(jw)) > 0$ with $H(jw) = -{\mathbf A}_{21}(jw - {\mathbf A}_{11})^{-1}{\mathbf A}_{12} - {\mathbf A}_{22}$ for all $w \in {\mathbb R}$.
\end{theorem}

Taking into account that $({\mathbf A}^{-1}|_{n-1})^{-1}$ is a rank-one perturbation of ${\mathbf A}|_{n-1}$, they got the following result.

\begin{theorem} Given the following partition of an $n \times n$ matrix $\mathbf A$:
$${\mathbf A} = \begin{pmatrix}{\mathbf A}|_{n-1} & \overline{a}_{n} \\ (\underline{a}_{n})^T & a_{nn} \end{pmatrix},$$
where ${\mathbf A}|_{n-1} \in {\mathcal M}^{(n-1) \times (n-1)}$, $\overline{a}_{n}, \ \underline{a}_{n}^T \in {\mathbb R}^{n-1}$, $a_{nn} \in {\mathbb R}$, $a_{nn} < 0$. Let
$$\widetilde{{\mathbf A}}|_{n-1} : = {\mathbf A}|_{n-1} - \frac{\overline{a}_{n}\underline{a}_{n}^T }{a_{nn}}.$$
Then, $\mathbf A$ is diagonally stable if and only if both ${\mathbf A}|_{n-1}$ and $\widetilde{{\mathbf A}}|_{n-1}$ have a common diagonal Lyapunov solution, i.e. there is a positive diagonal matrix ${\mathbf D}_{n-1} \in {\mathcal M}^{(n-1) \times (n-1)}$ such that ${\mathbf D}_{n-1}{\mathbf A}|_{n-1} + ({\mathbf A}|_{n-1})^T{\mathbf D}_{n-1} \prec 0 $ and ${\mathbf D}_{n-1}\widetilde{{\mathbf A}}|_{n-1} + (\widetilde{{\mathbf A}}|_{n-1})^T{\mathbf D}_{n-1} \prec 0 $.
\end{theorem}

Application to some well-known matrix classes (symmetric, Metzler, cyclic) are given. Basing on the above result, simple test for the stability of Metzler matrices was obtained in \cite{SHN4}. It was established earlier (see \cite{BARBERPL}, \cite{AR}), that for Metzler matrices, stability is equivalent to diagonal stability and to being a Hicksian matrix. However, similar results of Shorten and Narendra allows us to study rank-one perturbations of Metzler matrices.

{\bf 2008-2010 --- Kosov.} The Kharitonov criterion approach was further developed in \cite{KOS1} in order to provide finitely verified sufficient conditions of additive $D$-stability. Though the application of Kharitonov criterion looks a very promising method of the study of $D$-stability, the main question which arises here is: "How to reduce a matrix problem to a polynomial one?" Kosov suggested the following general way (see \cite{KOS1}, p. 767, Theorem 1).

\begin{theorem}\label{KO} Let an $n \times n$ matrix $\mathbf A$ and a map $\Phi: {\mathcal M}^{n \times n} \rightarrow {\mathcal M}^{n \times n}$ satisfy the following conditions.
\begin{enumerate}
\item[\rm (i)] The stability of the matrix $\Phi(\mathbf A) + {\mathbf D}$ implies the stability of ${\mathbf A} + {\mathbf D}$ for any positive diagonal matrix $\mathbf D$.
\item[\rm (ii)] For any $i = 1, \ \ldots, \ n$ there exists a matrix $${\mathbf D}_i^+ = {\rm diag}\{0, \ \ldots, \ 0, \ d_{ii}^+, \ 0, \ \ldots, \ 0\},$$ with the $i$th principal diagonal entry $d_{ii}^+ > 0$ while the rest of the entries are zeroes, such that the matrix
    $\Phi(\mathbf A) + {\mathbf D}_i^+$ is additive $D$-stable.
\item[\rm (iii)] The matrix ${\mathbf A} + {\mathbf D}$ remains stable for any positive diagonal matrix $\mathbf D$ which satisfies $d_{ii} \in [0, d_{ii}^+]$, $i = 1, \ \ldots, \ n$.
\end{enumerate}
Then the matrix $\mathbf A$ is additive $D$-stable.
\end{theorem}

Note, that for a map $\Phi: {\mathcal M}^{n \times n} \rightarrow {\mathcal M}^{n \times n}$, no properties like linearity or continuity are assumed.
As an example of a map $\Phi$, which satisfies the conditions of Theorem \ref{KO}, the author provides the map $W: {\mathcal M}^{n \times n} \rightarrow {\mathcal M}^{n \times n}$, which is defined by the following rule: $$W({\mathbf A}) = {\mathbf A}^{W} = \{a_{ij}^W\}_{i,j = 1}^n, $$
where $$a_{ij}^W = \left\{\begin{array}{cc}a_{ij} & i = j \\ |a_{ij}| & i \neq j \\ \end{array}\right.$$

For any $\mathbf A$ with negative principal diagonal, $-W({\mathbf A})$ gives the comparison matrix of ${\mathbf A}$ (see Appendix).

In \cite{KOS}, Kosov proved the inclusion relations between diagonally stable, partially stable and $P_0$-matrices, as well as properties of principal submatrices of multiplicative (additive) $D$-stable matrices first proved by Cross (see \cite{CROSS}), using the transition to the corresponding systems of differential equations. The author not only disproved the conjecture of \cite{WANGL}, but also showed that {\it a stronger condition of a matrix to be partially stable is also not sufficient for additive $D$-stability}. He gives sufficient condition for diagonal stability of $\mathbf A$, closely connected to those obtained in \cite{HERSS1} for $H$-matrices.
\begin{theorem} For an $n \times n$ matrix $\mathbf A$ to be diagonally stable it suffices that either ${\mathbf A}^{W}$ or $({\mathbf A}^{-1})^{W}$ be stable.
\end{theorem}
Basing on Kharitonov criterion, Kosov proved sufficient conditions for multiplicative and additive interval $D$-stability, introduced in \cite{RS1} (see the definition in Subsection 1.5). The construction he used is as follows.

\begin{enumerate}
\item[\rm Step 1.] For a $P_0$-matrix ${\mathbf A}$ and a matrix parallelepiped of the form
 $$ \Theta = \diag\{d_{ii}, \ \ 0 < d_{ii}^{min} < d_{ii} < d_{ii}^{max} < + \infty, \ \  i  = 1, \ \ldots, \ n\},$$
 we consider the parameter-dependent family $f_{{\mathbf D}{\mathbf A}}(\lambda) = f(d_{11}, \ \ldots, \ d_{nn}, \ \lambda)$. The coefficients of $f_{{\mathbf D}{\mathbf A}}(\lambda)$ are increasing functions of positive parameters $d_{ii}$ $(i = 1, \ \ldots, \ n)$.
 \item[\rm Step 2.] For positive diagonal matrices $${\mathbf D}_{min} = {\rm diag}\{ d_{11}^{min}, \ \ldots, \  d_{nn}^{min}\}$$ and $${\mathbf D}_{max} = {\rm diag}\{ d_{11}^{max}, \ \ldots, \  d_{nn}^{max}\},$$ we construct two corresponding characteristic polynomials $$f_{min}(\lambda) := f_{{\mathbf D}_{min}{\mathbf A}}(\lambda) = \lambda^n + a_1^{min}\lambda^{n-1} + \ldots + a_{n}^{min};$$
 $$f_{max}(\lambda) := f_{{\mathbf D}_{max}{\mathbf A}}(\lambda) = \lambda^n + a_1^{max}\lambda^{n-1} + \ldots + a_{n}^{max}.$$
\item[\rm Step 3.]  Define the interval polynomial
$$ F(\lambda) = \lambda^n + \sum_{i=1}^n[a_i^{min},a_i^{max}]\lambda^{n-i}.$$ It is easy to see, that for each positive diagonal matrix ${\mathbf D} \in \Theta$, the characteristic polynomial $f_{{\mathbf D}{\mathbf A}}(\lambda)$ belongs to $F(\lambda)$.
\end{enumerate}
Recall, that an $n \times n$ matrix $\mathbf A$ is called {\it $D$-stable with respect to $\Theta \subset {\mathcal M}^{n \times n}$} if ${\mathbf D}{\mathbf A}$ is stable for every matrix ${\mathbf D} \in \Theta$.
\begin{theorem}\cite{KOS} For a $P_0$-matrix $\mathbf A$ to be $D$-stable with respect to $\Theta$, it suffices the four Kharitonov polynomials corresponding to $ F(\lambda)$, be stable.
\end{theorem}

An analogous construction was considered for additive $D$-stability.

However, Kosov pointed out that the direct application of Kharitonov criterion "may lead to nonconstructive "hypersufficient" conditions" and provided examples of additive $D$-stable matrices, for which Kharitonov-based criterion fails.
To avoid too rough conditions, he suggested a method of decomposition of the interval polynomial into parts.

{\bf 2009 --- Shorten, Mason, King.} A short proof of the diagonal stability characterization from \cite{BARBERPL} (see Theorem \eqref{BBP}) was given in \cite{SMK}.

{\bf 2009 --- Arcak, Ge.} Generalization and development of results of \cite{ARCS1} on cyclic matrix structures was provided in \cite{GEA}. New sufficient condition for additive $D$-stability was presented, based on the stability criterion proved in \cite{LIW}, in terms of additive compound matrices.

The key idea was to show that, for certain structure of ${\mathbf A}$, its second compound matrix would be Metzler.

\begin{theorem} Let ${\mathbf A} \in {\mathcal M}^{n \times n}$. Then ${\mathbf A}^{[2]}$ is Metzler if and only if ${\mathbf A}$ has the following sign structure:
\begin{equation}\label{STR1}{\mathbf A} = \begin{pmatrix}* & + & 0 & \ldots & 0 & - \\
+ & * & + & \ddots & \ddots & 0 \\
0 & + & * & \ddots & \ddots & \vdots \\
\vdots & \ddots & \ddots & \ddots & \ddots & \vdots \\
0 & \ddots & \ddots & \ddots & \ddots & + \\
- & 0 & \ldots & 0 & + & * \\
 \end{pmatrix},\end{equation}
 where $"+"$ denotes a non-negative entry, $"-"$ denotes a non-positive entry, $"*"$ denotes an entry of an arbitrary sign.
\end{theorem}

\begin{theorem} Let a matrix ${\mathbf A} \in {\mathcal M}^{n \times n}$ is Hurwitz stable and satisfy the following conditions
\begin{enumerate}
\item[\rm 1.] $(-1)^n\det({\mathbf A} - {\mathbf D})> 0$ for every non-negative diagonal matrix ${\mathbf D}$;
\item[\rm 2.] ${\mathbf P}^{(-1)}{\mathbf A}{\mathbf P}$ satisfies sign structure \eqref{STR1} for some $n \times n$ invertible matrix ${\mathbf P}$ with the property that ${\mathbf P}^{(-1)}{\mathbf D}{\mathbf P}$ is a nonnegative diagonal matrix for any nonnegative diagonal matrix ${\mathbf D}$.
Then $\mathbf A$ is additively $D$-stable.
\end{enumerate}
\end{theorem}

The following result (the relaxation of the secant criterion) was provided for cyclic matrices.
\begin{theorem} A cyclic matrix of the form \eqref{STR2} with $a_i > 0$, (while $b_i$ have arbitrary signs) $i = 1, \ \ldots, \ n$ is additively $D$-stable if and only if it is stable.
\end{theorem}

A criterion of additive $D$-stability was proved for matrices of a special block form, where one of the blocks is a cyclic matrix.

{\bf 2009 --- Wimmer.} The paper \cite{WIM3} united the results of \cite{ARCS1} with the matrix forms of Small Gain Theorem proved in \cite{DARW} and provided a unified way of the proof.

 {\bf 2009 --- Burlakova.} The study of small-dimensional cases was continued in \cite{BURL2}. Some necessary and some sufficient conditions of $D$-stability based on the Routh--Hurwitz criterion were obtained in \cite{BURL2} for $n = 5$. These conditions were given in the form of a big number of nonlinear inequalities for the principal minors of a $5 \times 5$ matrix $\mathbf A$.

\subsection{2010s. Recent studies}

Here, we collect recently published results on diagonal and $D$-stability. This decade we characterize by the development of modern and classical approaches to the study of multiplicative and additive $D$-stability. New applications lead to some special cases of robust $D$-stability and to the question when $D$-stability is preserved under rank-one perturbations.

{\bf 2011 --- Arcak.} Continuing the work of \cite{ARCS1}, the following generalization of the secant criterion was obtained in \cite{ARC2}.
\begin{enumerate}
\item[\rm Step 1.] Given a matrix ${\mathbf A} = \{a_{ij}\}_{i,j = 1}^n$ (without loss the generality we assume $a_{ii}= -1$), its principal diagonal entries were excluded and the off-diagonal entries were associated with the weighted digraph $G({\mathbf A})$.
\item[\rm Step 2.] The general case of reducible matrix ${\mathbf A}$ was replaced with the irreducible case by the following result.
\begin{theorem} A reducible matrix ${\mathbf A}$ is diagonally stable if and only if the principal submatrices $\widetilde{{\mathbf A}}_{11}, \ \ldots, \ \widetilde{{\mathbf A}}_{ss}$ in its representation $\widetilde{{\mathbf A}}$ given by
$$\widetilde{{\mathbf A}} = {\mathbf P}{\mathbf A}{\mathbf P}^T = \begin{pmatrix} \widetilde{{\mathbf A}}_{11} & \widetilde{{\mathbf A}}_{12} & \ldots & \widetilde{{\mathbf A}}_{1s} \\
0 & \widetilde{{\mathbf A}}_{22} & \ldots & \widetilde{{\mathbf A}}_{2s} \\
\vdots & \ddots & \ddots & \vdots \\
0 & \ldots & 0 & \widetilde{{\mathbf A}}_{ss} \\ \end{pmatrix}$$
are all diagonally stable.
\end{theorem}
Note that each principal submatrix $\widetilde{{\mathbf A}}_{11}, \ \ldots, \ \widetilde{{\mathbf A}}_{ss}$ corresponds to a strongly connected component of a graph $G({\mathbf A})$.
\item[\rm Step 3.] The case of a single-circuit graph was considered. When $G({\mathbf A})$ consists of a single circuit, ${\mathbf A}$ could be written as
$$ {\mathbf P}{\mathbf A}{\mathbf P}^T = \begin{pmatrix} -1 & 0 & \ldots & \widetilde{a}_{1n} \\
 \widetilde{a}_{2n} & -1 & \ddots & \vdots \\
\vdots & \ddots & \ddots & 0 \\
0 & \ldots &  \widetilde{a}_{n,n-1} & -1 \\ \end{pmatrix},$$
which gave a partial case of the form \eqref{Sect}.
\begin{theorem} A matrix $\mathbf A$ which graph $G({\mathbf A})$ consists of a single circuit is diagonally stable if and only if
$$|\gamma|\Phi({\rm sgn}(\gamma),n)<1 ,$$
where
$$\gamma =\widetilde{a}_{2n} \ldots \widetilde{a}_{n,n-1}\widetilde{a}_{1n}\neq 0 $$
$$\Phi({\rm sgn}(\gamma),n):=\left\{\begin{array}{cc}\cos^n(\frac{\pi}{n}), & \ {\mbox if} \ \gamma<0 \\
1, & \ {\mbox if} \ \gamma>0   \end{array}\right.$$
\end{theorem}
\item[\rm Step 4.] Finally, a more general case was considered. $G({\mathbf A})$ assumed to be a strongly connected graph in which a pair of distinct simple circuits have at most one common vertex (so called {\it cactus structure}).
\begin{theorem} A matrix $\mathbf A$ which graph $G({\mathbf A})$ has a cactus structure with $l$ simple circuits, jth circuit of length $n_j$ traverse the set of vertices $I_j = \{i^j_1, \ \ldots, \ i^j_{n_j}\}$, $1 \leq i_1^j < \ldots <i^j_{n_j} \leq n $, is diagonally stable if and only if there exist constants $\theta_i^j>0$, $i \in I_j$, $j = 1, \ \ldots, \ l$, satisfying
$$|\gamma_j|\Phi({\rm sgn}(\gamma_j),n_j)< \prod_{i \in I_j}\theta_i^j,$$
where
$$\sum_{j \in J_i}\theta_i^j =1, i = 1, \ \ldots, \ n,$$
$J_i =\{j \in [l]: i \in I_j\}$, i.e. the set of circuits the vertex $i$ belongs to.
\end{theorem}
\end{enumerate}

{\bf 2012 --- Kim, Braatz.} The concept of {\it joint $D$-stability} (with respect to an ordered set of matrices) was introduced in \cite{KIMB}. Its relation to diagonal stability is studied.

{\bf 2013 --- Altafini.} In \cite{ALT}, the class of {\it diagonally equipotent} matrices ${\mathbf A}=\{a_{ij}\} \in {\mathcal M}^{n \times n}$, which lie on the boundary of diagonally dominant matrices and is defined by the equalities $$|a_{ii}| = \sum_{i \neq j} |a_{ij}|, \qquad i = 1, \ \ldots, \ n$$ was introduced. Diagonally equipotent matrices were shown to be $H$-matrices. The following criterion of diagonal stability of such matrices was obtained.
\begin{theorem} Let ${\mathbf A} \in {\mathcal M}^{n \times n}$ be irreducible diagonally equipotent with $a_{ii} < 0$, $i = 1, \ \ldots, \ n$. The following conditions are equivalent
\begin{enumerate}
\item[\rm 1.] ${\mathbf A}$ is nonsingular;
\item[\rm 2.] $\Gamma({\mathbf A})$ has at least one negative cycle of length $> 1$;
\item[\rm 3.] $\mathbf A$ is diagonally stable.
\end{enumerate}
\end{theorem}
The same conditions were shown to be equivalent to sign non-singularity and qualitative stability of a diagonally equipotent matrix ${\mathbf A}$.

{\bf 2013-2018 --- Pavani.} In \cite{PAV1}, the following method of checking $D$-stability, based on the results of Johnson \eqref{JO} and R. Johnson and Tesi \eqref{RJO}, was proposed.
\begin{enumerate}
\item[\rm Step 1.] An $n \times n$ matrix $\mathbf A$ was shown to be $D$-stable if and only if $\det({\mathbf A}{\mathbf D}^{-1} + {\mathbf D}{\mathbf A}^{-1}) \neq 0$ for any positive diagonal matrix $\mathbf D$.
\item[\rm Step 2.] After writing $\mathbf D$ in symbolic values (i.e. ${\mathbf D} = \{d_{11}, \ \ldots, \ d_{nn}\} $), the $LU$-decomposition of ${\mathbf A}{\mathbf D}^{-1} + {\mathbf D}{\mathbf A}^{-1}$ was calculated.
\item[\rm Step 3.] All the polynomial factors on $d_{11}, \ \ldots, \ d_{nn}$ that appear on the principal diagonal of the matrix $\mathbf U$ were checked if they have no real roots.
\end{enumerate}
 In \cite{IJP}, a procedure of checking robust $D$-stability of $4 \times 4$ matrices was presented, basing on results from \cite{JOHNT}. This procedure has the same starting point as in \cite{KL}, namely, checking positivity of some cubic polynomials of three variables. It also can be used to determine diagonal stability. The paper \cite{PAV2} was based on the same ideas. The definition of $D$-stability obviously implies the following criterion: {\it a stable matrix $\mathbf A$ is $D$-stable if and only if ${\mathbf A}{\mathbf D}^{-1}$ is not divisible by $x^2 +1$} for any positive diagonal matrix ${\mathbf D}$. Numerical algorithm for calculating the characteristic polynomial of matrix ${\mathbf A}{\mathbf D}$ and its remainder by division by $x^2+1$ was presented as well as numerical examples for $n = 5$.

{\bf 2014 --- Bierkens, Ran.} In \cite{BIR}, the following matrix class, that lies on the boundary of the set of $M$-matrices was considered:
${\mathbf A} \in {\mathcal M}^{n \times n}$ is called a {\it singular $M$-matrix} if ${\mathbf A} = \rho({\mathbf B}){\mathbf I} - {\mathbf B}$, for some (entry-wise) nonnegative matrix $\mathbf B$. For singular $M$-matrices, the following problems were set.

{\bf Problem 1.} {\it When ${\mathbf A} + x\otimes y$ is positive stable?} This problem deals with a partial case of additive ${\mathcal G}$-stability, where the matrix class ${\mathcal G} \subset {\mathcal M}^{n \times n}$ is the class of all matrices of rank one, or one of its subclasses.

{\bf Problem 2.} {\it When ${\mathbf G}({\mathbf A} + x\otimes y)$ is positive stable for all ${\mathbf G} \in {\mathcal G}$?}

The primary interest of the authors was the case when rank-one perturbation is positive.

{\bf Main Problem.} Given a singular $M$-matrix ${\mathbf A}$. Under what conditions ${\mathbf A} + x\otimes y$ is $D$-stable?
In the case of a singular $M$-matrix ${\mathbf A}$, $D$-stability of a rank-one perturbation ${\mathbf A} + x\otimes y$ was shown to be equivalent to the stability of ${\mathbf A} + x\otimes y$.

The methods used in \cite{BIR} for the analysis of the above problems were typical for studying rank-one perturbations.
The algebraic simplicity of zero eigenvalue of ${\mathbf A}$ was shown to be a necessary condition for the stability of ${\mathbf A} + x\otimes y$. Certain sufficient conditions, as well as special classes of singular $M$-matrices ($2$-dimensional, normal, symmetric, etc.) were considered. The conditions when the rank-one perturbation ${\mathbf A} + x\otimes y$ will be a $P$-matrix, were established.

This kind of problems lead to more general problems on robust stability of $D$-stable and diagonally stable matrices under special types of perturbations. One of them is the study of rank-one perturbations of diagonally semistable matrices.

{\bf 2015 --- Giorgi, Zuccotti.} While the review paper \cite{LOG} is based on the problems of mathematical ecology, another review paper on $D$-stability \cite{GIZ} deals with the classical economic motivation of the study of $D$-stable matrices based on \cite{AM}. The paper included not well-known results of Magnani, published in Italian.

{\bf 2016 --- Kushel.} In \cite{KU1}, multiplicative $D$-stability of some known and new matrix classes was established basing on the following criterion.
\begin{theorem} Let an $n \times n$ matrix $\mathbf A$ be a $P$-matrix and $({\mathbf D}{\mathbf A})^2$ be a $Q$-matrix for every positive diagonal matrix $\mathbf D$.
Then $\mathbf A$ is $D$-stable.
\end{theorem}
\begin{corollary} Let $\mathbf A$ be a $P$-matrix. If $\mathbf A$ is strictly row (column) square diagonally dominant for every order of minors, then $\mathbf A$ is $D$-stable.
\end{corollary}

\section{Open problems in $({\mathfrak D}, {\mathcal G}, \circ)$-stability}
\subsection{Preliminaries on binary operations}
Let us recall the following definitions and properties we will use later.
Consider a binary operation $\circ$ on a matrix class ${\mathcal G}_0 \subseteq {\mathcal M}^{n \times n}$:
$$\circ:{\mathcal G}_0 \times {\mathcal G}_0 \rightarrow {\mathcal M}^{n \times n}.$$For the further study, it would be convenient to assume that the class ${\mathcal G}$ belongs to ${\mathcal G}_0$ to avoid the question of spreading the binary operation $\circ$ to the matrices from ${\mathcal G}$. Let us mention the following operation properties (see, for example, \cite{CURT}).
\begin{enumerate}
\item[\rm 1.] Associativity
$${\mathbf A}\circ({\mathbf B}\circ{\mathbf C}) = ({\mathbf A}\circ{\mathbf B})\circ{\mathbf C}$$
for every ${\mathbf A}, {\mathbf B}, {\mathbf C} \in {\mathcal G}_0$.
\item[\rm 2.] There exists an {\it identity element} ${\mathbf L} \in {\mathcal G}_0$:
$${\mathbf L}\circ{\mathbf A} = {\mathbf A}\circ{\mathbf L} = {\mathbf A}$$
for every ${\mathbf A} \in {\mathcal G}_0$.
\item[\rm 3.] There exist {\it inverses}: for every ${\mathbf A} \in {\mathcal G}_0$, there is ${\mathbf A}^{-1} \in {\mathcal G}_0$ such that
$${\mathbf A}\circ{\mathbf A}^{-1} = {\mathbf A}^{-1}\circ{\mathbf A} = {\mathbf L}.$$
\item[\rm 4.] Commutativity
$${\mathbf A}\circ{\mathbf B} = {\mathbf B}\circ{\mathbf A}$$
for every ${\mathbf A}, {\mathbf B} \in {\mathcal G}_0$.
\end{enumerate}

A {\it group} $({\mathcal G}_0, \circ)$ is a set ${\mathcal G}_0$ equipped with a binary operation $\circ$, which satisfies Properties 1-3. If, in addition, the operation $\circ$ satisfies Property 4, a group $({\mathcal G}_0, \circ)$ is called {\it abelian}. If ${\mathcal G} \subset {\mathcal G_0} \subseteq {\mathcal M}^{n \times n}$ and is closed with respect to $\circ$ then $({\mathcal G}, \circ)$ is a subgroup of $({\mathcal G}_0, \circ)$ if $\circ$ satisfies Properties 1-3 on $\mathcal G$. Any matrix property, that is preserved under operation $\circ$, the identity element and inverses with respect to $\circ$ can form a subgroup with respect to $\circ$.

 A group $({\mathcal G}_0, \circ)$ is called {\it topological} if it is a topological space and the group operation $\circ$ is continuous in this topological space. For the definition and theory of topological groups, we refer to \cite{PONT}, for more detailed study of the question see \cite{ART}). In some cases, we will also assume that the class ${\mathcal G}$ forms a subgroup of the topological group ${\mathcal G}_0$ (i.e. a subgroup, which is a closed subspace in the topological space ${\mathcal G}_0$).

We may also consider more theoretical examples, which arises mostly in control theory, i.e. Lyapunov operator (see \cite{BHAT2}) and its generalizations, block Hadamard product (see \cite{HMN}, \cite{CHO}), Redheffer product (see \cite{TIM}), Hurwitz product (for the definition see, for example, \cite{ALAL}), the max-algebra operations (see \cite{BUT}), sub-direct sums (see \cite{FAJ}).

\subsection{Relations to matrix addition and matrix multiplication}

First, let us consider operation of {\bf matrix addition}. Later, we need the following cases.
\begin{enumerate}
\item[\rm 1.] The operation $\circ$ and matrix addition $+$ are connected with the rule of associativity
\begin{equation}\label{+ass}({\mathbf A} \circ {\mathbf B}) + {\mathbf C} = {\mathbf A} \circ ({\mathbf B} + {\mathbf C}).\end{equation}
For example, $\circ$ is also matrix addition.
\item[\rm 2.] The operation $\circ$ is distributive over $+$
\begin{equation}\label{+dist1}({\mathbf A} + {\mathbf B})\circ{\mathbf C} = ({\mathbf A} \circ {\mathbf C}) + ({\mathbf B}\circ{\mathbf C}). \end{equation}
Here, we consider the operations of usual and Hadamard matrix multiplication.
\item[\rm 3.] The operation $+$ is distributive over $\circ$
\begin{equation}\label{+dist2}{\mathbf A} +({\mathbf B}\circ {\mathbf C}) = ({\mathbf A} +{\mathbf B}) \circ ({\mathbf A} +{\mathbf C}).\end{equation}
As an example, we consider the operation of entry-wise maximum.
\end{enumerate}
Note, that for an arbitrary operation $\circ$ none of the above formulae may hold.

{\bf Scalar multiplication}
\begin{enumerate}
\item[\rm 1.] The operation of multiplication by a scalar $\alpha \in {\mathbb R}$ is connected to the operation $\circ$ by the rules of associativity and commutativity:
\begin{equation}\label{scaas} \alpha({\mathbf A}\circ{\mathbf B}) = (\alpha{\mathbf A})\circ{\mathbf B} = {\mathbf A}\circ(\alpha{\mathbf B})\end{equation}
for every ${\mathbf A}, {\mathbf B} \in {\mathcal M}^{n \times n}$.

{\bf Examples.} Matrix multiplication, Hadamard matrix multiplication.
\item[\rm 2.] The operation of scalar multiplication is connected to the operation $\circ$ by the rule of distributivity:
\begin{equation}\label{scadist}\alpha({\mathbf A}\circ{\mathbf B}) = (\alpha{\mathbf A})\circ(\alpha{\mathbf B})\end{equation}
for every ${\mathbf A}, {\mathbf B} \in {\mathcal M}^{n \times n}$.

{\bf Examples.} Matrix addition, matrix maximum.
\end{enumerate}

{\bf Matrix multiplication.}
Now let us consider the relations between $\circ$ and matrix multiplication. Later, we consider one of the following cases:
\begin{enumerate}
\item[\rm 1.] Associativity:
\begin{equation}\label{*ass}{\mathbf A}\circ({\mathbf B}{\mathbf C}) = ({\mathbf A}{\mathbf B})\circ{\mathbf C}\end{equation}
Examples: matrix multiplication.
\item[\rm 2.] Distributivity (the operation $\circ$ is distributive over matrix multiplication):
\begin{equation}\label{*dist1}{\mathbf A}\circ({\mathbf B}{\mathbf C}) = ({\mathbf A}\circ{\mathbf B})({\mathbf A}\circ{\mathbf C})\end{equation}
\item[\rm 3.] Distributivity (matrix multiplication is distributive over $\circ$):
\begin{equation}\label{*dist2}{\mathbf A}({\mathbf B}\circ{\mathbf C}) = ({\mathbf A}{\mathbf B})\circ({\mathbf A}{\mathbf C})\end{equation}
Examples: matrix addition.
\end{enumerate}

\subsection{Open problems}
The problem of defining and studying different cases of $({\mathfrak D},{\mathcal G},\circ)$-stability mainly deals with the properties of the corresponding binary operation $\circ$. Here, we consider the following questions and problems, connected to elementary properties of $({\mathfrak D},{\mathcal G},\circ)$-stable matrices, which will be studied later.

{\bf Problem 1}. Given a binary operation $\circ$ on the set ${\mathcal M}^{n \times n}$ of matrices with real entries, when the equality
$$\sigma({\mathbf A}\circ{\mathbf B}) = \sigma({\mathbf B}\circ{\mathbf A})$$
holds for every ${\mathbf A}, \ {\mathbf B} \in {\mathcal M}^{n \times n}$?

Here, we have the following most obvious cases.
\begin{enumerate}
\item[\rm 1.] When the operation $\circ$ is commutative, we have ${\mathbf A}\circ{\mathbf B} = {\mathbf B}\circ{\mathbf A}$ which implies $\sigma({\mathbf A}\circ{\mathbf B}) = \sigma({\mathbf B}\circ{\mathbf A})$.
\item[\rm 2.] When $\circ$ is matrix multiplication, defined on the set of nonsingular matrices. Then ${\mathbf A}{\mathbf B} = {\mathbf B}^{-1}({\mathbf B}{\mathbf A}){\mathbf B}$ implies $\sigma({\mathbf A}{\mathbf B}) = \sigma({\mathbf B}{\mathbf A})$.
\end{enumerate}

{\bf Problem 2.} Given a binary operation $\circ$ on the set ${\mathcal M}^{n \times n}$, when the equality
       $$({\mathbf A}\circ{\mathbf B})^T = {\mathbf B}^T\circ{\mathbf A}^T$$
holds for every ${\mathbf A}, \ {\mathbf B} \in {\mathcal M}^{n \times n}$?

The above equality obviously holds for matrix addition, matrix multiplication and Hadamard matrix multiplication.

{\bf Problem 3.} Let the operation $\circ$ on ${\mathcal M}^{n \times n}$ be associative and invertible. Given a matrix $\mathbf A$, we have an operation inverse $(\circ{\mathbf A})^{-1}$. Assume, we know the localization of $\sigma(\mathbf A)$ inside a stability region $\mathfrak D$: $$\sigma(\mathbf A) \subset {\mathfrak D}.$$ When we can find a stability region $\widetilde{\mathfrak D}$, dependent on $\mathfrak D$, such that $$\sigma(\circ{\mathbf A})^{-1} \subset \widetilde{\mathfrak D}?$$

More strictly, when we can find a bijective mapping $\varphi:\overline{{\mathbb C}} \rightarrow \overline{{\mathbb C}}$, which connects $\sigma(\mathbf A)$ and $\sigma(\circ{\mathbf A})^{-1}$? Such mappings are well-known for the operations of matrix addition and matrix multiplication.

{\bf Problem 4.} Given a binary operation $\circ$ on the set ${\mathcal M}^{n \times n}$, can we find a rule, connecting $\circ$ to the "usual" operations of matrix multiplication and matrix addition?

As an example, we mention {\it mixed-product property} (see \cite{TSOY}), which connects the operations of Kronecker multiplication $\otimes$ and "usual" matrix multiplication by the equality
$$({\mathbf A} \otimes {\mathbf B})({\mathbf C} \otimes {\mathbf D}) = ({\mathbf A}{\mathbf C})\otimes({\mathbf B}{\mathbf D}), $$
which holds for every ${\mathbf A}, \ {\mathbf B}, \ {\mathbf C}, \ {\mathbf D}   \in {\mathcal M}^{n \times n}$.

\section{Matrix classes and their properties} Here, we collect and analyze the most studied cases of matrix classes $\mathcal G$. We are especially interested in the following basic facts:
\begin{enumerate}
\item[\rm -] the inclusion relations between the studied matrix classes;
\item[\rm -] for a given group operation $\circ$ on ${\mathcal M}^{n \times n}$, if ${\mathcal G}$ is closed with respect to this operation, moreover, if $({\mathcal G}, \circ)$ form a subgroup;
\item[\rm -] commutators of the class ${\mathcal G}$ and transformations that leave this class invariant.
\end{enumerate}

 \begin{enumerate}
\item[\rm 1.] Class ${\mathcal S}$ of symmetric matrices from ${\mathcal M}^{n \times n}$. This matrix class, as well as all its subclasses, is closed with respect to matrix transposition. Class ${\mathcal S}$ equipped with the operation of matrix addition forms a group, nonsingular symmetric matrices form a group with respect to matrix multiplication, and matrices without zero entries form a group with respect to Hadamard matrix multiplication.
To study commutators, we need the following lemma (see \cite{HOJ}, p. 172).
\begin{lemma}\label{COMMU} Let $\mathbf A$, ${\mathbf B}$ be symmetric matrices. Then ${\mathbf A}{\mathbf B}$ is also symmetric if and only if $\mathbf A$ and $\mathbf B$ commute.
\end{lemma}
\item[\rm 2.] Class ${\mathcal H}$ of symmetric positive definite matrices. Here, we mention the following characterization (see \cite{BHAT2}, p. 2): {\it a matrix ${\mathbf A} \in {\mathcal M}^{n \times n}$ is symmetric positive definite if and only if ${\mathbf A} = {\mathbf B}^T{\mathbf B}$ for some matrix ${\mathbf B} \in {\mathcal M}^{n \times n}$}. The class of symmetric positive definite matrices is closed under Hadamard multiplication (first proved in \cite{SHU}, see also \cite{BHAT2}) and under matrix addition. However, for ${\mathbf A}, {\mathbf B} \in {\mathcal H}$, the usual matrix product ${\mathbf A}{\mathbf B}$ belongs to ${\mathcal H}$ if and only if $\mathbf A$ and $\mathbf B$ commute (\cite{BHAT2}). Class ${\mathcal H}$ is also closed with respect to multiplicative inverse. Later, we will use one more equivalent characterization of positive (negative) definiteness: {\it a matrix ${\mathbf A} \in {\mathcal M}^{n \times n}$ is symmetric positive definite if and only if its all its eigenvalues are positive (respectively, negative)}.
\item[\rm 3.] Class ${\mathcal H}_{\alpha}$ of symmetric $\alpha$-diagonal matrices, for a given partition $\alpha = (\alpha_1, \ \ldots \ \alpha_p)$ of $[n]$. Recall, that an $n \times n$ matrix $\mathbf H$ is called an $\alpha$-diagonal matrix if its principal submatrices ${\mathbf H}[\alpha_i]$  (formed by rows and columns with indices from $\alpha_i$, $i = 1, \ \ldots, \ p$) are nonzero, while the rest of the entries is zero. This class forms a group with respect to matrix addition. Later, we consider ${\mathcal H}_{\alpha}$ as a subclass of the class $\mathcal H$ of symmetric positive definite matrices, assuming both block structure and positive definiteness.
\item[\rm 4.] Class $\mathcal D$ of diagonal matrices. This class forms a group with respect to matrix addition. Nonsingular diagonal matrices also form an abelian group with respect to matrix multiplication.
\item[\rm 5.] Sign pattern classes ${\mathcal D}_S$. First, define a {\it sign pattern} ${\rm Sign}({\mathbf D})$ of a diagonal matrix $\mathbf D$ as follows: $${\rm Sign}({\mathbf D}) := {\rm diag}\{{\rm sign}(d_{11}), \ \ldots, \ {\rm sign}(d_{nn})\}.$$
Two diagonal matrices ${\mathbf D}_1$ and ${\mathbf D}_2$ are said to belong to the same sign pattern class if ${\rm Sign}({\mathbf D}_1) = {\rm Sign}({\mathbf D}_2)$. For a given sign pattern $S$, define ${\mathcal D}_S$ as a sign pattern class of diagonal matrices.
 The set of all sign pattern classes covers the set of all diagonal matrices. So we have a proper decomposition
 $${\mathcal G} = \bigcup_{S}{\mathcal G}(S).$$
 Sign pattern classes are studied in connection with matrix inertia properties, $D$-hyperbolicity (see Section 8) and Schur $D$-stability.
 \item[\rm 6.] Class ${\mathcal D}^+$ of positive diagonal matrices. This class forms abelian group with respect to matrix multiplication and is closed with respect to matrix addition.
 \item[\rm 7.] Class ${\mathcal D}_{\alpha}$ of $\alpha$-scalar matrices (resp. ${\mathcal D}^+_{\alpha}$ of positive $\alpha$-scalar matrices). Recall that, for a given partition $\alpha=(\alpha_1, \ \ldots \ \alpha_p)$ of $[n]$, $1 \leq p \leq n$, a diagonal matrix $\mathbf D$ is called an {\it $\alpha$-scalar matrix} if ${\mathbf D}[\alpha_k]$ is a scalar matrix for every $k = 1, \ \ldots, \ p$, i.e.
    $${\mathbf D} = \diag\{d_{11} {\mathbf I}_{\alpha_1}, \ \ldots, \ d_{pp} {\mathbf I}_{\alpha_p}\}.$$
    ${\mathbf D}$ is called a {\it positive $\alpha$-scalar matrix} if, in addition, $d_{ii} > 0$, $i = 1, \ \ldots, \ p$.
    For studying this matrix class, see \cite{HERM}, \cite{WAN}, \cite{KHAK2} and \cite{AB1}, \cite{AB2}.
    For a fixed $\alpha$, the class ${\mathcal D}_{\alpha}$ is closed with respect to matrix addition and forms abelian group with respect to matrix multiplication.
\item[\rm 8.] Class ${\mathcal D}_{\tau}$ of positive diagonal matrices ordered with respect to a given permutation $\tau \in \Theta_{[n]}$ and a union ${\mathcal D}_{\tau_1, \ldots, \tau_k}$ of $k$ classes, defined by the permutations $\tau_1, \ \ldots, \ \tau_k$. (Recall that a positive diagonal matrix ${\mathbf D} = \diag\{d_{11}, \ \ldots, \ d_{nn}\}$ is called {\it ordered with respect to a permutation $\tau = (\tau(1), \ \ldots, \ \tau(n))$ of $[n]$}, or {\it $\tau$-ordered}, if it satisfies the inequalities $$d_{\tau(i)\tau(i)} \geq d_{\tau(i+1)\tau(i+1)}, \qquad i = 1, \ \ldots, \ n-1. $$
     In what follows, $\Theta_{[n]}$ denotes, as usual, the set of all the permutations of $[n]$. Obviously, $${\mathcal D} = \bigcup_{\tau \in \Theta_{[n]}}{\mathcal D}_{\tau}.$$
     This class is closed with respect to matrix addition and matrix multiplication, however, does not contain multiplicative inverses.
\item[\rm 9.] Class ${\mathcal D}_{\Theta}$ of diagonal matrices satisfying the inequalities
$$\Theta = \prod(d_{ii}^{min}, \ d_{ii}^{max})= $$ $$ \diag (d_{ii}, \ 0 < d_{ii}^{min} < d_{ii} < d_{ii}^{max} < +\infty, \ \ i = 1, \ 2, \ \ldots \ n); $$
and class ${\mathcal D}_{\Theta_0}$, where $$\Theta_0 = \prod(0, \ d_{ii}^{max})= $$ $$ \diag (d_{ii}, \ 0 < d_{ii} < d_{ii}^{max} < +\infty, \ \ i = 1, \ 2, \ \ldots \ n). $$
\item[\rm 10.] Class ${\mathcal D}_V$ of vertex diagonal matrices. Recall that a real diagonal matrix $\mathbf D$ is called {\it vertex diagonal} if $|{\mathbf D}|=1 $, i.e. $|d_{ii}| = 1$ for any $i = 1, \ \ldots, \ n$ (\cite{BHK}).
\item[\rm 11.]
 Now let us consider the following characterizations of a $\tau$-ordered matrix $D$:
 $$d_{max}({\mathbf D}) := \max_i\dfrac{d_{\tau(i)\tau(i)}}{d_{\tau(i+1)\tau(i+1)}}, \qquad i = 1, \ \ldots, \ n.$$
$$d_{min}({\mathbf D}) := \min_i\dfrac{d_{\tau(i)\tau(i)}}{d_{\tau(i+1)\tau(i+1)}}, \qquad i = 1, \ \ldots, \ n.$$
According to the definition, $1 \leq d_{min}({\mathbf D}) \leq d_{max}({\mathbf D})$ for every positive diagonal matrix $\mathbf D$.
\end{enumerate}
Here, we have the following chains of inclusions:
\begin{equation}\label{chain}
{\mathcal D}_{\alpha}^+ \subset {\mathcal D}^+ \subset {\mathcal H}_{\alpha} \subset {\mathcal H}.
\end{equation}

Now let us consider the following pairwise commuting subclasses of ${\mathcal H}$:
\begin{enumerate}
\item[\rm 1.] $({\mathcal H}, {\mathcal I})$, where the class ${\mathcal I}$ consists of the only one identity matrix $\mathbf I$.
\item[\rm 2.] $({\mathcal H}_{\alpha}, {\mathcal D}_{\alpha})$ ($\alpha$-scalar diagonal matrices commute with $\alpha$-block symmetric matrices).
\item[\rm 3.] $({\mathcal D}, {\mathcal D})$ (diagonal matrices commute within themselves).
\end{enumerate}
\part{Methods of analysis}

\section{Inclusion relations and topological properties} Here, we collect basic statements describing the class of $({\mathfrak D}, \ {\mathcal G}, \ \circ)$-stable matrices.
\subsection{General statements}
First, let us mention a topological property, which shows if the definition of $({\mathfrak D}, \ {\mathcal G}, \ \circ)$-stability is meaningful, i.e. if the defined class of $({\mathfrak D}, \ {\mathcal G}, \ \circ)$-stable matrices is nonempty.

\begin{theorem} Given a bounded (in absolute value) stability region ${\mathfrak D} \subseteq {\mathbb C}$, a matrix class ${\mathcal G} \subset {\mathcal M}^{n \times n}$ and a continuous binary operation $\circ$ on ${\mathcal M}^{n \times n} \times {\mathcal M}^{n \times n}$. Then, for the class of $({\mathfrak D},{\mathcal G},\circ)$-stable matrices to be nonempty it is necessary the class ${\mathcal G}$ be bounded in ${\mathcal M}^{n \times n}$.
\end{theorem}
{\bf Proof.} Let the matrix class $\mathcal G$ be unbounded, i.e. there exists a sequence $\{{\mathbf G}_i\}_{i=1}^{\infty}$ from ${\mathcal G}$ with $\|{\mathbf G}_i\| \rightarrow \infty$. Let there exist at least one $({\mathfrak D}, \ {\mathcal G}, \ \circ)$-stable matrix $\mathbf A$. By definition,  $\sigma({\mathbf G}\circ {\mathbf A}) \subset {\mathfrak D}$ for all ${\mathbf G} \in {\mathcal G}$. Since ${\mathfrak D}$ is bounded in $\mathbb C$, there is a positive value $R$ such that $|\lambda| \leq R$ for all $\lambda \in {\mathfrak D}$. Thus the spectral radius $\rho({\mathbf G}\circ {\mathbf A}) \leq R$ for all ${\mathbf G} \in {\mathcal G}$. By the continuity of the operation $\circ$, we have
 $\|{\mathbf G}_i\circ {\mathbf A}\|\rightarrow\infty $ as $\|{\mathbf G}_i\| \rightarrow \infty$ for any fixed $\mathbf A$, which implies
 $\rho({\mathbf G}_i\circ {\mathbf A}) \rightarrow \infty$ as well. Thus we came to a contradiction. $\square$

Note that all the operations on ${\mathcal M}^{n \times n}$ considered in Section 1 (namely, matrix addition, matrix multiplication and Hadamard matrix multiplication as well as block Hadamard multiplication) are continuous with respect to the usual topology of ${\mathcal M}^{n \times n}$. Thus, considering $({\mathfrak D},{\mathcal G},\circ)$-stability with respect to any bounded stability region ${\mathfrak D} \subseteq {\mathbb C}$ and any of the above operation, we must consider a bounded matrix class $\mathcal G$. As an example here, we may consider Schur $D$-stability. At the same time, if ${\mathfrak D}$ is unbounded, we may consider an unbounded matrix class $\mathcal G$ as well as bounded.

The next results deal with basic inclusion relations between different classes of $({\mathfrak D},{\mathcal G},\circ)$-stable matrices.

\begin{theorem}\label{mainob1} Let ${\mathfrak D}\subset {\mathbb C}$ be a stability region of the complex plane, ${\mathcal G} \subset {\mathcal M}^{n \times n}$ be a matrix class and $\circ$ be a binary operation on ${\mathcal M}^{n \times n}$. Then
     \begin{enumerate}
\item[\rm 1.] For any subset ${\mathfrak D}_1$ of the complex plane such that ${\mathfrak D}_1 \subseteq {\mathfrak D}$, the class of $({\mathfrak D}_1, \ {\mathcal G}, \ \circ)$-stable matrices belongs to the class of $({\mathfrak D}, \ {\mathcal G}, \ \circ)$-stable matrices.
\item[\rm 2.] For any matrix class ${\mathcal G}_1$ such that ${\mathcal G}_1 \subseteq {\mathcal G}$, conversely, the class of $({\mathfrak D}, \ {\mathcal G}, \ \circ)$-stable matrices belongs to the class of $({\mathfrak D}, \ {\mathcal G}_1, \ \circ)$-stable matrices. In particular, if $${\mathcal G} = \bigcup_{i \in I}{\mathcal G}_{i},$$
a matrix $\mathbf A$ is $({\mathfrak D}, \ {\mathcal G}, \ \circ)$-stable if and only if it is $({\mathfrak D}, \ {\mathcal G}_i, \ \circ)$-stable for any $i \in I$.
\end{enumerate}
     \end{theorem}
  {\bf Proof.} The proof obviously follows from the definition of $({\mathfrak D}, \ {\mathcal G}_1, \ \circ)$-stability.

Now we make the most common observation which may be used describing the class of $({\mathfrak D},{\mathcal G},\circ)$-stable matrices.
 \begin{theorem}\label{mob} Given a stability region ${\mathfrak D} \subset {\mathbb C}$, a matrix class ${\mathcal G}\subset{\mathcal M}^{n \times n}$ and a binary operation $\circ$ on ${\mathcal M}^{n \times n}$,
any condition on matrices which implies ${\mathfrak D}$-stability and which is preserved under ${\mathbf G}\circ(\cdot)$ for any $\mathbf G$ from the class ${\mathcal G}$, is sufficient for $({\mathfrak D},{\mathcal G},\circ)$-stability.
\end{theorem}
The proof is obvious.
 This is a generalization of a simple observation, made by Johnson for the case of multiplicative $D$-stability (see \cite{JOHN1}, p. 54, Observation (i)).

However, in each special case, there may be other classes, not covered by this reasoning.

\subsection{Inclusion relations between stability regions} Let ${\mathcal G}$ be the class of positive diagonal matrices, $\circ$ be matrix multiplication. The following inclusion relations between $({\mathfrak D}, \ {\mathcal G}, \ \circ)$-stability classes are based on the inclusion relations between the corresponding stability regions (the positive direction of the real axes ${\mathbb R}^+$ belongs to the open right half plane of the complex plane ${\mathbb C}^+$ which belongs to the complex plane without the imaginary axes ${\mathbb C}\setminus{\mathbb I}$):
\begin{equation}\label{Incl}\mbox{$D$-positive matrices} \subset \mbox{$D$-stable matrices} \subset \mbox{$D$-hyperbolic matrices}.\end{equation}
Let us consider a sector around the positive direction of the real axes with an inner angle $2\theta$, $0 < \theta < \frac{\pi}{2}$:
$${\mathbb C}^+_{\theta}:=\{z \in {\mathbb C}: |\arg(z)| < \theta\}.$$
The following inclusions between stability regions
$${\mathbb R}^+ \subset {\mathbb C}^+_{\theta} \subset {\mathbb C}^+ \subset {\mathbb C}\setminus(- {\mathbb C}^+_{\theta}) \subset {\mathbb C}\setminus{\mathbb R}^-$$
imply the corresponding inclusions between stability classes, where $D$ denotes the class of positive diagonal matrices:
$$\mbox{$D$-positive matrices} \subset \mbox{$({\mathbb C}^+_{\theta}, D)$-stable matrices} \subset \mbox{$D$-stable matrices} \subset$$
$$ \mbox{$(({\mathbb C}\setminus(- {\mathbb C}^+_{\theta})), D)$-stable matrices} \subset \mbox{$({\mathbb C}\setminus{\mathbb R}^-, D)$-stable matrices}.$$

It is easy to see that the last class of $({\mathbb C}\setminus{\mathbb R}^-, D)$-stable matrices coincides with the class of $P_0^+$-matrices.

Let us introduce the following notations for the matrix classes.
\begin{enumerate}
\item[\rm 1.] ${\mathbb D}_{H}$ for the class of $D$-hyperbolic matrices;
\item[\rm 2.] ${\mathbb D}^-_{\theta}$ for $(({\mathbb C}\setminus(- {\mathbb C}^+_{\theta})), D)$-stable matrices;
\item[\rm 3.] ${\mathbb D}_C$ for $D$-stable matrices;
\item[\rm 4.] ${\mathbb D}^+_{\theta}$ for $({\mathbb C}^+_{\theta}, D)$-stable matrices;
\item[\rm 5.] ${\mathbb D}^+$ for $D$-positive matrices;
\item[\rm 6.] ${\mathcal P}_0^+$ for $P_0^+$-matrices.
\end{enumerate}

\begin{center}
\begin{tikzpicture}
  \begin{scope}[]
    \draw (120:1.3) circle (3);
     \draw (110:1.3) circle (2.2);
    \draw  (60:1.3) circle (2.5);
    \draw ( 90:1.2) circle (1.5);
    \draw ( 90:1.2) circle (0.9);
    \draw ( 90:1.2) circle (0.4);
 \end{scope}
  \node at ( 60:-1.6)    {${\mathcal P}_0^+$};
  \node at (160:2.3)    {${\mathbb D}^-_{\theta}$};
  \node at (20:3)    {${\mathbb D}_{H}$};
 \node at (90:2.3)    {${\mathbb D}_C$};
 \node at (90:0.5)    {${\mathbb D}^+_{\theta}$};
 \node at (90:1.2)    { ${\mathbb D}^+$};
\end{tikzpicture}
\end{center}
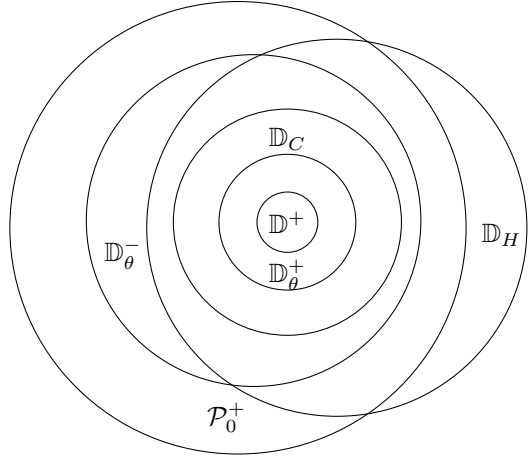
\begin{center} Figure 1. Inclusion relations between $({\mathfrak D}, D)$-stability classes, determined by inclusion relations between stability regions $\mathfrak D$.
\end{center}

Now consider the following generalization of Sequence \eqref{Incl} of inclusion relations. Let ${\mathfrak D} \subset {\mathbb C}$ be an open stability region, denote its boundary by $\partial({\mathfrak D})$. The following concept was introduced in \cite{BHPM}: a matrix ${\mathbf A} \in {\mathcal M}^{n \times n}$ is called {\it $\partial({\mathfrak D})$-singular} if $\sigma(\mathbf A) \bigcap \partial({\mathfrak D}) \neq \emptyset$ and {\it $\partial({\mathfrak D})$-regular} if $\sigma(\mathbf A) \bigcap \partial({\mathfrak D}) = \emptyset$ (for the related studies, see \cite{ASCM}, \cite{BAP} -\cite{BHPM}). We define the concept of $(\partial({\mathfrak D}), {\mathcal G}, \circ)$-regularity basing on $\partial({\mathfrak D})$-regularity: a matrix ${\mathbf A} \in {\mathcal M}^{n \times n}$ is called {\it $(\partial({\mathfrak D}), {\mathcal G}, \circ)$-regular} if $\sigma({\mathbf G}\circ{\mathbf A})\bigcap \partial({\mathfrak D}) = \emptyset$ for any matrix ${\mathbf G}$ from the given matrix class $\mathcal G$.

 For any subset ${\mathfrak D}_1 \subseteq {\mathfrak D}$ we have
$${\mathfrak D}_1 \subseteq {\mathfrak D} \subset ({\mathbb C} \setminus \partial(\mathfrak D)).$$
Thus the following inclusions hold:
$$\mbox{$({\mathfrak D}_1, {\mathcal G}, \circ)$-stable matrices} \subset \mbox{$({\mathfrak D}, {\mathcal G}, \circ)$-stable matrices} \subset \mbox{$(\partial({\mathfrak D}), {\mathcal G}, \circ)$-regular matrices}.$$

\subsection{Inclusion relations between matrix classes} The following relations are based on Inclusion chain \eqref{chain} between matrix classes.
$$\mbox{$H$-stable matrices} \subset \mbox{$H_{\alpha}$-stable matrices} \subset \mbox{$D$-stable matrices} \subset $$ $$\mbox{$D_{\alpha}$-stable matrices} \subset \mbox{stable matrices}.$$
The relation $$\mbox{$H$-stable matrices} \subset \mbox{$D$-stable matrices}$$ is pointed out in \cite{AM}.

 Let $\alpha$ and $\beta$ be two partitions of the set $[n]$, such that $\beta \subseteq \alpha$. Then
$$\mbox{$H_{\beta}$-stable matrices} \subset \mbox{$H_{\alpha}$-stable matrices} $$
and
$$\mbox{$D_{\alpha}$-stable matrices} \subset \mbox{$D_{\beta}$-stable matrices}. $$

The inclusion $$\mbox{additive $H_{\beta}$-stable matrices} \subset \mbox{additive $H_{\alpha}$-stable matrices} $$
was pointed out in \cite{GUH1} (see \cite{GUH1}, p. 327, Theorem 2.1).

Further,
$$\mbox{$D_{\tau}$-stable matrices} \subset \mbox{$D$-stable matrices},$$
for any permutation $\tau$ of the set $[n]$. Let us consider the decompositions of the class of $D$-stable matrices. The following statement was proved in \cite{KU1}.
 \begin{lemma}\label{Dst}
A matrix $\mathbf A$ is $D$-stable if and only if it is $D_\tau$-stable for any $\tau \in \Theta_{[n]}$.
\end{lemma}

Let us introduce the following notations for the matrix classes.
\begin{enumerate}
\item[\rm 1.] $\mathbb S$ for the class of stable matrices;
\item[\rm 2.] ${\mathbb D}_{\tau}$ for the class of $D_{\tau}$-stable matrices;
\item[\rm 3.] ${\mathbb D}_{\alpha}$ for $D_{\alpha}$-stable matrices;
\item[\rm 4.] ${\mathbb H}_{\alpha}$ for $H_{\alpha}$-stable matrices;
\item[\rm 5.] ${\mathbb H}_C$ for $H$-stable matrices;
\end{enumerate}

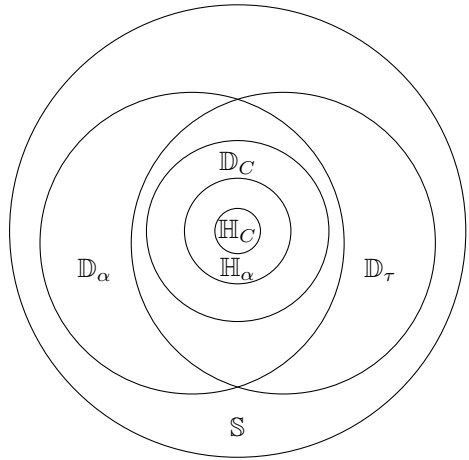
\begin{center}
\begin{tikzpicture}
  \begin{scope}[]
      \draw ( 90:1.2) circle (3);
    \draw (120:1.2) circle (2);
    \draw  (60:1.2) circle (2);
    \draw ( 90:1.2) circle (1.2);
    \draw ( 90:1.2) circle (0.7);
    \draw ( 90:1.2) circle (0.3);
 \end{scope}
  \node at ( 90:-1.4)    {$\mathbb S$};
  \node at (160:2)    {${\mathbb D}_{\alpha}$};
  \node at (20:2)    {${\mathbb D}_{\tau}$};
 \node at (90:2.1)    {${\mathbb D}_C$};
 \node at (90:0.7)    {${\mathbb H}_{\alpha}$};
 \node at (90:1.2)    {${\mathbb H}_C$};
\end{tikzpicture}
\end{center}
\begin{center} Figure 2. Inclusion relations between $({\mathbb C}^-, {\mathcal G})$-stability classes, determined by inclusion relations between matrix classes ${\mathcal G}$.
\end{center}

The class of Schur $D$-stable matrices belongs to the class of vertex $D$-stable matrices. In some cases ($n =3$ (see \cite{FL1}, p. 20, Theorem 2.7), tridiagonal matrices (see \cite {FL2}, p. 46, Theorem 4.2), some others) these classes coincide.
The idea of using the transition from the studying a convex matrix polyhedron spanned by a finite number of vertices to the study of its vertices and edges is widely used for analyzing different stability types (e.g. the stability of interval matrices). It goes back to Kharitonov theorem and its generalizations.
 Consider the following generalization of the concept of vertex $D$-stability. Given a stability region ${\mathfrak D}$, a binary matrix operation $\circ$ and a convex polyhedron ${\mathcal B} \subset {\mathcal M}^{n \times n}$, spanned by $N$ vertices ${\mathbf B}_1, \ \ldots, \ {\mathbf B}_N$. A matrix ${\mathbf A} \in {\mathcal M}^{n \times n}$ is called {\it vertex $({\mathfrak D}, \circ)$-stable} if ${\mathbf B}_i\circ{\mathbf A}$ is ${\mathfrak D}$-stable for each $i = 1, \ldots, n$ and is called {\it edge $({\mathfrak D}, \circ)$-stable} if ${\mathbf B}\circ{\mathbf A}$ is ${\mathfrak D}$-stable for each ${\mathbf B}$ belongs to one of the edges of ${\mathcal B}$. Thus
$$\mbox{vertex $({\mathfrak D}, \circ)$-stable matrices} \subset \mbox{edge $({\mathfrak D}, \circ)$-stable matrices} \subset $$ $$\mbox{ $({\mathfrak D}, {\mathcal B}, \circ)$-stable matrices}. $$
This concept is also applied for studying multiplicative and additive $D$-stability. For this, the class of positive diagonal matrices ${\mathcal D}^+$ is either approximated by a convex polyhedron \eqref{Polytope} or replaced with a bounded class of diagonal interval matrices (see the concept of interval $D$-stability).

Finally, let us consider the following result on the belonging to the class of $({\mathfrak D},{\mathcal G},\circ)$-stable matrices and the connection to ${\mathfrak D}$-stability.

\begin{theorem}\label{Dstab} Let $\mathfrak D$ be an arbitrary stability region and $\circ$ be a group operation on ${\mathcal M}^{n \times n}$. Then
 \begin{enumerate}
 \item[\rm 1.] If the matrix class $\mathcal G \subset {\mathcal M}^{n \times n}$ is closed with respect to the operation $\circ$, then $\mathbf A$ is $({\mathfrak D},{\mathcal G},\circ)$-stable implies ${\mathbf G} \circ {\mathbf A}$ is $({\mathfrak D},{\mathcal G},\circ)$-stable for any ${\mathbf A} \in {\mathcal M}^{n \times n}$ and any ${\mathbf G} \in {\mathcal G}$.
 \item[\rm 2.] If the matrix class $\mathcal G$ includes the identity element $\mathbf L$ (with respect to the operation $\circ$), then every $({\mathfrak D},{\mathcal G}, \circ)$-stable matrix is ${\mathfrak D}$-stable.
\item[\rm 3.] If $\mathcal G$ forms a subgroup with respect to the operation $\circ$, then ${\mathbf G}_0\circ{\mathbf A}$ is $({\mathfrak D},{\mathcal G},\circ)$-stable (for at least one matrix ${\mathbf G}_0 \in {\mathcal G}$) implies that ${\mathbf A}$ is $\mathfrak D$-stable and $({\mathfrak D},{\mathcal G},\circ)$-stable.
\end{enumerate}
\end{theorem}
{\bf Proof}. \begin{enumerate}
 \item[\rm 1.] Let $\mathbf A$ be $({\mathfrak D},{\mathcal G},\circ)$-stable. Consider ${\mathbf G}\circ{\mathbf A}$, for an arbitrary ${\mathbf G} \in {\mathcal G}$. Then, taking an arbitrary ${\mathbf G}_0 \in {\mathcal G}$, we obtain
$${\mathbf G}_0\circ({\mathbf G}\circ{\mathbf A}) = [\mbox{associativity}] = ({\mathbf G}_0\circ{\mathbf G})\circ{\mathbf A} =$$
$$ = {\mathbf G_1}\circ{\mathbf A},$$
where ${\mathbf G_1}:={\mathbf G}_0\circ{\mathbf G} \in {\mathcal G}$ due to its closeness. Thus $\sigma({\mathbf G}_0\circ({\mathbf G}\circ{\mathbf A})) = \sigma({\mathbf G_1}\circ{\mathbf A}) \subset {\mathfrak D}$ for any ${\mathbf G} \in {\mathcal G}$.
\item[\rm 2.] Let ${\mathbf A}$ be a $({\mathfrak D},{\mathcal G}, \circ)$-stable matrix. Then obviously $\sigma({\mathbf A}) = \sigma({\mathbf L}\circ{\mathbf A}) \subset {\mathfrak D}$.
\item[\rm 3.] Let ${\mathbf G}_0\circ{\mathbf A}$ be $({\mathfrak D},{\mathcal G},\circ)$-stable. Consider $\mathbf A$. Then
$${\mathbf G}\circ{\mathbf A} = {\mathbf G}\circ{\mathbf L}\circ{\mathbf A} = {\mathbf G}\circ({\mathbf G}_0^{-1}{\mathbf G}_0)\circ{\mathbf A} =$$
$$ = [\mbox{associativity}] = ({\mathbf G}\circ{\mathbf G}_0^{-1})\circ({\mathbf G}_0\circ{\mathbf A}) = {\mathbf G}_1\circ({\mathbf G}_0\circ{\mathbf A}), $$ where ${\mathbf G}_1 := {\mathbf G}\circ{\mathbf G}_0^{-1} \in {\mathcal G}$.
Thus $\sigma({\mathbf G}\circ{\mathbf A}) = \sigma({\mathbf G}_1\circ({\mathbf G}_0\circ{\mathbf A})) \subset {\mathfrak D}$.
\end{enumerate}
$\square$

Consider the operation of matrix addition, then the identity element is a zero matrix $\mathbf O$. For matrix multiplication, it is an identity matrix ${\mathbf I}$, for Hadamard multiplication, it is a matrix ${\mathbf E}$ with all the entries $e_{ij} = 1$.

\begin{corollary}(See, for example, \cite{LOG}, p. 79.)
\begin{enumerate}
 \item[\rm 1.] The set of (multiplicative) $D$-stable matrices is invariant under multiplication by a positive diagonal matrix $\mathbf D$.
\item[\rm 2.] Any $D$-stable matrix is stable.
\end{enumerate}
\end{corollary}
As a nontrivial example, we may consider a group of commutators of a given symmetric positive definite matrix ${\mathbf G}$.

\section{Basic properties of $({\mathfrak D}, \ {\mathcal G}, \ \circ)$-stable matrices}
Here, we consider some basic matrix operations, which preserve the class of $({\mathfrak D}, \ {\mathcal G}, \ \circ)$-stable matrices for some specified stability regions ${\mathfrak D}$, matrix classes ${\mathcal G}$ and binary operations $\circ$. We provide a unified way to prove the elementary properties for the partial cases, described in Section 1, mostly focusing on the properties of the binary operation $\circ$.

\subsection{Transposition} First, consider ${\mathbf A}^T$ (the transpose of ${\mathbf A}$).
\begin{theorem}\label{trans} Let ${\mathfrak D}\subset{\mathbb C}$ be an arbitrary (symmetric with respect to the real axes) stability region, ${\mathcal G}\subset {\mathcal M}^{n \times n}$ be an arbitrary matrix class and $\circ$ be a binary operation, satisfying the following property:
\begin{equation}\label{prop1}({\mathbf G}\circ{\mathbf A})^T = {\mathbf A}^T\circ{\mathbf G}^T \end{equation} for any matrix ${\mathbf A} \in {\mathcal M}^{n \times n}$ and any ${\mathbf G} \in {\mathcal G}$.
Then:
\begin{enumerate}
\item[\rm 1.] a matrix $\mathbf A$ is left $({\mathfrak D}, \ {\mathcal G}, \ \circ)$-stable if and only if ${\mathbf A}^T$ is right $({\mathfrak D}, \ {\mathcal G}^T, \ \circ)$-stable.
\item[\rm 2.] If, in addition, $\mathcal G$ is closed with respect to the matrix transposition (i.e. ${\mathbf G} \in {\mathcal G}$ if and only if ${\mathbf G}^T \in {\mathcal G}$) and the equality
\begin{equation}\label{comm}
\sigma({\mathbf G}\circ{\mathbf A}) = \sigma({\mathbf A}\circ{\mathbf G})
\end{equation}
holds,
 then $\mathbf A$ is $({\mathfrak D}, \ {\mathcal G}, \ \circ)$-stable if and only if ${\mathbf A}^T$ is $({\mathfrak D}, \ {\mathcal G}, \ \circ)$-stable.
\end{enumerate}
\end{theorem}
{\bf Proof.}  Let ${\mathbf A}$ be left $({\mathfrak D}, \ {\mathcal G}, \ \circ)$-stable. Consider ${\mathbf A}^T$. Take arbitrary $\widetilde{\mathbf G}$ from the class ${\mathcal G}^T$. Applying Property \eqref{prop1} to $\widetilde{{\mathbf G}}\circ{\mathbf A}^T$, we obtain
$${\mathbf A}^T \circ \widetilde{{\mathbf G}} = (\widetilde{{\mathbf G}}^T\circ{\mathbf A})^T, $$
which implies
$$ \sigma({\mathbf A}^T \circ\widetilde{ {\mathbf G}}) = \sigma(\widetilde{{\mathbf G}}^T\circ{\mathbf A})^T = \sigma(\widetilde{{\mathbf G}}^T\circ{\mathbf A})$$
since the spectra of real-valued matrices are symmetric with respect to the real axes. The inclusion ${\mathbf G}^T \in {\mathcal G}$ implies $\sigma({\mathbf G}^T\circ{\mathbf A}) \subset {\mathfrak D}$.

The second part of Theorem \ref{trans} obviously follows from Property \eqref{comm}. $\square$

The operations of matrix multiplication, matrix addition, Hadamard and block Hadamard matrix multiplication all satisfy Property \eqref{prop1} and Property \eqref{comm}. Thus, the first part of Theorem \ref{trans} holds for multiplicative $({\mathfrak D},{\mathcal G})$-stable, additive $({\mathfrak D},{\mathcal G})$-stable and Hadamard $({\mathfrak D},{\mathcal G})$-stable matrices independently of stability region ${\mathfrak D}$. For all the partial cases, described in Subsections 1.3-1.5, the corresponding matrix classes ${\mathcal G}$ are closed with respect to matrix transposition. Thus the second part of Theorem \ref{trans} holds for all this cases. For many of them, it was pointed out before, see, for example, the following corollaries.
\begin{corollary}
 $\mathbf A$ is multiplicative $D$-stable if and only if ${\mathbf A}^T$ is multiplicative $D$-stable (see \cite{JOHN2}).
 \end{corollary}
 \begin{corollary}
 $\mathbf A$ is multiplicative $H$-stable if and only if ${\mathbf A}^T$ is multiplicative $H$-stable (see \cite{JOHN2}).
 \end{corollary}
\begin{corollary}
 $\mathbf A$ is multiplicative $D(\alpha)$-stable if and only if ${\mathbf A}^T$ is multiplicative $D(\alpha)$-stable (see \cite{BAO}).
\end{corollary}
\begin{corollary}
 $\mathbf A$ is multiplicative $H(\alpha)$-stable if and only if ${\mathbf A}^T$ is multiplicative $H(\alpha)$-stable (see \cite{BAO}).
\end{corollary}
\begin{corollary}
 $\mathbf A$ is ordered $D$-stable if and only if ${\mathbf A}^T$ is ordered $D$-stable (see \cite{JOHN2}).
 \end{corollary}
 \begin{corollary}
 $\mathbf A$ is interval $D$-stable if and only if ${\mathbf A}^T$ is interval $D$-stable (see \cite{JOHN2}).
 \end{corollary}
\begin{corollary}
 $\mathbf A$ is Schur $D$-stable if and only if ${\mathbf A}^T$ is Schur $D$-stable.
\end{corollary}
\begin{corollary}
 $\mathbf A$ is multiplicative $D$-positive ($D$-aperiodic) if and only if ${\mathbf A}^T$ is multiplicative $D$-positive (respectively, $D$-aperiodic) (see \cite{BAO}).
\end{corollary}
\begin{corollary}
 $\mathbf A$ is multiplicative $D$-hyperbolic if and only if ${\mathbf A}^T$ is multiplicative $D$-hyperbolic (see \cite{BAO}).
\end{corollary}
\begin{corollary}
 $\mathbf A$ is additive $D$-stable if and only if ${\mathbf A}^T$ is additive $D$-stable (see \cite{BAO}).
\end{corollary}
\begin{corollary}
 $\mathbf A$ is Hadamard $H$-stable if and only if ${\mathbf A}^T$ is Hadamard $H$-stable (see \cite{JOHD}).
\end{corollary}
\begin{corollary}
 $\mathbf A$ is $B_k$-stable if and only if ${\mathbf A}^T$ is $B_k$-stable (see \cite{JOHD}).
\end{corollary}

\subsection{Inversion} Now, suppose $\circ$ be associative and invertible and consider $(\circ{\mathbf A})^{-1}$ (the inverse of ${\mathbf A}$ with respect to the binary operation $\circ$).

 Assume that there are two stability regions $\mathfrak D, \widetilde{\mathfrak D} \subset {\mathbf C}$ such that
\begin{equation}\label{in}
\sigma(\mathbf A) \subset {\mathfrak D} \ \mbox{implies} \ \sigma((\circ{\mathbf A})^{-1})\subset\widetilde{\mathfrak D}.
 \end{equation}
In particular cases, we know a one-to-one map $\varphi: {\mathbb C} \rightarrow {\mathbb C}$ which connects the spectra of $\mathbf A$ and $(\circ{\mathbf A})^{-1}$:
 $$\sigma((\circ{\mathbf A})^{-1}) = \varphi(\sigma(\mathbf A)).$$
 The following statement holds.
    \begin{theorem}\label{inver}
Let $\circ$ be an associative and invertible matrix operation, ${\mathfrak D}$ and $\widetilde{\mathfrak D} \subset{\mathbb C}$ be two stability regions connected by Property \eqref{in}, and ${\mathcal G}\subset {\mathcal M}^{n \times n}$ be a class of invertible with respect to $\circ$ matrices. Then
 \begin{enumerate}
 \item[\rm 1.] A matrix $\mathbf A$ is left $({\mathfrak D}, \ {\mathcal G}, \ \circ)$-stable implies $(\circ{\mathbf A})^{-1}$ is right $(\widetilde{{\mathfrak D}}, \ {\mathcal G}^{-1}, \ \circ)$-stable.

 \item[\rm 2.] If, in addition, the region $\widetilde{{\mathfrak D}}$ is also connected to the region ${\mathfrak D}$ by Property \eqref{in}, i.e.
$$ \sigma((\circ{\mathbf A})^{-1})\subset\widetilde{\mathfrak D} \ \mbox{implies} \ \sigma(\mathbf A) \subset {\mathfrak D},$$
 matrix class ${\mathcal G}\subset {\mathcal M}^{n \times n}$ is closed with respect to $\circ$-inversion and Property \eqref{comm} holds, then $\mathbf A$ is $({\mathfrak D}, \ {\mathcal G}, \ \circ)$-stable if and only if $(\circ{\mathbf A})^{-1}$ is $({\mathfrak D}, \ {\mathcal G}, \ \circ)$-stable.
 \end{enumerate}
\end{theorem}
{\bf Proof.} For the proof of the first part, let ${\mathbf A}$ be $({\mathfrak D}, \ {\mathcal G}, \ \circ)$-stable. Consider $(\circ{\mathbf A})^{-1}$. Taking arbitrary $(\circ\mathbf G)^{-1}$ from the class ${\mathcal G}^{-1}$, by associativity and invertibility we obtain:
$$(\circ{\mathbf A})^{-1}\circ(\circ\mathbf G)^{-1} = (\circ({\mathbf G}\circ{\mathbf A}^{-1}))^{-1}. $$
Since $\sigma({\mathbf G}\circ{\mathbf A}) \in {\mathfrak D}$, we have $\sigma((\circ{\mathbf A})^{-1}\circ(\circ\mathbf G)^{-1}) = \sigma((\circ({\mathbf G}\circ{\mathbf A}^{-1}))^{-1})\subset \widetilde{{\mathfrak D}}$.

The second part of the Theorem is proved by applying the same reasoning and Property \eqref{comm} to $(\circ{\mathbf A})^{-1}$. $\square$

For matrix addition and matrix multiplication, we know the concrete functions $\varphi_{\circ}(\lambda)=\frac{1}{\lambda}$ and $\varphi_+(\lambda)= -\lambda$. Thus, given a stability region ${\mathfrak D}$, we consider its transformations $-{\mathfrak D}$ and ${\mathfrak D}^{-1}$.

A region ${\mathfrak D}$ is invariant with respect to $\varphi_{\circ}$ if and only if $${\mathfrak D}={\mathfrak D}^{-1} \qquad \mbox{where} \  {\mathfrak D}^{-1}:=\{\lambda \in {\mathbb C}: \frac{1}{\lambda} \in {\mathfrak D}\}.$$ The examples of such regions are:
\begin{enumerate}
\item[-] the unit circle $\{\lambda \in {\mathbb C}: |\lambda|=1\}$;
\item[-] conic regions (both positive and negative directions of the real axes, the imaginary axes, sector regions, left and right half-planes of the complex plane, the complex plane without the imaginary axes, etc), without the origin.
\end{enumerate}

A region ${\mathfrak D}$ is invariant with respect to $\varphi_{+}$ if and only if ${\mathfrak D}= - {\mathfrak D}$. The examples of such regions are:
\begin{enumerate}
\item[-] the unit disk;
\item[-] line-consisting regions (the real and imaginary axes,  the complex plane without the imaginary axes, etc).
\end{enumerate}

Now let us consider an LMI region ${\mathfrak D}$, defined by its characteristic function \eqref{LMI}. Note, that ${\mathfrak D}$ is necessarily convex. It is easy to see that the region $- {\mathfrak D}$ is also an LMI region.
However, for ${\mathfrak D}^{-1}$ it is easy to show that ${\mathfrak D}^{-1}$ is an EMI region.

The matrix classes ${\mathcal G}\subset {\mathcal M}^{n \times n}$ which are closed with respect to taking multiplicative inverses are symmetric and symmetric positive definite matrices, $\alpha$-block diagonal matrices, diagonal matrices and any fixed sign pattern class of diagonal matrices, $\alpha$-scalar matrices and vertex diagonal matrices. The classes which are closed with respect to taking additive inverses are: symmetric, diagonal and vertex diagonal matrices.

Remind, that the property $\sigma({\mathbf A}\circ{\mathbf B}) = \sigma({\mathbf B}\circ{\mathbf A})$ means that left  $({\mathfrak D}, \ {\mathcal G}, \ \circ)$-stability coincides with right  $({\mathfrak D}, \ {\mathcal G}, \ \circ)$-stability, and it holds for multiplication, addition, Hadamard and block Hadamard multiplication. However, here we will not consider Hadamard products because the connection between the spectra of a matrix and its Hadamard inverse is not so trivial.

Now let us consider the known partial cases of $({\mathfrak D}, \ {\mathcal G}, \ \circ)$-stability which satisfy the second part of Theorem, i.e. for which ${\mathfrak D} = {\mathfrak D}^{-1}$ and ${\mathcal G} = {\mathcal G}^{-1}$.
\begin{corollary} $\mathbf A$ is multiplicative $D$-stable if and only if ${\mathbf A}^{-1}$ is multiplicative $D$-stable (see \cite{JOHN2}).
\end{corollary}
\begin{corollary} $\mathbf A$ is multiplicative $H$-stable if and only if ${\mathbf A}^{-1}$ is multiplicative $H$-stable (see \cite{JOHN2}).
\end{corollary}
\begin{corollary} $\mathbf A$ is multiplicative $H(\alpha)$-stable if and only if ${\mathbf A}^{-1}$ is multiplicative $H(\alpha)$-stable (see \cite{JOHN2}).
\end{corollary}
\begin{corollary} $\mathbf A$ is multiplicative $D(\alpha)$-stable if and only if ${\mathbf A}^{-1}$ is multiplicative $D(\alpha)$-stable (see \cite{JOHN2}).
\end{corollary}
\begin{corollary} $\mathbf A$ is multiplicative $D$-positive ($D$-aperiodic) if and only if ${\mathbf A}^{-1}$ is multiplicative $D$-positive (respectively, $D$-aperiodic) (see \cite{BAO}).
\end{corollary}
\begin{corollary} $\mathbf A$ is multiplicative $D$-hyperbolic if and only if ${\mathbf A}^{-1}$ is multiplicative $D$-hyperbolic (see \cite{BAO}).
\end{corollary}

For some cases, we have though ${\mathfrak D} = {\mathfrak D}^{-1}$, the matrix class ${\mathcal G}^{-1}$ do not coincide with ${\mathcal G}$ but can be easily described by matrix inequalities.
\begin{corollary} $\mathbf A$ is ordered $D$-stable with respect to a permutation $\tau \in \Theta$ if and only if ${\mathbf A}^{-1}$ is ordered $D$-stable with respect to $\tau^{-1}$ (see \cite{JOHN2}).
\end{corollary}
\begin{corollary} $\mathbf A$ is interval $D$-stable with respect to a parallelepiped of the form $$ \Theta = \diag\{d_{ii}, \ \ 0 < d_{ii}^{min} < d_{ii} < d_{ii}^{max} < + \infty, \ \  i  = 1, \ \ldots, \ n\}.$$
  if and only if ${\mathbf A}^{-1}$ is interval $D$-stable with respect to the parallelipiped $$ \Theta^{-1} = \diag\{d_{ii}, \ \ 0 < \frac{1}{d_{ii}^{min}} < d_{ii} < \frac{1}{d_{ii}^{max}} < + \infty, \ \  i  = 1, \ \ldots, \ n\}.$$
\end{corollary}

For the case of Schur $D$-stability, we have the following statement:
\begin{corollary} $\mathbf A$ is Schur $D$-stable if and only if ${\mathbf A}^{-1}$ is multiplicative $({\mathbb C}\setminus \overline{D}(0,1), \widetilde{D})$-stable, where ${\mathbb C} \setminus \overline{D}(0,1)$ is the exterior of the closed unit disk, $\widetilde{D}$ is the class of diagonal matrices with $|d_{ii}| > 1$.
\end{corollary}
Note that the class of vertices remains the same under inversion, thus we have

\begin{corollary} $\mathbf A$ is vertex Schur stable if and only if ${\mathbf A}^{-1}$ is vertex ${\mathfrak D}$-stable, where ${\mathfrak D}$ is the exterior of the closed unit disk.
\end{corollary}

The corresponding statements for the case of additive ${\mathcal G}$-stability can be easily obtained.

\subsection{Multiplication by a scalar} Multiplication by a scalar is particularly useful for the transition from an unbounded matrix class ${\mathcal G}$ to some bounded class $\widetilde{{\mathcal G}}$. Given a finite or infinite interval $(\underline{\alpha},\overline{\alpha})$ of the real line, we refer to Properties \eqref{scaas} and \eqref{scadist} (see Section 4) connecting a binary operation $\circ$ to the operation of scalar multiplication. The following statement holds.

\begin{theorem} Let $(\underline{\alpha},\overline{\alpha}) \subseteq {\mathbb R}$ be a finite or infinite interval, ${\mathfrak D}\subset{\mathbb C}$ be a stability region satisfying $\alpha\lambda \in {\mathfrak D}$ for any $\lambda \in {\mathfrak D}$, $\alpha \in (\underline{\alpha},\overline{\alpha})$, ${\mathcal G} \subset {\mathcal M}^{n \times n}$ be an arbitrary matrix class and $\circ$ be a binary operation. If one of the following cases holds
 \begin{enumerate} \item[\rm 1.] the operation $\circ$ is connected to scalar multiplication by Property \eqref{scaas}
 \item[\rm 2.] $\circ$ is connected to scalar multiplication by Property \eqref{scadist} and the matrix class ${\mathcal G} \subset {\mathcal M}^{n \times n}$ satisfies $\frac{1}{\alpha}{\mathbf G} \in {\mathcal G}$ for any ${\mathbf G} \in {\mathcal G}$ and any $\alpha \in (\underline{\alpha},\overline{\alpha})$.
 \end{enumerate}
 then a matrix $\mathbf A$ is $({\mathfrak D}, \ {\mathcal G}, \ \circ)$-stable implies $\alpha{\mathbf A}$ is $({\mathfrak D}, \ {\mathcal G}, \ \circ)$-stable for any $\alpha \in (\underline{\alpha},\overline{\alpha})$.
\end{theorem}
{\bf Proof.}\begin{enumerate}
\item[\rm 1.] Since ${\mathbf G}\circ(\alpha{\mathbf A}) = \alpha({\mathbf G}\circ{\mathbf A})$ and $\alpha{\mathfrak D} \subseteq {\mathfrak D}$, we have
$$\sigma({\mathbf G}\circ(\alpha{\mathbf A})) = \alpha\sigma({\mathbf G}\circ{\mathbf A}) \subset {\mathfrak D}. $$
\item[\rm 2.] Since ${\mathbf G}\circ(\alpha{\mathbf A}) = \alpha((\frac{1}{\alpha}{\mathbf G})\circ{\mathbf A})$ and $\frac{1}{\alpha}{\mathbf G} \in {\mathcal G}$ and $\alpha{\mathfrak D} \subseteq {\mathfrak D}$, we have
$$\sigma({\mathbf G}\circ(\alpha{\mathbf A})) = \alpha\sigma((\frac{1}{\alpha}{\mathbf G})\circ{\mathbf A}) \subset {\mathfrak D}. $$
\end{enumerate} $\square$

Considering $(\underline{\alpha},\overline{\alpha}) = {\mathbb R}$ (except, possibly, zero), we obtain that the stability region ${\mathfrak D}$ consists of lines coming through the origin. Thus if $\mathbf A$ is a $D$-hyperbolic matrix, $\alpha\mathbf A$ is also $D$-hyperbolic for any $\alpha \in {\mathbb R}$. In its turn, considering $(\underline{\alpha},\overline{\alpha}) = (0, +\infty)$ (except, possibly, zero), we obtain that the stability region ${\mathfrak D}$ consists of half-lines coming from the origin. As examples, we may consider positive direction of the real axes, open and closed right (left) halfplane and so on. Thus if $\mathbf A$ is $D$-positive ($D$-stable), $\alpha\mathbf A$ is also $D$-positive (respectively, $D$-stable) for any $\alpha >0$. Considering $(\underline{\alpha},\overline{\alpha}) =(-1,1)$, we obtain, that the multiplication by $\alpha$ maps the unit disk ${\mathfrak D} = \{z \in {\mathbb C}: |z|< 1\}$ into itself. Thus if $\mathbf A$ is Schur $D$-stable, $\alpha\mathbf A$ is also Schur $D$-stable for any $\alpha \in (-1, 1)$.

\subsection{Similarity transformations} Here, we consider the following question: given a nonsingular matrix $\mathbf S$ and a $({\mathfrak D},{\mathcal G}, \circ)$-stable matrix $\mathbf A$, when ${\mathbf S}{\mathbf A}{\mathbf S}^{-1}$ is again $({\mathfrak D},{\mathcal G}, \circ)$-stable? I.e. we describe the class of similarity transformation that preserve $({\mathfrak D},{\mathcal G}, \circ)$-stability.

\begin{theorem}
Given a class of nonsingular matrices ${\mathcal S} \subset {\mathcal M}^{n \times n}$ (closed with respect to multiplicative inversion), a stability region ${\mathfrak D}\subset{\mathbb C}$, a matrix class ${\mathcal G}\subset {\mathcal M}^{n \times n}$ and a binary matrix operation $\circ$. If the matrix class ${\mathcal S}$ commutes with the  matrix class ${\mathcal G}$ and one of the following cases holds
\begin{enumerate}
\item[\rm 1.] The operation $\circ$ is matrix multiplication;
\item[\rm 2.] The operation $\circ$ is connected to matrix multiplication by Property \eqref{*dist2}
\end{enumerate}
then a matrix ${\mathbf S}{\mathbf A}{\mathbf S}^{-1}$ is $({\mathfrak D},{\mathcal G}, \circ)$-stable if and only if ${\mathbf A}$ is $({\mathfrak D},{\mathcal G}, \circ)$-stable.
\end{theorem}
{\bf Proof.} Case 1. Let ${\mathbf A}$ be multiplicative $({\mathfrak D},{\mathcal G})$-stable. Consider ${\mathbf S}{\mathbf A}{\mathbf S}^{-1}$. Since an arbitrary ${\mathbf G} \in {\mathcal G}$ commutes with an arbitrary $\mathbf S \in {\mathcal S}$, we have
$${\mathbf G}{\mathbf S}{\mathbf A}{\mathbf S}^{-1} = {\mathbf S}{\mathbf G}{\mathbf A}{\mathbf S}^{-1} = {\mathbf S}({\mathbf G}{\mathbf A}){\mathbf S}^{-1}. $$
Since $$\sigma({\mathbf S}({\mathbf G}{\mathbf A}){\mathbf S}^{-1}) = \sigma({\mathbf G}{\mathbf A}) \subset {\mathfrak D}$$
we obtain that ${\mathbf S}{\mathbf A}{\mathbf S}^{-1}$ is also multiplicative $({\mathfrak D},{\mathcal G})$-stable. For the inverse direction is enough to notice that if ${\mathbf B} := {\mathbf S}{\mathbf A}{\mathbf S}^{-1}$ then ${\mathbf A} = {\mathbf S}^{-1}{\mathbf B}{\mathbf S}$ and $\mathbf G$ commutes with $\mathbf S$ if and only if $\mathbf G$ commutes with ${\mathbf S}^{-1}$.

Case 2. Commutativity between ${\mathcal G}$ and ${\mathcal S}$ implies that the matrix class ${\mathcal G}$ is invariant with respect to the linear transformations from ${\mathcal S}$. Let ${\mathbf A}$ be $({\mathfrak D},{\mathcal G},\circ)$-stable. Consider ${\mathbf S}{\mathbf A}{\mathbf S}^{-1}$. Applying Property \eqref{*dist2}, we obtain
$${\mathbf G}\circ{\mathbf S}{\mathbf A}{\mathbf S}^{-1} = ({\mathbf S}{\mathbf S}^{-1}{\mathbf G})\circ({\mathbf S}{\mathbf A}{\mathbf S}^{-1})= {\mathbf S}({\mathbf S}^{-1}{\mathbf G}{\mathbf S}{\mathbf S}^{-1})\circ({\mathbf A}{\mathbf S}^{-1})= $$
$$= {\mathbf S}({\mathbf S}^{-1}{\mathbf G}{\mathbf S}\circ{\mathbf A}){\mathbf S}^{-1} ={\mathbf S}(\widetilde{{\mathbf G}}\circ{\mathbf A}){\mathbf S}^{-1},$$
where $\widetilde{{\mathbf G}} = {\mathbf S}^{-1}{\mathbf G}{\mathbf S} \in {\mathcal G}$. Since $$\sigma({\mathbf S}(\widetilde{{\mathbf G}}\circ{\mathbf A}){\mathbf S}^{-1}) = \sigma(\widetilde{{\mathbf G}}\circ{\mathbf A}) \subset {\mathfrak D},$$
we obtain that ${\mathbf S}{\mathbf A}{\mathbf S}^{-1}$ is $({\mathfrak D},{\mathcal G},\circ)$-stable. The proof for the inverse direction copies the same reasoning.
$\square$

Here, let us consider several matrix classes and their commutators. As it is known, the class ${\mathcal G}$ of diagonal matrices commutes with itself and with the class of permutation matrices. Thus we obtain the following statement (see, for example, \cite{BAO}, p. 68 for the case of $D$-positive matrices, \cite{AM}, p. 450, Theorem 2 and \cite{JOHN1}, p. 54, Observation (ii), for the case of $D$-stable matrices).
\begin{corollary} Let $\mathbf A$ belong to one of the following classes: $D$-stable matrices, $D$-positive matrices, Schur $D$-stable matrices or $D$-hyperbolic matrices. Then the matrices ${\mathbf D}{\mathbf A}{\mathbf D}^{-1}$, where ${\mathbf D}$ is a diagonal matrix and ${\mathbf P}{\mathbf A}{\mathbf P}^{-1}$, where ${\mathbf P}$ is a permutation matrix, also belong to the same class.
\end{corollary}

\subsection{Example. Symmetrized $({\mathcal G})$-negativity} Now we consider some well-known matrix problems from the point of view of $({\mathfrak D},{\mathcal G},\circ)$-stability. Let the stability region ${\mathfrak D}$ be the negative direction of the real axes (with zero excluded). Let the binary operation $\circ$ be the {\it symmetrized matrix product} (see \cite{BHAT2}), defined as follows:
 $${\mathbf A}{\circ}_S {\mathbf B}:={\mathbf B}{\mathbf A} + {\mathbf A}^T{\mathbf B}^T, \qquad {\mathbf A}, {\mathbf B} \in {\mathcal M}^{n \times n}. $$

 {\bf Simultaneous stability.} The problem of simultaneous stability is rather old (for the beginning, see \cite{TAU}, \cite{BAR}). Recall, that two matrices ${\mathbf A}_1$, ${\mathbf A}_2$ are called {\it simultaneously stable} if they admit a common positive definite solution of the Lyapunov equation \eqref{lyap}.

  Assume a fixed matrix ${\mathbf B}$ be symmetric positive definite. Consider the Lyapunov equation \eqref{lyap}:
 $$ {\mathbf A}{\circ}_S {\mathbf B}:={\mathbf B}{\mathbf A} + {\mathbf A}^T{\mathbf B} = {\mathbf W},$$
where ${\mathbf W}$ is a symmetric negative definite matrix.
Using an equivalent characterization of negative definiteness, we obtain that ${\mathbf W}$ is negative definite if and only if $\sigma({\mathbf W}) \subset {\mathbb R}^-$. Thus the problem of describing the matrix class ${\mathcal G} \subset {\mathcal M}^{n \times n}$ such that
$$ \sigma({\mathbf G}{\circ}_S {\mathbf B}) = \sigma({\mathbf W})\subset {\mathbb R}^-,$$
for every ${\mathbf G} \in {\mathcal G}$ is equivalent to establishing simultaneous stability of all matrices from ${\mathcal G}$, with the common Lyapunov factor ${\mathbf B}$.

{\bf Preserving definiteness.} Now assume that ${\mathbf A}$ is fixed. Consider the following type of problems, which often arises in studying robust stability. Describe the subclass ${\mathcal G}$ of symmetric positive definite matrices such that for each ${\mathbf G} \in {\mathcal G}$
$$ {\mathbf A}{\circ}_S {\mathbf G}:={\mathbf G}{\mathbf A} + {\mathbf A}^T{\mathbf G} = {\mathbf W},$$
where ${\mathbf W}$ is a symmetric negative definite matrix.
I.e. $$ \sigma({\mathbf A}{\circ}_S {\mathbf G}) \subset {\mathbb R}^-,$$
for every ${\mathbf G} \in {\mathcal G}$.

\section{Generalized diagonal stability and sufficient conditions of $({\mathfrak D},{\mathcal G},\circ)$-stability}
In this section, we generalize the widely studied concept of diagonal stability for the case of different stability regions $\mathfrak D$.

\subsection{Volterra--Lyapunov $({\mathfrak D},{\mathcal P})$-stability}
Here, we consider a stability region ${\mathfrak D}$, defined by generalized Lyapunov equation \eqref{GenLyap} and provide the following definition.

Given a stability region ${\mathfrak D}$, defined by Equation \eqref{GenLyap} and a subclass $\mathcal P$ of the class of symmetric positive definite matrices $\mathcal H$, we call a matrix $\mathbf A$ {\it Volterra--Lyapunov $({\mathfrak D},{\mathcal P})$-stable}, if Equation \eqref{GenLyap} admits a solution in the matrix class $\mathcal P$, i.e. if there exists a matrix ${\mathbf P} \in {\mathcal P}$ such that
$$ {\mathbf W}:=\sum_{i,j = 0}^{n-1}c_{ij}({\mathbf A}^T)^i{\mathbf P}{\mathbf A}^j$$
is positive definite.

As partial cases, we mention the following matrix classes.

\begin{enumerate}
\item[\rm 1.] Remind, that an $n \times n$ real matrix $\mathbf A$ is called {\it diagonally stable} if there exists a positive diagonal matrix $\mathbf D$ such that ${\mathbf D}{\mathbf A} + {\mathbf A}^T{\mathbf D}$ is positive definite. In this case, the matrix ${\mathbf D}$ is called a {\it Lyapunov scaling factor} of ${\mathbf A}$. The concept of diagonal stability arises in \cite{QR}, referring \cite{AM} as a characterization of multiplicative $D$-stability. The property of diagonal stability is studied in \cite{CROSS} as {\it Volterra--Lyapunov stability} and in \cite{LOG} as {\it dissipativity}. For other references and tittles of this property, see \cite{LOG}, p. 82. Here, the stability region is the left hand side of the complex plane, described by classical Lyapunov equation \eqref{lyap}, the matrix class ${\mathcal P}$ is the class of positive diagonal matrices.
\item[\rm 2.] An $n \times n$ real matrix $\mathbf A$ is called {\it $\alpha$-scalarly stable} if there exists a positive $\alpha$-scalar matrix ${\mathbf D}_{\alpha}$ such that ${\mathbf D}_{\alpha}{\mathbf A} + {\mathbf A}^T{\mathbf D}_{\alpha}$ is positive definite. This matrix class was introduced in \cite{HERM} under the name of Lyapunov $\alpha$-scalar stability and then studied in \cite{GUH1}, \cite{WAN}. Here, the stability region is again the left hand side of the complex plane, the matrix class ${\mathcal P}$ is the class ${\mathcal D}_{\alpha}$ of $\alpha$-scalar matrices.
\item[\rm 3.] An $n \times n$ real matrix $\mathbf A$ is called {\it $\alpha$-diagonally stable} if there exists a symmetric positive definite $\alpha$-diagonal matrix ${\mathbf H}_{\alpha}$ such that ${\mathbf H}_{\alpha}{\mathbf A} + {\mathbf A}^T{\mathbf H}_{\alpha}$ is positive definite. This matrix class was mentioned in \cite{KHAK2} in connection with the study of $D_{\alpha}$-stable matrices and studied in \cite{AB1} in connection with robust stability properties. For the applications, see also \cite{KHAL2}. Here, the stability region is again the left hand side of the complex plane, the matrix class ${\mathcal P}$ is the class ${\mathcal H}_{\alpha}$ of symmetric $\alpha$-diagonal positive definite matrices.
 \item[\rm 4.]  An $n \times n$ real (not necessarily symmetric) matrix $\mathbf A$ is called {\it positive definite} if its symmetric part ${\mathbf A} + {\mathbf A}^T$ is positive definite. This matrix class was introduced in \cite{JOHN6} as a generalization of positive definiteness to non-symmetric matrices. As an equivalent characterization, it was stated that ${\mathbf A} + {\mathbf A}^T$ is positive definite if and only if $x^T{\mathbf A}x > 0$ for every nonzero vector $x \in {\mathbb R}^n$. For such matrices, the term $L$-stability is also used (see \cite{GUH1}). Here, the stability region is again the left hand side of the complex plane, the matrix class ${\mathcal P}$ consists of the only one identity matrix $\mathbf I$.
 \item[\rm 5.] An $n \times n$ real matrix $\mathbf A$ is called {\it Schur diagonally stable} if there exists a positive diagonal matrix ${\mathbf D}$ such that ${\mathbf D} - {\mathbf A}^T{\mathbf D}{\mathbf A}$ is positive definite. This definition was given in \cite{BHK}. Here, the stability region is the unit disk, defined by the Stein equation \eqref{ste}, the matrix class ${\mathcal P}$ is the class of positive diagonal matrices.
\end{enumerate}

In connection with the definition of generalized Volterra--Lyapunov stability, the following crucial question arises:

{\bf Problem 5.} Given a Lyapunov stability region $\mathfrak D$, two matrix classes ${\mathcal P} \subset {\mathcal H}$ and ${\mathcal G} \subset {\mathcal M}^{n \times n}$, and a binary operation $\circ$ on ${\mathcal M}^{n \times n}$, how the class of Volterra--Lyapunov $({\mathfrak D},{\mathcal P})$-stable matrices is connected to the class of $({\mathfrak D}, \ {\mathcal G}, \ \circ)$-stable matrices?

Taking the matrix class ${\mathcal P}$ to be the class of positive diagonal matrices, we obtain the following generalizations of diagonal stability: given a stability region ${\mathfrak D}$, defined by Equation \eqref{GenLyap}, we call a matrix $\mathbf A$ {\it diagonally ${\mathfrak D}$-stable} (with respect to a polynomial region $\mathfrak D$), if Equation \eqref{GenLyap} admits a solution in the class ${\mathcal D}^+$, i.e. if there exists a positive diagonal matrix ${\mathbf D}$ such that
$$ {\mathbf W}:=\sum_{i,j = 0}^{n-1}c_{ij}({\mathbf A}^T)^i{\mathbf D}{\mathbf A}^j$$
is positive definite. The obtained concept of $\mathfrak D$-diagonal stability includes the known concepts of (Lyapunov) diagonal stability and Schur diagonal stability.

Using different generalizations of Lyapunov theorem stated in Subsection 2.2, we get the following definitions of diagonal ${\mathfrak D}$-stability for different kinds of regions.

\begin{enumerate}
\item[\rm 1.] {\bf Diagonal hyperbolicity.} A matrix $\mathbf A$ is called {\it diagonally hyperbolic} if the Lyapunov equation \ref{lyap} admits a diagonal solution, i.e if there is a diagonal matrix $\mathbf D$ such that the matrix ${\mathbf D}{\mathbf A} + {\mathbf A}^{T}{\mathbf D}$ is positive definite. This generalization is based on the criterion of hyperbolicity by Ostrowski and Schneider (see Theorem \ref{OS}).
\item[\rm 2.] {\bf Diagonal stability for LMI regions.} Given an LMI region $\mathfrak D$, defined by \eqref{LMI}, a matrix $\mathbf A$ is called {\it diagonally $\mathfrak D$-stable} (with respect to an LMI region $\mathfrak D$) if the generalized Lyapunov equation \eqref{LMIeq} admits a positive diagonal solution, i.e. if there is a positive diagonal matrix $\mathbf D$ such that the matrix
   $${\mathbf W} :=\ {\mathbf L}\otimes {\mathbf H} + {\mathbf M}\otimes({\mathbf H}{\mathbf A})+{\mathbf M}^T\otimes({\mathbf A}^T{\mathbf H})$$
    is negative definite.
\item[\rm 3.]{\bf Diagonal stability for EMI regions.} Given an EMI region $\mathfrak D$, defined by \eqref{EMI}, a matrix $\mathbf A$ is called {\it diagonally $\mathfrak D$-stable} (with respect to an EMI region $\mathfrak D$) if the generalized Lyapunov equation \eqref{EMIeq} admits a positive diagonal solution, i.e. if there is a positive diagonal matrix $\mathbf D$ such that the matrix
    $${\mathbf W} :=\ {\mathbf R}_{11}\otimes {\mathbf H} + {\mathbf R}_{12}\otimes({\mathbf H}{\mathbf A})+{\mathbf R}_{12}^T\otimes({\mathbf A}^T{\mathbf H})+{\mathbf R}_{22}\otimes({\mathbf A}^T{\mathbf H}{\mathbf A})$$
    is negative definite.
\end{enumerate}

Other cases, based, for example, on the results of Gutman and Jury (see \cite{GUJU}) \cite{GUT}, \cite{GUT2}) can also be considered.

{\bf Problem 6.} To describe the classes of ${\mathfrak D}$-diagonally stable matrices and Volterra-Lyapunov $({\mathfrak D}, {\mathcal P})$-stable matrices, generalizing the results which describe the classes of diagonally stable matrices (see Section 3).

For the study of diagonal stability for LMI regions, we refer to the results in \cite{KU2}.

\subsection{Solutions of the Lyapunov equation and $({\mathbb C}^+,{\mathcal G})$-stability}
Here, we consider the most simple cases. For the convenience of the proof, we take ${\mathfrak D} = {\mathbb C}^+$ and study {\it positive stability} of matrices. The classical stability case is studied by analogy.

\begin{theorem}\label{Commt} Let ${\mathfrak D} = {\mathbb C}^+$, ${\mathcal P}, {\mathcal G} \subseteq {\mathcal H}$ be two commuting subclasses of symmetric positive definite matrices, and $\circ$ be matrix multiplication or matrix addition. Then an $n \times n$ matrix ${\mathbf A}$ is both multiplicative and additive $({\mathfrak D}, {\mathcal G})$-stable if there exist a matrix ${\mathbf P}\in {\mathcal P}$ such that
\begin{equation}\label{lyap1} {\mathbf W} : ={\mathbf P}{\mathbf A} + {\mathbf A}^T{\mathbf P} \end{equation}
       is positive definite.
\end{theorem}
{\bf Proof.} First, let us consider the case of matrix multiplication. Let ${\mathbf W} : ={\mathbf P}{\mathbf A} + {\mathbf A}^T{\mathbf P}$ be symmetric positive definite for some ${\mathbf P}\in {\mathcal P}$. Then, multiplying Equality \eqref{lyap1} from the both sides on arbitrary ${\mathbf G} \in {\mathcal G}$, we obtain $$ {\mathbf G}{\mathbf W}{\mathbf G} : ={\mathbf G}{\mathbf P}{\mathbf A}{\mathbf G} + {\mathbf G}{\mathbf A}^T{\mathbf P}{\mathbf G}$$
     $$ {\mathbf G}{\mathbf W}{\mathbf G} : =({\mathbf G}{\mathbf P})({\mathbf A}{\mathbf G}) + ({\mathbf A}{\mathbf G})^T({\mathbf G}{\mathbf P}).$$
     From the properties of positive definite matrices (see Section 5, also see \cite{BHAT2}) we obtain that ${\mathbf G}{\mathbf W}{\mathbf G}$ and ${\mathbf G}{\mathbf P}$ are both symmetric positive definite. Thus ${\mathbf A}{\mathbf G}$ is (positive) stable by Lyapunov theorem.

Now let $\circ$ be the operation of matrix addition. Again, let ${\mathbf W} : ={\mathbf P}{\mathbf A} + {\mathbf A}^T{\mathbf P}$ be symmetric positive definite for some ${\mathbf P}\in {\mathcal P}$. Consider $\widetilde{{\mathbf W}} : ={\mathbf P}({\mathbf A}+{\mathbf G}) + ({\mathbf A} + {\mathbf G})^T{\mathbf P}$. Then
$$\widetilde{{\mathbf W}} : ={\mathbf P}{\mathbf A}+{\mathbf P}{\mathbf G} + {\mathbf A}^T{\mathbf P} + {\mathbf G}{\mathbf P}=$$
$${\mathbf W} + ({\mathbf P}{\mathbf G}+{\mathbf G}{\mathbf P}).$$
It follows from the commutativity of positive definite classes ${\mathcal P}$ and ${\mathcal G}$ that ${\mathbf P}{\mathbf G}+{\mathbf G}{\mathbf P}$ is also symmetric positive definite (see Section 5). Thus $\widetilde{{\mathbf W}}$ is symmetric positive definite as a sum of positive definite matrices.
      $\square$
\begin{corollary} Diagonally stable matrices are multiplicative $D$-stable (see \cite{AM}).
\end{corollary}
\begin{corollary} Positive definite (not necessarily symmetric) matrices are $H$-stable (see \cite{AM}, p. 449, Theorem 1, also \cite{OSS}, p. 82).
\end{corollary}
\begin{corollary} $\alpha$-scalar diagonally stable matrices are $H_{\alpha}$-stable (see \cite{HERM}, p. 45, Theorem 4.4).
\end{corollary}
\begin{corollary} $\alpha$-block diagonally stable matrices are $D_{\alpha}$-stable (see \cite{KHAK2}, also \cite{AB1}).
\end{corollary}
      For the class of positive diagonal matrices, the existence of positive diagonal solution of Lyapunov equation \eqref{lyap} is sufficient, but not necessary for $D$-stability (see, for example, \cite{JOHN1}). Johnson pointed the Lyapunov diagonal stability as the oldest sufficient condition for $D$-stability (referring \cite{QR}).  In the case of $H$-stable matrices, Ostrowski and Schneider in \cite{OSS} proved the following sufficient condition for $H$-stability: {\it an $n \times n$ matrix ${\mathbf A}$ is $H$-stable if ${\mathbf A}$ + ${\mathbf A}^T$ is positive definite} (see \cite{OSS}, p. 82). Considering also the case of positive semidefinite and singular matrix ${\mathbf A}$ + ${\mathbf A}^T$, they provide the complete characterization of $H$-stable matrices (see \cite{OSS}, p. 82, Theorem 4, also p. 81 Theorem 3 for $H$-semistability), which shows {\it the proper inclusion of the class of positive definite matrices to the class of $H$-stable matrices}. Analogically, {\it Schur diagonally stable matrices form a proper subclass in the class of Schur $D$-stable matrices} (see \cite{KAB}).

 The following result can be easily deduced from Theorem \ref{Commt}.
      \begin{theorem}
      Let an $n \times n$ matrix $\mathbf A$ be (positive) stable. Then it is multiplicative and additive $({\mathbb C}^+,{\mathcal G})$-stable, where the matrix class ${\mathcal G}$ is a subclass of all symmetric positive definite matrices which commute to the symmetric positive definite solution $\mathbf P$ of the Lyapunov equation for the matrix $\mathbf A$.
      \end{theorem}
 {\bf Proof}. For the proof, it is enough to put in the statement of Theorem \ref{Commt} the commuting classes ${\mathcal P} = \{{\mathbf P}\}$ and ${\mathcal G} = \{{\mathbf G} \in {\mathcal H}: {\mathbf G}{\mathbf P} = {\mathbf P}{\mathbf G}\}$.
      
\subsection{Solutions of the Lyapunov equation and ${\mathcal G}$-hyperbolicity}
The following statements are based on Ostrowski and Schneider result (see Theorem \ref{OS}).

Given a subclass $\mathcal G$ of the class of nonsingular symmetric matrices, we call an $n \times n$ matrix $\mathbf A$ {\it multiplicative ${\mathcal G}$-hyperbolic} if all the eigenvalues of ${\mathbf G}{\mathbf A}$ have nonzero real parts for every $n\times n$ matrix ${\mathbf G} \in {\mathcal G}$. This definition generalizes the concept of $D$-hyperbolicity (see Subsection 1.3).

We call an $n \times n $ matrix $\mathbf A$ {\it diagonally hyperbolic} if Lyapunov equation \eqref{lyap} admits a nonsingular diagonal solution, i.e. there is a diagonal matrix $\mathbf D$ such that ${\mathbf D}{\mathbf A} + {\mathbf A}^T{\mathbf D}$ is negative definite.

Unlike the case of stability, we consider the cases of multiplicative and additive ${\mathcal G}$-hyperbolicity separately.

\begin{theorem}\label{Commt3} Let ${\mathfrak D} = \{\lambda \in {\mathbb C}: {\rm Re}(\lambda)\neq 0\}$, ${\mathcal P}, {\mathcal G} \subseteq {\mathcal H}$ be two commuting subclasses of nonsingular symmetric matrices. Then an $n \times n$ matrix ${\mathbf A}$ is  multiplicative ${\mathcal G}$-hyperbolic if there exist a matrix ${\mathbf P}\in {\mathcal P}$ such that
\begin{equation} {\mathbf W} : ={\mathbf P}{\mathbf A} + {\mathbf A}^T{\mathbf P} \end{equation}
       is positive definite.
\end{theorem}
{\bf Proof.}  Let ${\mathbf W} : ={\mathbf P}{\mathbf A} + {\mathbf A}^T{\mathbf P}$ be symmetric positive definite for some ${\mathbf P}\in {\mathcal P}$. Then, multiplying Equality \eqref{lyap1} from the both sides on arbitrary ${\mathbf G} \in {\mathcal G}$, we obtain $$ {\mathbf G}{\mathbf W}{\mathbf G} : ={\mathbf G}{\mathbf P}{\mathbf A}{\mathbf G} + {\mathbf G}{\mathbf A}^T{\mathbf P}{\mathbf G}$$
     $$ {\mathbf G}{\mathbf W}{\mathbf G} : =({\mathbf G}{\mathbf P})({\mathbf A}{\mathbf G}) + ({\mathbf A}{\mathbf G})^T({\mathbf G}{\mathbf P})$$
     By Lemma \ref{COMMU} (see Section 5) and the properties of positive definite matrices, we obtain that ${\mathbf G}{\mathbf W}{\mathbf G}$ is symmetric positive definite and ${\mathbf G}{\mathbf P}$ is also symmetric. Thus ${\mathbf A}{\mathbf G}$ is stable by Ostrowski--Schneider theorem.
      $\square$
\begin{corollary} Diagonally hyperbolic matrices are multiplicative $D$-hyperbolic.
\end{corollary}

The following result follows from the above reasoning.

\begin{theorem}
      Let an $n \times n$ matrix $\mathbf A$ have no pure imaginary eigenvalues (i.e. with zero real parts). Then it is multiplicative ${\mathcal G}$-hyperbolic, where the matrix class ${\mathcal G}$ is a subclass of all nonsingular symmetric matrices which commute to the nonsingular symmetric solution $\mathbf P$ of the Lyapunov equation for the matrix $\mathbf A$.
      \end{theorem}

Now let us consider the additive ${\mathcal G}$-hyperbolicity. Recall, that a {\it sign pattern} ${\rm Sign}({\mathbf D})$ of a diagonal matrix $\mathbf D$ is defined as follows: $${\rm Sign}({\mathbf D}) := {\rm diag}\{{\rm sign}(d_{11}), \ \ldots, \ {\rm sign}(d_{nn})\}.$$
Two diagonal matrices ${\mathbf D}_1$ and ${\mathbf D}_2$ are said to belong to the same sign pattern class if ${\rm Sign}({\mathbf D}_1) = {\rm Sign}({\mathbf D}_2)$. For a given sign pattern $S$, ${\mathcal D}_S$ denotes a sign pattern class of diagonal matrices (see Section 5).

\begin{theorem}\label{Commt4} Let ${\mathfrak D} = \{\lambda \in {\mathbb C}: {\rm Re}(\lambda)\neq 0\}$, ${\mathcal D}_S$ be a sign pattern class of nonsingular diagonal matrices. Then an $n \times n$ matrix ${\mathbf A}$ is additive ${\mathcal D}_S$-hyperbolic if there exist a matrix ${\mathbf D}_0\in {\mathcal D}_S$ such that
\begin{equation} {\mathbf W} : ={\mathbf D}_0{\mathbf A} + {\mathbf A}^T{\mathbf D}_0 \end{equation}
       is positive definite.
\end{theorem}
{\bf Proof.} let ${\mathbf W} : ={\mathbf D}_0{\mathbf A} + {\mathbf A}^T{\mathbf D}_0$ be symmetric positive definite for some ${\mathbf D}_0\in {\mathcal D}_S$. Consider $\widetilde{{\mathbf W}} : ={\mathbf D}_0({\mathbf A}+{\mathbf D}) + ({\mathbf A} + {\mathbf D})^T{\mathbf D}_0$. Then
$$\widetilde{{\mathbf W}} : ={\mathbf D}_0{\mathbf A}+{\mathbf D}_0{\mathbf D} + {\mathbf A}^T{\mathbf D}_0 + {\mathbf D}{\mathbf D}_0=$$
$${\mathbf W} + ({\mathbf D}_0{\mathbf D}+{\mathbf D}{\mathbf D}_0).$$
Since all diagonal matrices commute and ${\mathbf D}$, $\mathbf D_0$ belong to the same sign pattern class, we obtain ${\mathbf D}_0{\mathbf D}+{\mathbf D}{\mathbf D}_0$ is also symmetric positive definite. Thus $\widetilde{{\mathbf W}}$ is symmetric positive definite as a sum of positive definite matrices.
      $\square$
\section{Methods of study}
As we can see through the review in Part I, the most important method of studying multiplicative and additive stability is the analysis of Lyapunov equation \eqref{lyap} and the concept of diagonal stability which is of independent interest due to a lot of applications. The generalization of this concept and its application to the concept of $({\mathfrak D}, {\mathcal G}, \circ)$-stability, we considered in the previous section. Now we make a brief analysis of another approach from the point of view of the possible applications to different kinds of $({\mathfrak D}, {\mathcal G}, \circ)$-stability, unifying and showing the perspectives.
\subsection{Qualitative approach}
The main idea of the {\it generalized qualitative stability} is as follows. Let us introduce the following partition of the real line ${\mathbb R}$:
$${\mathbb R} = \bigcup_{i=1}^7{{\mathbb R}_i}, $$
where $$ {\mathbb R}_1 := (-\infty; -1); \qquad {\mathbb R}_2 := \{-1\};$$
$$ {\mathbb R}_3 := (-1, 0); \qquad {\mathbb R}_4 := \{0\};$$
$$ {\mathbb R}_5 := (0,1); \qquad {\mathbb R}_6 := \{1\}; \qquad {\mathbb R}_7 := (1, + \infty)$$
In some partial cases, the number of the sets can be reduced. The choice of the partition is motivated by the convenience of the description of the set products $R_iR_j$ which arises in the study of the products of matrices.

Given two matrices ${\mathbf A}, {\mathbf B} \in {\mathcal M}^{n \times n}$, we say that ${\mathbf A}$ is {\it $m$-sign equivalent} to ${\mathbf B}$ if for every pair of indices $i,j$, $a_{ij}$ and $b_{ij}$ belongs to the same class ${\mathbb R}_k$, $k = 1, \ \ldots, \ 7$. Now we define {\it $m$-sign pattern class} as the set of all matrices from ${\mathcal M}^{n \times n}$, $m$-sign-equivalent to a given one. The $m$-sign pattern class, generated by ${\mathbf A}$, we will denote as $m({\mathbf A})$. Given a stability region ${\mathfrak D}$, we say that a given $m$-sign pattern class {\it requires ${\mathfrak D}$-stability}, if all matrices from this class are ${\mathfrak D}$-stable and {\it allows ${\mathfrak D}$-stability}, if at least one matrix from this class is ${\mathfrak D}$-stable. Given a matrix class ${\mathcal G}$ (described by its $m$-sign pattern), and the binary operation $\circ$ (assumed to be matrix multiplication, Hadamard or block Hadamard matrix multiplication) we say that a given $m$-sign pattern class {\it requires $({\mathfrak D}, {\mathcal G}, \circ)$-stability}, if all matrices from this class are $({\mathfrak D}, {\mathcal G}, \circ)$-stable and {\it allows $({\mathfrak D}, {\mathcal G}, \circ)$-stability}, if at least one matrix from this class is $({\mathfrak D}, {\mathcal G}, \circ)$-stable.

{\bf Example.} Besides qualitative stability (see the paper by Quirk and Ruppert \cite{QR}), which was shown to be a sufficient condition for $D$-stability, the following construction for the case of Schur stability was introduced in \cite{KAB1} and studied in \cite{PRY1}, \cite{PRY2}. Given two matrices ${\mathbf A} = \{a_{ij}\}_{i,j =1}^n$ and ${\mathbf B}=\{b_{ij}\}_{i,j =1}^n$, they are called {\it modulus equivalent} if each $a_{ij}$ and $b_{ij}$ belongs to the same set $C_1 = R_3\cup R_5$, $C_2 = R_1 \cup R_7$, $C_3 = R_2 \cup R_6$ or $C_4 = R_4$. The set of matrices, modulus equivalent to a given one is called a {\it modulus pattern class}. A matrix ${\mathbf A}$ is called {\it qualitative Schur stable}, if all the matrices from the modulus pattern class $m(A)$ are Schur stable. The set of qualitative Schur stable matrices is rather small, it consists of diagonal Schur stable matrices and their permutations, and, as it is easy to see, every qualitative Schur stable matrix is Schur $D$-stable and moreover, Schur diagonally stable (see \cite{KAB}, p. 74).

This approach is potentially useful for the applications, since it does not require exact knowledge of the matrix entries, but only the localization in some prescribed intervals. In literature, it is mostly studied by graph-theoretic methods.

Another possible approach to qualitative stability generalization will be as follows. Given a stability region ${\mathfrak D}$, assume that ${\mathfrak D}\cap{\mathbb R} \neq \emptyset$. Then we introduce the following partition of the real line ${\mathbb R}$:
$${\mathbb R} = \bigcup_{i=1}^4{{\mathbb R}_i}, $$
where $$ {\mathbb R}_1 := {\mathbb R}\cap {\mathfrak D}; \qquad {\mathbb R}_2 := {\mathbb R}\setminus{\mathfrak D};$$
$$ {\mathbb R}_3 := {\mathbb R}\cap \partial({\mathfrak D}); \qquad {\mathbb R}_4 := \{0\}.$$

This kind of partition would allow us to describe qualitatively ${\mathfrak D}$-stable matrices using Gershgorin theorem for special types of regions.

Finally, the concept of the "sign pattern set", where we consider the location not only of the entries of $\mathbf A$ but also of its minors, would be of interest.
Given a $n \times n$ matrix $\mathbf A$, it determines a sequence of sign patterns $\{{\mathcal A}, \ {\mathcal A}^{(2)}, \ \ldots, \ {\mathcal A}^{(n)}\}$, where ${\mathcal A}^{(j)}$ is the sign pattern of the $j$th compound matrix ${\mathbf A}^{(j)}$. We say that a set of sign patterns $\{{\mathcal A}, \ {\mathcal A}^{(2)}, \ \ldots, \ {\mathcal A}^{(n)}\}$ allows ${\mathfrak D}$-stability (or is potentially ${\mathfrak D}$-stable) if there exists at least one ${\mathfrak D}$-stable matrix which set of sign patterns coincides with  $\{{\mathcal A}, \ {\mathcal A}^{(2)}, \ \ldots, \ {\mathcal A}^{(n)}\}$. We say that a set of sign patterns $\{{\mathcal A}, \ {\mathcal A}^{(2)}, \ \ldots, \ {\mathcal A}^{(n)}\}$ requires ${\mathfrak D}$-stability if all matrices with such set of sign patterns are ${\mathfrak D}$-stable.

Note, that not any set of sign patterns of necessary sizes can be determined by a matrix. Some links between sign patterns are required.

\section{General $({\mathfrak D},{\mathcal G},\circ)$-stability theory: further development and open problems}

\subsection{Characterization of $({\mathfrak D},{\mathcal G},\circ)$-stability: open problems}
Now we consider the main problems, connected to the class of $({\mathfrak D},{\mathcal G},\circ)$-stable matrices.
\paragraph{Checking $({\mathfrak D}, \ {\mathcal G}, \ \circ)$-stability}
The following two approaches, as well as any of their combinations are often used for establishing $({\mathfrak D}, \ {\mathcal G}, \ \circ)$-stability.
 \begin{enumerate}
\item[\rm 1.] Imposing some additional conditions on a matrix $\mathbf A$. For some important cases, $\mathbf A$ is assumed to belong to a specific matrix class, defined by determinantal inequalities.
\item[\rm 2.] Considering some more wide or more narrow stability region ${\mathfrak D}$ or matrix class ${\mathcal G}$, to make a transition to studying another stability type which would be easier to characterize.
\end{enumerate}

 We start with the problem of major importance: given a stability region ${\mathfrak D}$, a matrix class ${\mathcal G}$ and an operation $\circ$, how to verify if a given $n \times n$ matrix $\mathbf A$ is $({\mathfrak D},{\mathcal G},\circ)$-stable, using just a finite number of steps? Note that we deal with the classes ${\mathcal G}$, that contains an infinite number of matrices.

Let us observe the modern state of the characterization problem for the most important partial cases, listed in Section 1.

\begin{enumerate}
\item[\rm 1.] {\bf Multiplicative $D$-stable matrices}. The problem of matrix $D$-stability characterization is one of the most important old problems of matrix stability. However, it still remains open. Besides the general characterization problem, easy-to-verify sufficient conditions for $D$-stability are of great interest.
\item[\rm 2.]{\bf Multiplicative $H$-stable matrices.} In spite of the characterization problem of multiplicative $D$-stability is still unsolved, the characterization problem of multiplicative $H$-stability has been solved (see \cite{CARL3}, \cite{CAS}).
\item[\rm 3.]{\bf Multiplicative and additive $H(\alpha)$-stable matrices.} Lying "between" $H$-stable and $D$-stable matrices, this class is not characterized yet. However, for some special partitions $\alpha$, a full characterization may be provided.
\item[\rm 4.]{\bf $D(\alpha)$-stable matrices.} The above is true also for this class, which lies "between" stable and $D$-stable matrices.
\item[\rm 5.]{\bf $D$-positive and $D$-aperiodic matrices.} Though some necessary conditions as well as some classes of $D$-positive matrices were studied in \cite{BAO}, this characterization problem is not solved and even has not been studied in full volume.
\item[\rm 6.]{\bf Schur $D$-stable matrices.} While the study of continuous-time linear systems leads to the multiplicative $D$-stability problem the study of a discrete-time case leads to Schur $D$-stability. The problem of characterizing Schur $D$-stable matrices is also not solved yet.
\item[\rm 7.]{\bf $D$-hyperbolic matrices.} This new matrix class, introduced in \cite{AB2}, have not been studied in full volume, though some examples and applications are considered.
\item[\rm 8.]{\bf Additive $D$-stable matrices.} This matrix class is widely studied by the same methods that are used for studying multiplicative $D$-stablity. However, attempts to characterize additive $D$-stability are also not succeed yet.
\item[\rm 9.]{\bf Hadamard $H$-stable matrices.} Since Hadamard products are used to characterize (multiplicative) $D$-stability  and diagonal stability, the study of different kinds of Hadamard ${\mathcal G}$-stability is a matter of further development.
\item[\rm 10.]{\bf $B_k$-stable and $B_k$-nonsingular matrices.} The characterization of these matrix classes is also an open problem. For some study, see \cite{DJD}.
\end{enumerate}

Together with the most important characterization problem, we should mention the following connected subproblems.

\paragraph{Describing new classes of $({\mathfrak D},{\mathcal G},\circ)$-stable matrices} Such classes are supposed to be characterized by some collection of easy-to-verify conditions. To obtain the description of a new class for the most general case, we are particularly interested in some easy-to-verify conditions of ${\mathfrak D}$-stability (for a given stability region ${\mathfrak D}$). For the exception of some well-known partial cases, this is a hard problem as is. For special cases of stability regions ${\mathfrak D}$, such as the left (right)-hand side of the complex plane, unit disk and real axes, a number of such conditions is obtained and used as a base of various $({\mathfrak D},{\mathcal G},\circ)$-stability criteria. For some classes of multiplicative $D$-stable matrices, see, for example, \cite{JOHN1}, \cite{DATTA}, for Schur $D$-stable see \cite{BHK} and \cite{PRY1}, for $D$-positive see \cite{BAO}, for additive $D$-stable see \cite{GEA}.

\paragraph{Proving $({\mathfrak D},{\mathcal G},\circ)$-stability of a given matrix class} Using general results (even if they are known) usually requires a huge amount of computations. That is why finding sufficient conditions is particularly useful. The matrices we study arise in analyzing specific mathematical models, thus they are likely to have some specific properties (e.g. symmetric positive definite, oscillatory, stochastic, $M$-matrices). The problem of proving $({\mathfrak D},{\mathcal G},\circ)$-stability of a naturally arisen matrix class characterized by its determinantal properties leads to a variety of unsolved matrix problems connected to the problems of stability of dynamical systems. We can express them as embedding relations between the class of stable matrices and other matrix classes. The most important are the question of the stability of $P^2$-matrices, asked by Hershkowitz and Johnson in \cite{HERJ1} and the question of the stability of strictly GKK $\tau$-matrices by Holtz and Schneider (see \cite{HOLS}).

\subsection{Further development of $({\mathfrak D},{\mathcal G},\circ)$-stability theory}
By analogy with already highly developed theory for partial cases (multiplicative and additive $D$-stability, Schur $D$-stability), here we provide some concept closely related to $({\mathfrak D},{\mathcal G},\circ)$-stability with the description of related problems.
\paragraph{Total $({\mathfrak D},{\mathcal G},\circ)$-stability} Here, we recall the following definition (see \cite{KAB}, p. 35). A property of an $n \times n$ matrix ${\mathbf A}$ is called {\it hereditary} if every principal submatrix of ${\mathbf A}$ shares it. The property of $({\mathfrak D},{\mathcal G},\circ)$-stability is not hereditary even in the classical case of multiplicative $D$-stability (see \cite{LOG}). Thus we introduce the following class. Given a stability region ${\mathfrak D}$, a matrix class ${\mathcal G} \subset {\mathcal M}^{k \times k}$, $k = 1, \ \ldots, \ n$ and a binary operation $\circ$ defined on ${\mathcal M}^{k \times k}$, $k = 1, \ \ldots, \ n$, a matrix $\mathbf A$ is called {\it totally $({\mathfrak D},{\mathcal G},\circ)$-stable} if it is $({\mathfrak D},{\mathcal G},\circ)$-stable and every its principal submatrix is also $({\mathfrak D},{\mathcal G},\circ)$-stable. As in the case of multiplicative $D$-stability, this matrix class may be used for studying properties of principal submatrices of $({\mathfrak D},{\mathcal G},\circ)$-stable matrices and for establishing necessary conditions for $({\mathfrak D},{\mathcal G},\circ)$-stability. Special cases of triples $({\mathfrak D},{\mathcal G},\circ)$, for which $({\mathfrak D},{\mathcal G},\circ)$-stability implies total $({\mathfrak D},{\mathcal G},\circ)$-stability are also of interest.

The class of (multiplicative) totally stable matrices was introduced in \cite{QR} (see \cite{QR}, p. 314), referring \cite{ME}, where a necessary condition for total stability was given. For the definition and study of this class see also \cite{KAB}. This class also arises in connection with further defined robust $D$-stability (see, for example, \cite{HART}, p. 205).
\paragraph{Inertia and inertia preservers} Here, we are restricted to studying specific stability regions $\mathfrak D$ with
${\rm int}({\mathfrak D}) \neq \emptyset$ and $\overline{{\mathfrak D}} \neq {\mathbb C}$. So we have three nonempty sets: ${\rm int}({\mathfrak D})$, $\partial({\mathfrak D})$ and ${\rm int}({\mathfrak D}^c) = {\mathbb C}\setminus \overline{{\mathfrak D}}$. The {\it inertia} of a square matrix ${\mathbf A}$ (with respect to a given domain ${\mathfrak D}$) is defined as a triple $(i_+({\mathbf A}), \ i_0({\mathbf A}), \ i_-({\mathbf A}))$, where $i_+({\mathbf A})$ $(i_-({\mathbf A}))$ is the number of the eigenvalues of ${\mathbf A}$ inside (respectively, outside) ${\mathfrak D}$, $i_0({\mathbf A})$ is the number of the eigenvalues on the boundary of ${\mathfrak D}$. Counting the number of eigenvalues in a given domain is also a problem of great importance in engineering. An $n \times n$ real matrix $\mathbf A$ is called {\it $({\mathfrak D},{\mathcal G},\circ)$-inertia preserving} if
     $$(i_+({\mathbf G}\circ{\mathbf A}), \ i_0({\mathbf G}\circ{\mathbf A}), \ i_-({\mathbf G}\circ{\mathbf A})) = (i_+({\mathbf G}), \ i_0({\mathbf G}), \ i_-({\mathbf G}))$$
     for every matrix ${\mathbf G} \in {\mathcal G}$.
Let us consider the partial cases.

 In the case, when ${\mathfrak D} = \{z \in {\mathbb C}: {\rm Re}(z) > 0\}$, we consider $i_+({\mathbf A})$ $(i_-({\mathbf A}))$ to be the number of the eigenvalues with positive (respectively, negative) real parts, $i_0({\mathbf A})$ to be the number of the eigenvalues with zero real parts, i.e. on the imaginary axes. The study of inertia preservers under a multiplication by a symmetric matrix $\mathbf H$ ($\mathcal G$ to be the class of symmetric matrices and $\circ$ to be matrix multiplication) was started by Sylvester and continued by Ostrowski and Schneider \cite{OSS} (see \cite{OSS}, p. 76, Theorem 1), where the key results, connecting inertia and stability were presented. These results were used to characterize the class of $H$-stable matrices. Classical results on this theme were obtained by Taussky \cite{TAU1}, Carlson and Schneider \cite{CAS}. An overview of this topic is presented in \cite{DATTA1}, where the inertia with respect to the unit disk is also considered. The inertia is used for the characterization of the class of $D$-stable matrices (see \cite{DATTA1}, p. 582 and references therein).

 For the generalized stability region $\mathfrak D$ and the same class of symmetric matrices $\mathcal H$, the characterization of inertia preservers is posed as an open problem in \cite{DATTA1} (p. 593, Problem 1). The tridiagonal case was considered in \cite{CARL1}, further generalization was provided in \cite{CHEN}.

\paragraph{Robustness of $({\mathfrak D},{\mathcal G},\circ)$-stability} Now we introduce one more concept of great importance in system theory. Here, we again consider a specific type of stability regions ${\mathfrak D}$, so-called {\it Kharitonov regions} (for the definitions and properties see, for example, \cite{SOF}). A matrix $\mathbf A$ is said to be {\it robustly $({\mathfrak D},{\mathcal G},\circ)$-stable} if it is $({\mathfrak D},{\mathcal G},\circ)$-stable and remains $({\mathfrak D},{\mathcal G},\circ)$-stable for sufficiently small perturbations of ${\mathbf A}$. In other words, ${\mathbf A}$ is robustly $({\mathfrak D},{\mathcal G},\circ)$-stable if ${\mathbf A}$ is $({\mathfrak D},{\mathcal G},\circ)$-stable and there exists an $\epsilon > 0$ such that for any real-valued matrix ${\mathbf \Delta}$ with $\|{\mathbf \Delta}\| < \epsilon$, the matrix ${\mathbf A} + {\mathbf \Delta}$ is $({\mathfrak D},{\mathcal G},\circ)$-stable.

Note, that in general, $({\mathfrak D},{\mathcal G},\circ)$-stability is not a robust property, even in the classical case of multiplicative $D$-stability (see \cite{AB1} for the corresponding examples). Thus discovering sufficient conditions which lead to the classes of robustly $({\mathfrak D},{\mathcal G},\circ)$-stable matrices is of great importance.

 The class of robust $D(\alpha)$-stable matrices was analyzed in \cite{AB1} (p. 3, Definition 3), see also \cite{AB3}.
 In the same paper \cite{AB1} robustly $D$-hyperbolic and $D(\alpha)$-hyperbolic classes are analyzed.

\paragraph{$\mathfrak D$-stability measurement and general $\mathfrak D$-stabilization problem}
Here, we introduce the concept and state some problems that are connected to robust $\mathfrak D$-stability. We start with the following question, asked by Hershkowitz (see \cite{HER1}, p. 162).

Given a $P$-matrix $\mathbf A$, how far is it from being stable?

 He outlined two directions for giving an answer:

\begin{enumerate}
\item[-] in terms of the width of a wedge around the negative direction of the real axes, which is free from eigenvalues;
\item[-] in terms of the inertia of $\mathbf A$ (how much eigenvalues are located in the closed left-hand side of the complex plane).
\end{enumerate}
The combination of this two approaches was used in \cite{HERB}, \cite{KEL}.

Here, we state the following more general problem.

{\bf Problem 7}. Given an arbitrary stability region ${\mathfrak D} \subset \mathbb{C}$, and a matrix $\mathbf A$ from ${\mathcal M}^{n \times n}$, how far is $\mathbf A$ from being ${\mathfrak D}$-stable?

The answer may use the combination of the following approaches:
\begin{enumerate}
\item[-] description of the new stability region ${\mathfrak D}_1$ such that ${\mathfrak D} \subseteq {\mathfrak D}_1$ and $\sigma(\mathbf A) \subset {\mathfrak D}_1$;
\item[-] counting the inertia of $\mathbf A$ with respect to the stability region $\mathfrak D$.
\end{enumerate}

Another problem, mentioned in \cite{HER1} is {\it multiplicative $D$-stabilization problem} (see \cite{HER1} p. 162, then p.170): given a square real-valued matrix $\mathbf A$, can we find a diagonal matrix $\mathbf D$ such that ${\mathbf D}{\mathbf A}$ is positive stable? Simple example with a circulant matrix shows that it is not always possible. For the results on a stabilization of matrices using a diagonal matrix, we refer to \cite{BAL}, \cite{YAR}, \cite{LOC}.

In full generality, we state this problem as follows:

{\bf Problem 8}. Given a matrix $\mathbf A$ from ${\mathcal M}^{n \times n}$, an arbitrary stability region ${\mathfrak D} \subset \mathbb{C}$, a matrix class ${\mathcal G} \subset {\mathcal M}^{n \times n}$ and a binary matrix operation $\circ$, when it is possible to find a matrix ${\mathbf G}_0 \in {\mathcal G}$ such that $\sigma({\mathbf G}_0 \circ{\mathbf A}) \subset {\mathfrak D}$?

A matrix $\mathbf A$ is called {\it $({\mathfrak D},{\mathcal G},\circ)$-stabilizable} if the answer to Problem 8 is affirmative. As it follows from the definition, the class of $({\mathfrak D},{\mathcal G},\circ)$-stable matrices belongs to the class of $({\mathfrak D},{\mathcal G},\circ)$-stabilizable matrices.
\paragraph{$({\mathfrak D},{\mathcal G},\circ)$ stability measurement and $({\mathfrak D},{\mathcal G},\circ)$-stabilization problem}
Here, we ask the following more specific question.

{\bf Problem 9.} Given a ${\mathfrak D}$-stable matrix $\mathbf A$, a matrix class ${\mathcal G} \subset {\mathcal M}^{n \times n}$ and a binary matrix operation $\circ$, how far $\mathbf A$ is from being $({\mathfrak D},{\mathcal G},\circ)$-stable? The directions of giving the answer to this question are as follows.
 \begin{enumerate}
 \item[-] Describing subclasses ${\mathcal G}_1$ of the class $\mathcal G$, such that ${\mathcal G_1} \subseteq {\mathcal G}$ (or conversely) and $\sigma({\mathbf G}\circ{\mathbf A}) \subset {\mathfrak D}$ for every ${\mathbf G} \in {\mathcal G}_1$. Note, that every $({\mathfrak D},{\mathcal G},\circ)$-stabilizable matrix can be considered as $({\mathfrak D},{\mathcal G_1},\circ)$-stable for some nonempty class ${\mathcal G_1} \subseteq {\mathcal G}$.
\item[-] Describing the new stability region ${\mathfrak D}_1$ such that ${\mathfrak D} \subseteq {\mathfrak D}_1$ and $\sigma({\mathbf G}\circ{\mathbf A}) \subset {\mathfrak D}_1$ for every ${\mathbf G} \in {\mathcal G}$.
\item[-] Counting the inertia of ${\mathbf G}\circ{\mathbf A}$ with respect to the stability region $\mathfrak D$ while ${\mathbf G}$ is varying along the class ${\mathcal G}$.
\end{enumerate}
We may also use the combinations of the described above approaches.

As examples of partial multiplicative $D$-stability, we mention the classes of $D(\alpha)$-stable matrices and $D_{\tau}$-stable matrices.
\paragraph{Relations between different classes of $({\mathfrak D},{\mathcal G},\circ)$-stable matrices}
Besides of relations between different $({\mathfrak D},{\mathcal G},\circ)$-stability classes, described in Section 2, based on inclusion relations between stability regions and matrix classes, relations between classes, defined by different binary operations are of interest. In general form, the problem is stated as follows.

{\bf Problem 10.} Given two triples $({\mathfrak D}_1,{\mathcal G}_1,\circ)$ and $({\mathfrak D}_2,{\mathcal G}_2,\star)$, do the corresponding classes of $({\mathfrak D},{\mathcal G},\circ)$-stable matrices intersect? For which classes of matrices $\mathbf A$ do they coincide?

The relations between matrix classes, defined in Section 1 are investigated in various ways. We do not provide any diagrams here, just refer to the following papers. The relations between Lyapunov diagonally stable, multiplicative and additive $D$-stable matrices were first studied in \cite{CROSS}. In \cite{HER1}, p. 173, the diagram showing relations between Lyapunov diagonally stable, multiplicative and additive $D$-stable matrices, is provided. For some matrix types, different stability types, namely multiplicative and additive $D$-stability classes coincide (\cite{HER1}, p. 174). For the relations between matrix classes, we refer to \cite{LOG}, where multicomponent diagrams are presented, see also \cite{CADHJ}, p. 154, Fig 1, \cite{BHK3}, \cite{CALN}. The book \cite{KAB} provides a lot of information on this topic. The relations between Hadamard $H$-stability, Lyapunov diagonal stability and $D$-stability were first considered in \cite{JOHN3}, p. 304.
\paragraph{Further development: from matrices to other objects}
Here, we briefly mention natural generalizations of $D$-stability which arises during study of nonlinear systems (see \cite{BEF} and references therein), theory of $D$-stability for polynomial matrices (see \cite{HENBS}), recent studies of multidimensional matrices (tenzors) and so on.

\part{Applications}
\section{Robustness of linear systems}
\subsection{Continuous-time case}
Given a continuous-time linear system
\begin{equation}\label{sys}\dot{x}(t) = {\mathbf A}x(t), \end{equation}
where $x(t) \in {\mathbb R}^n$ is the system state vector and ${\mathbf A} \in {\mathcal M}^{n \times n}$ is the time-invariant system matrix. Recall, that the linear system \eqref{sys} is called {\it asymptotically stable} if every finite initial state excites a bounded response, which, in addition, converges to $0$ as $t \rightarrow 0$. The system \eqref{sys} is asymptotically stable if and only if all eigenvalues of ${\mathbf A}$ has negative real parts (see \cite{CHEN2}).

Let the system matrix $\mathbf A$ be $({\mathfrak D},{\mathcal G},\circ)$-stable with respect to the stability region $\mathfrak D ={\mathbb C}^-$, a matrix class ${\mathcal G} \subset {\mathcal M}^{n \times n}$ and a binary matrix operation $\circ$. Then each system of the perturbed family
\begin{equation}\label{sysper}\dot{x}(t) = ({\mathbf G}\circ{\mathbf A})x(t), \ \qquad {\mathbf G}\in {\mathcal G} \end{equation}
is asymptotically stable.

The above statement directly follows from the definition of $({\mathbb C}^-,{\mathcal G},\circ)$-stability. It is particularly useful for the cases, when the entries the system matrix ${\mathbf A}$ are known (or should be considered) up to a perturbation of a specified form.

{\bf Example.} Consider the following class of biological models studied in \cite{BFG}:
\begin{equation}\label{Bio}
\dot{x} = {\mathbf S}g(x)+ {\mathbf V}u,
\end{equation}
where $x \in {\mathbb R}^{n}$ represents the concentration of each biological species in the system, $u$ is the vector of constant influxes or outfluxes, $g(x)$ is a vector of reaction rates and each component $g_i(\cdot)$, $i = 1, \ \ldots, \ n$ is a positive monotone function, $\mathbf S$ is a system stoichiometry matrix. This kind of systems can be described by a set of biochemical reactions (for more details see \cite{DOP}).

The Jacobian of System \eqref{Bio} can be represented as a product of the form
$${\mathbf J} = {\mathbf B}{\mathbf D}{\mathbf C},$$
where ${\mathbf D}$ is a diagonal matrix, whose diagonal entries are partial derivatives of the reaction rates (assumed to be positive).

Analysis of robust stability of matrix $J$ leads to the study of the matrix  $\widetilde{{\mathbf J}} =( {\mathbf C}{\mathbf B}){\mathbf D},$
and to establishing its stability for all positive diagonal matrices ${\mathbf D}$, which is exactly the property of multiplicative $D$-stability.

\subsection{Hopf bifurcation phenomena}
Considering an unstable continuous-time systems, we get the following two cases
\begin{enumerate}
\item[\rm 1.] Exponential instability (the system matrix $\mathbf A$ has a real positive eigenvalue).
\item[\rm 2.] Oscillatory instability (the system matrix $\mathbf A$ has a pair of complex eigenvalues with nonnegative real parts).
\end{enumerate}

Different types of $({\mathfrak D},{\mathcal G}, \circ)$-stability are used to exclude the cases of exponential instability and oscillatory instability in systems of general form.

{\bf Example.} To exclude the case of oscillatory instability for the systems of the form \eqref{Abbb}
\begin{equation}\label{Abbb}\left\{\begin{array}{cc} \dot{x} = {\mathbf A}x + {\mathbf B}y & \\
{\mathbf E}(\epsilon)\dot{y} = {\mathbf C} x + {\mathbf D} y, \end{array}\right.
\end{equation}
where ${\mathbf E}(\epsilon)$ is an $\alpha$-scalar matrix of the form
$${\mathbf E}(\epsilon) = {\rm diag}\{\epsilon_1{\mathbf I}_{\alpha_1}, \ \ldots, \ \epsilon_m{\mathbf I}_{\alpha_m}\}, $$
the concept of $D$-hyperbolicity is used (see \cite{AB2}).

\subsection{Discrete-time case}
Given a system of difference equations
\begin{equation}\label{sysd}x[k+1] = {\mathbf A}x[k]. \end{equation}
The stability concept for difference systems is analogical to the one for continuous systems. The linear system \eqref{sysd} is called {\it asymptotically stable} if every finite initial state excites a bounded response, which, in addition, approaches $0$ as $k \rightarrow \infty$. The system \eqref{sysd} is asymptotically stable if and only if all eigenvalues $\lambda$ of ${\mathbf A}$ satisfy $|\lambda| < 1$ (see \cite{CHEN2}).

Let the system matrix $\mathbf A$ be $({\mathfrak D},{\mathcal G},\circ)$-stable with respect to the stability region $\mathfrak D = {\mathcal D}(0,1)$, a matrix class ${\mathcal G} \subset {\mathcal M}^{n \times n}$ and a binary matrix operation $\circ$. Then each system of the perturbed family
\begin{equation}\label{sysper}x[k+1] = ({\mathbf G}\circ{\mathbf A})x[k], \ \qquad {\mathbf G}\in {\mathcal G} \end{equation}
is asymptotically stable.

\subsection{Fractional differential systems} Consider the following concept which appeared recently and have applications in viscoelasticity, acoustics, polymeric chemistry, etc (see \cite{MAT}, \cite{MAT1} for the theory and references therein for the applications).

Given a linear system in the following form:
\begin{equation}\label{sysf}d^{\theta}x(t) = {\mathbf A}x(t), \end{equation}
with $0 < \theta \leq 1$, $x(0) = x_0$. As in the case of System \eqref{sys} (which corresponds to $\theta=1$), the linear system \eqref{sysf} is called {\it asymptotically stable} if every finite initial state excites a bounded response, which, in addition, converges to $0$ as $t \rightarrow 0$. The system \eqref{sysf} is asymptotically stable if and only if all eigenvalues $\lambda$ of ${\mathbf A}$ satisfy $|\arg(\lambda)|> \theta\dfrac{\pi}{2}$ (see \cite{MAT}). This corresponds to ${\mathfrak D}$-stability with respect to the stability region ${\mathfrak D} = {\mathbb C}\setminus {\mathbb C}^+_{\theta}$.

Let the system matrix $\mathbf A$ be $({\mathfrak D},{\mathcal G},\circ)$-stable with respect to ${\mathfrak D} = {\mathbb C}\setminus {\mathbb C}^+_{\theta}$, a matrix class ${\mathcal G} \subset {\mathcal M}^{n \times n}$ and a binary matrix operation $\circ$. Then each system of the perturbed family
\begin{equation}\label{sysfper}d^{\theta}x(t)  = ({\mathbf G}\circ{\mathbf A})x(t), \ \qquad {\mathbf G}\in {\mathcal G} \end{equation}
is asymptotically stable.

 \subsection{Robust eigenvalue localization: general problem} In practice, studying dynamic systems, some perturbations of a system matrix may occur, and, in general, the matrix entries may be known up to some small values (for example, caused by linearization error). One of the most important system dynamics problems (see \cite{BAR}) is as follows. Given a stability region ${\mathfrak D} \subset {\mathbb C}$ and a matrix ${\mathbf A} \in {\mathcal M}^{n \times n}$, when a perturbed matrix $\widetilde{\mathbf A} = {\mathbf A} + {\mathbf \Delta}$ is ${\mathfrak D}$-stable?
A number of papers are devoted to studying this property, called {\it robust stability}, with respect to different stability regions $\mathfrak D$ (see, for example, \cite{ASCM} for EMI regions).

Sometimes, special types of perturbations are considered or some information of ${\mathbf \Delta}$ is provided, and we have the following description of the {\it uncertain matrix} $\widetilde{\mathbf A}$:
$$\widetilde{\mathbf A} = {\mathbf A} + {\mathbf U}{\mathbf \Delta}{\mathbf V}, $$
where ${\mathbf U}$, ${\mathbf V}$ are known matrices and introduced to specify the {\it structure of uncertainty}, ${\mathbf \Delta}$ is bounded by its norm. In some cases, we can easily come from studying $({\mathfrak D},{\mathcal G},\circ)$-stability to studying robust $\mathfrak D$-stability problem with some special structure of uncertainty.

{\bf Example 1.} For the case of $({\mathfrak D},{\mathcal G}, +)$-stability (i.e. the operation $\circ$ is matrix addition), we have
to check, if $\sigma({\mathbf A} + {\mathbf G}) \subset {\mathfrak D}$ for every matrix ${\mathbf G} \in {\mathcal G}$. Thus, assuming the norm of ${\mathbf G}$ to be sufficiently small, we immediately obtain robust ${\mathfrak D}$-stability problem with a specified structure of uncertainty (from the class ${\mathcal G}$).

{\bf Example 2.} Considering multiplicative or Hadamard $({\mathfrak D},{\mathcal G})$-stability and using the distributivity law, we obtain:
$${\mathbf G}\circ{\mathbf A} = ({\mathbf I} + ({\mathbf G} - {\mathbf I}))\circ{\mathbf A} = {\mathbf A} + ({\mathbf G} - {\mathbf I})\circ{\mathbf A}.$$
Thus, assuming that $\|{\mathbf G} - {\mathbf I}\|$ is sufficiently small, we obtain that every multiplicative (respectively, Hadamard) $({\mathfrak D},{\mathcal G})$-stable matrix is robustly ${\mathfrak D}$-stable with the uncertainty structure $({\mathbf G} - {\mathbf I})\circ{\mathbf A}$.

{\bf Example 3.} Considering the operation of entry-wise maximum $\oplus_m$, we obtain the class of $({\mathfrak D},{\mathcal G}, \oplus_m)$-stable matrices, that for a specific choice of ${\mathcal G}$ can be considered as an interval matrix. Using the commutativity and distributivity laws:
$${\mathbf G}\oplus_m {\mathbf A} = {\mathbf A}\oplus_m{\mathbf G} = {\mathbf A} + (-{\mathbf A}) + ({\mathbf A}\oplus_m{\mathbf G}) = $$ $$
{\mathbf A} + ({\mathbf A} - {\mathbf A})\oplus_m({\mathbf G} - {\mathbf A}) = {\mathbf A} + {\mathbf O}\oplus_m({\mathbf G} - {\mathbf A}).$$
Thus, for small values of $\|{\mathbf G} - {\mathbf A}\|$, $({\mathfrak D},{\mathcal G}, \oplus_m)$-stability problem leads to robust ${\mathfrak D}$-stability problem with the uncertainty structure ${\mathbf O}\oplus_m({\mathbf G} - {\mathbf A})$.

\section{Global asymptotic stability and diagonal Lyapunov functions}
Here, we recall the following definition from the theory of differential equations (see, for example, \cite{BACR}). Consider a nonlinear system of ordinary differential equations of the form
\begin{equation}\label{Nonsyst}
\dot{x} = f(x),
\end{equation}
where $x \in {\mathbb R}^n$. The system \eqref{Nonsyst} is called {\it stable at the origin} if for each $\epsilon > 0$ there is $\delta > 0$ such that for each solution $x$ if $\|x (0)\| < \delta$ then $\|x (t)\| < \epsilon$ for every $t \geq 0$ and (locally) {\it asymptotically stable at the origin} if, in addition, there exists $\delta_0 > 0$ such that $\lim_{t \rightarrow + }\|x(t)\| = 0$ for each $x$ such that $\|x(0)\| < \delta_0$. The system \eqref{Nonsyst} is called {\it globally asymptotically stable at the origin} if $\delta_0$ can be taken arbitrarily large.

The system \eqref{Nonsyst} corresponds to the family of Jacobian matrices $${\mathcal J}:=\{{\mathbf J}_{f(x)}\}_{x \in {\mathbb R}^n}.$$ In general, stability of all the family ${\mathcal J}$ does not guarantee the global asymptotic stability of \eqref{Nonsyst}. However, the concept of diagonal stability helps to describe some families of nonlinearities, for which global asymptotic stability holds.

\subsection{Continuous-time case}
Consider perturbed systems of the form
\begin{equation}\label{sysnl} \dot{x}(t) = {\mathbf A}(f(x(t))) \end{equation}
where $x(t) \in {\mathbb R}^n$, $f(t)=(f_1(\cdot), \ \ldots, \ f_n(\cdot))$ is a vector-function with each coordinate $f_i: {\mathbb R} \rightarrow {\mathbb R}$ be a continuous function, satisfying the conditions:
$$f_i(\xi)\xi > 0, \qquad f_i(0)=0, \ \qquad i = 1, \ \ldots, \ n.$$
$$ \int_0^{x_i}f_i(\tau)d\tau\rightarrow \infty \ \mbox{as} \ |x_i|\rightarrow \infty \ \qquad i = 1, \ \ldots, \ n.$$

It was first shown by Persidskii \cite{PER} that {\it the equilibrium $x = 0$ of all the systems \eqref{sysnl} is globally asymptotically stable if the system matrix ${\mathbf A}$ is diagonally stable}. This result was extended by Kaszkurewicz and Bhaya \cite{KAB2} to a more general class of time-dependent system perturbation described by

\begin{equation}\label{sysnlper} \dot{x}(t) = ({\mathbf A}\circ {\mathbf F})(x,t), \end{equation}
where $\circ$ denotes Hadamard product, $x(t) \in {\mathbb R}^n$, ${\mathbf F}(\cdot,\cdot)$ is matrix function with each entry $f_{ij}: {\mathbb R}\times{\mathbb R}_{+} \rightarrow {\mathbb R}$ be a continuous function, satisfying the conditions:
$$f_{ij}(\xi,t)\xi > 0, \qquad f_{ij}(0,t)=0, \ \qquad i,j = 1, \ \ldots, \ n.$$
In \cite{KAB2}, the conditions of global stability of the equilibrium of \eqref{sysnlper} were given in terms of entry-wise diagonal dominance of the matrix-function ${\mathbf F}$ and diagonal stability of the matrix $W({\mathbf A})$. Recall, that $W({\mathbf A}) = {\mathbf A}^{W} = \{a_{ij}^W\}_{i,j = 1}^n, $
where $$a_{ij}^W = \left\{\begin{array}{cc}a_{ij} & i = j \\ |a_{ij}| & i \neq j \\ \end{array}\right.$$

This kind of results include perturbations of multiplicative and additive types as well as Hadamard and block Hadamard products.

The system of the form
$$\dot{x}(t) = f({\mathbf A}x(t)) $$
is studied by the same methods.
\subsection{Discrete-time case}
Discrete-time perturbed systems were considered in \cite{KAB2} analogically.

\begin{equation}\label{sysnldper} x(k+1) = ({\mathbf A}\circ {\mathbf \Phi})(x,k), \end{equation}
where $\circ$ denotes Hadamard product, $x(k) \in {\mathbb R}^n$, ${\mathbf \Phi}(\cdot,\cdot)$ is a matrix function with each entry $\phi_{ij}: {\mathbb R}\times{\mathbb Z}_{+} \rightarrow {\mathbb R}$ satisfy the conditions:
$$|\phi_{ij}(\xi,k)| \leq |\xi|, \qquad \phi_{ij}(0,t)=0, \ \qquad i,j = 1, \ \ldots, \ n, \ \ k = 0, \ 1, \ \ldots.$$
In \cite{KAB2}, the following statement was proved.
\begin{theorem}\cite{KAB2} The equilibrium $x = 0$ of all systems in the class \eqref{sysnldper} is globally asymptotically stable if the matrix $|{\mathbf A}| = \{|a_{ij}|\}$ is Schur diagonally stable.
\end{theorem}

\section{Passivity and network stability analysis}
First, let us recall some definitions and notations from nonlinear control theory (see \cite{SJK}). Given a (nonlinear) system $H$ of ordinary differential equations of the form

\begin{equation}\label{Pass}H: \left\{\begin{array}{cc} \dot{x} = f(x,u) & \\
y = h(x,u), \end{array}\right.
\end{equation}
where $x(t) \in {\mathbb R}^n$ is a {\it state vector}, $u(t) \in {\mathbb R}^m$ is an {\it input vector} and $y(t) \in {\mathbb R}^m$ is an {\it output vector}. Assume that the system $H$ has a stationary solution (equilibrium) $x^*$ for the corresponding $u^*$ and $y^*$.

We say that $H$ is {\it output strictly passive} with respect to the equilibrium $x^*$ if it is dissipative with respect to the function
$$w(u - u^*, y-y^*)= (u - u^*)^T(y - y^*) - \sigma^{-1}|y - y^*|^2, \ \qquad \sigma > 0, $$
i.e. there exists a {\it storage function} $V(x(t))$ which is differentiable by $t$ and satisfies the condition
$$\dot{V}(x(t)) \leq w(u(t), y(t)).$$ Here $\sigma > 0$ is called a {\it gain}.

 \subsection{Diagonal stability and passivity}
 Following \cite{ARCS1}, consider a sequence of biochemical reactions, where the end product drives the first reaction.
$$\dot{x}_1 = -f_1(x_1) + g_n(x_n), $$
$$\dot{x}_2 = -f_2(x_2) + g_1(x_1), $$
$$\ldots \ldots \ldots $$
$$\dot{x}_n = -f_n(x_n) + g_{n-1}(x_{n-1}),$$
where $f_i(\cdot)$, $i = 1, \ \ldots, \ n$ and $g_i(\cdot)$, $i = 1, \ \ldots, \ n-1$ are increasing functions and $g_n(\cdot)$ is a decreasing function (for the biological background see \cite{TYO}, \cite{THR1}, \cite{THR2}).
The cyclic interconnection structure of the sequence of systems $H_1, \ \ldots, \ H_n$, where each $H_i$ is illustrated as follows
$$ 
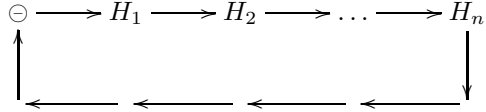
\xymatrix{
  \circleddash \ar[r]^{} & H_1  \ar[r]^{} & H_2  \ar[r]^{} & \ldots  \ar[r]^{} & H_n \ar[d]^{} \\
   \ar[u]^{} &  \ar[l]^{} &  \ar[l]^{} &  \ar[l]^{} &  \ar[l]^{}  }$$

The Jacobian linearization at the equilibrium is of form \eqref{STR2}, which is shown to be diagonally stable under certain conditions (see \cite{ARCS1} and the review part). Thus the global asymptotic stability of the above system is established.

Diagonal stability concept is used for construction a composite Lyapunov function for an interconnected system (see \cite{ARC2}).
Let us consider the {\it cascade interconnection} where each $H_i$ is output strictly passive:
$$\xymatrix{
  x_1 \ar[r]^{} & H_1 \ar[r]^{} & H_2 \ar[r]^{} & \ldots \ar[r]^{} & H_n \ar[r]^{} & x_n  } $$

Then the diagonal stability is used to prove an {\it input feedforward passivity (IFP)} property for the cascade, which quantifies the amount of feedforward gain required to re-establish passivity (see \cite{ARCS1}, p. 1536, Corollary 5).

In \cite{ARC2}, a more general problem of interconnection of dynamical systems $H_i$, $i = 1, \ \ldots, \ n$ is considered. The systems $H_i$ are connected according to a feedback law, defined by
$$u = ({\mathbf K}\otimes {\mathbf I}_m)y, $$
where ${\mathbf K} \in {\mathcal M}^{n \times n}$, ${\mathbf I}_m$ is $m \times m$ identity matrix, $u = (u_1^T, \ \ldots, \ u_n^T)^T$ is a general input vector constructed from the input vectors $u_i \in {\mathbb R}^m$ of the systems $H_i$, respectively, $y = (y_1^T, \ \ldots, \ y_n^T)^T$ is a general output vector. Suppose each component $H_i$ is output strictly passive relative to its equilibrium $x_i^*$.
It was proved in \cite{ARC2} that {\it the equilibrium $x^* = ((x_1^*)^T, \ \ldots, \ (x_n^*)^T)^T$ of the interconnected system is stable if the matrix
$${\mathbf E} = -{\mathbf I} + {\mathbf \Sigma}{\mathbf K}$$
is diagonally stable, where ${\mathbf \Sigma} = {\rm diag}\{\sigma_1, \ \ldots, \ \sigma_n\}$ is a diagonal matrix, constructed using the gains $\sigma_i$}.

In \cite{FER}, applications of diagonal stability to stochastic systems were considered.
\subsection{Schur diagonal stability and Small Gain Theorem}
Consider a linear interconnection of the systems $H_i$. The following stability criteria was proved in \cite{MOYH} (see \cite{MOYH}, p. 146, Theorem 5): {\it Let ${\mathbf H}$ be a matrix of the interconnected system, ${\mathbf \Sigma} = {\rm diag}\{\sigma_1, \ \ldots, \ \sigma_n\}$ is a diagonal matrix, constructed using the gains $\sigma_i$. Then the interconnected system is stable if the matrix ${\mathbf A}:={\mathbf \Sigma}{\mathbf H}$ is Schur diagonally stable}.

See also \cite{DARW} for the applications of Schur diagonal and $D$-stability to input-to state stability of large-scale interconnected systems.

\section{Classical dynamic models}
\subsection{Enthoven--Arrow dynamic model}
The concepts of multiplicative $D$- and $H$-stability are included in a number of books in economic theory (see, for example, \cite{BED}, \cite{KEK}, \cite{QS}, \cite{WOOD}). They are based on the following classical dynamic model (see \cite{ENT}, \cite{AM}).

\begin{equation}\label{PR} \dot{p} = ({\mathbf K}^{-1} - {\mathbf B}\eta)^{-1}({\mathbf Q} + {\mathbf B})(p-p_0),
\end{equation}
where $p \in {\mathbb R}^n$ is a vector of market prices (i.e. each its component $p_i$ is the price of the $i$th good exchanged in a competitive market).

Denoting ${\mathbf G}:= ({\mathbf K}^{-1} - {\mathbf B}\eta)^{-1}$ and ${\mathbf A} := {\mathbf Q} + {\mathbf B}$, and putting $p_0:= 0$, we obtain from the system \eqref{PR}

\begin{equation}\label{PR1} \dot{p} = {\mathbf G}{\mathbf A}p,
\end{equation}
where the matrix ${\mathbf G}$ is considered to be symmetric positive definite (multiplicative $H$-stability) or positive diagonal ($D$-stability).

\subsection{Lotka--Volterra model}

Consider the following model of population dynamics in a community of $n$ biologic species.

\begin{equation} \dot{x} = x({\mathbf E} - {\mathbf \Gamma}x),
\end{equation}
where $x$ shows the intrinsic rate of natural increase, ${\mathbf E}$ is a vector parameter, $-{\mathbf \Gamma}$ is a matrix which shows interacting with other species.

It was noted by Kosov \cite{KOS} that the diagonal matrix ${\mathbf \Gamma}$ is proportional to the velocity vector of species reproduction in isolated state in an absence of self-limited factors, thus it has upper and lower bounds.

The study of Lotka--Volterra model is based on the concept of diagonal stability.

\subsection{Reaction-diffusion system}
Consider the following reaction-diffusion system with Neumann boundary condition (see, for example, \cite{WANGL}).
\begin{equation}\label{REDiff} u_t = \left\{\begin{array}{ccc}{\mathbf D}\Delta u + f(u) \\
\dfrac {\partial u}{\partial v} = 0 \\
u(x, 0) = u_0(x)
\end{array} \right.\end{equation}
where ${\mathbf D}$ is the matrix of diffusion coefficients.

By linear approximation at $u = 0$, the system \eqref{REDiff} can be reduced to the following form

\begin{equation}\label{REDiff2}\left\{\begin{array}{ccc} v_t = {\mathbf D}\Delta v + {\mathbf A}v \\
\dfrac {\partial u}{\partial v} = 0 \\
u(x, 0) = u_0(x)\end{array} \right.
\end{equation}

The asymptotic stability of the above system is established through additive $D$-stablity of $\mathbf A$.

\subsection{Time-invariant multiparameter singular perturbation problem}

Consider the following system (see \cite{AB1})

\begin{equation}\label{Abbb} \left\{\begin{array}{cc} \dot{x} = {\mathbf A}x + {\mathbf B}y & \\
\epsilon_i\dot{y}_i = {\mathbf C}_i x + {\mathbf D}_i y & i = 1, \ \ldots, \ m. \end{array}\right.
\end{equation}
Here small parameters $\epsilon_i > 0$ for each $i = 1, \ \ldots, \ m$ and all the ratios $\frac{\epsilon_i}{\epsilon_j}$ are supposed to be bounded.

This kind of systems is connected to the concept of multiplicative $D(\alpha)$-stability.

The robust $D(\alpha)$-stability is applied to study boundary layer systems of the form:

\begin{equation}
{\mathbf E}(\epsilon)\dot{z} = {\mathbf D}z.
\end{equation}

\section*{Acknowlegements}
The author would like to thank Mikhail Tyaglov and Alexander Dyachenko for many helpful comments. The research was supported by the National Science Foundation of China grant number
11750110414 and partially supported by jointly funded by the National Natural Science Foundation of China and the Israel Science Foundation grant number 11561141001.

\part{Appendix. The dictionary of matrix classes}
\section*{Matrix constructions}
\begin{enumerate}
\item[\rm 1.] {\bf Sign pattern class and sign stability.} Given an $n \times n$ matrix ${\mathbf A} = \{a_{ij}\}_{i,j = 1}^n$, its {\it sign pattern} ${\mathbf S}({\mathbf A})$ is an $n \times n$ matrix defined by:
 $${\mathbf S}({\mathbf A}) = \{s_{ij}\}_{i,j = 1}^n, \ \mbox{where} \ s_{ij} = {\rm sgn}(a_{ij}), \ i,j = 1, \ \ldots, \ n.$$
 Two matrices ${\mathbf A}$ and ${\mathbf B}$ are called {\it sign similar} if ${\mathbf S}({\mathbf A}) = {\mathbf S}({\mathbf B})$. Denote ${\mathcal A}$ the set of all matrices, sign-similar to a given matrix ${\mathbf A}$. Then ${\mathbf A}$ is called {\it sign-stable} or {\it qualitative stable}, if any matrix from ${\mathcal A}$ is stable.
\item[\rm 2.] {\bf Comparison matrix.} A matrix $M({\mathbf A}) = \{\widetilde{a}_{ij}\}_{i,j = 1}^n$ is called a {\it comparison matrix} of a matrix ${\mathbf A} = \{a_{ij}\}_{i,j = 1}^n$, if $$\widetilde{a}_{ij} = \left\{\begin{array}{cc}|a_{ij}| & i = j \\ -|a_{ij}| & i \neq j \\ \end{array}\right.$$
\item[\rm 3.] {\bf Compound matrix.} The {\it $j$th compound matrix} ${\mathbf A}^{(j)}$ $(1 \leq j \leq n)$ of an $n \times n$ matrix $\mathbf A$ is a matrix that consists of all the minors
$A\begin{pmatrix}
  i_1 & \ldots & i_j \\
  k_1 & \ldots & k_j
\end{pmatrix}$, where $1 \leq i_1 < \ldots < i_j \leq n, \ 1 \leq \ k_1 < \ldots < k_j \leq n$,  of the initial matrix ${\mathbf A}$. The minors are listed in the lexicographic order. The matrix ${\mathbf A}^{(j)}$ is $\binom{n}j \times \binom{n}j$ dimensional, where $\binom{n}j = \dfrac{n!}{j!(n-j)!}$. The first compound matrix ${\mathbf A}^{(1)}$ is equal to ${\mathbf A}$.
\item[\rm 4.] {\bf Additive compound matrix.} Given an $n \times n$ matrix ${\mathbf A} = \{a_{ij}\}_{i,j = 1}^n$ and an $n \times n$ identity matrix ${\mathbf I} = \{\delta_{ij}\}_{i,j = 1}^n$, {\it the second additive compound matrix} ${\mathbf A}^{[2]}$ is a matrix that consists of the sums of minors of the following form:
$$a^{[2]}_{\alpha\beta} = \begin{vmatrix} a_{ik} & \delta_{il} \\ a_{jk} & \delta_{jl} \end{vmatrix} + \begin{vmatrix} \delta_{ik} & a_{il} \\ \delta_{jk} & a_{jl} \end{vmatrix}, $$
where $\alpha = (i,j), \  1 \leq i < j \leq n$, $\beta = (k,l), \  1 \leq k < l \leq n$,
listed in the lexicographic order. The matrix  ${\mathbf A}^{[2]}$ is $\binom{n}2 \times \binom{n}2$ dimensional (see \cite{FID}).
\end{enumerate}
\section*{Matrix classes}
Here, we list the definitions of the main matrix classes used above. An $n \times n$ real matrix ${\mathbf A}= \{a_{ij}\}_{i,j = 1}^n$ is called:
 \begin{enumerate}
\item[\rm 1.] {\it symmetric positive definite (semidefinite)} if $a_{ij} = a_{ji}$, $i,j = 1, \ \ldots, \ n$ and $x^T{\mathbf A}x > 0$ (respectively, $\geq 0$) for every nonzero vector $x \in {\mathbb R^n}$. 
\item[\rm 2.] {\it diagonal}, if it has with nonzero entries on its principal diagonal, while the rest are zero.
\item[\rm 3.]  {\it  row diagonally dominant} if the following inequalities hold:
 \begin{equation}\label{ROW}|a_{ii}| \geq \sum_{i\neq j}|a_{ij}| \qquad i = 1, \ \ldots, \ n.\end{equation} and {\it strictly row diagonally dominant} if Inequalities \eqref{ROW} are strict.
 A matrix $\mathbf A$ is called {\it (strictly) column diagonally dominant} if ${\mathbf A}^T$ is (strictly) row diagonally dominant.
 \item[\rm 4.] {\it generalized diagonally dominant} if there exist positive scalars (weights) $m_i$, $i = 1, \ \ldots, \ n$ such that
$$m_i|a_{ii}|>\sum_{j\neq i }m_j|a_{ij}|, \qquad i = 1, \ \ldots, \ n. $$
If, in addition, $a_{ii} < 0$ ($a_{ii} > 0$), $i = 1, \ \ldots, \ n,$ then ${\mathbf A}$ is called {\it negative diagonally dominant (NDD)} (respectively, {\it positive diagonally dominant (PDD)}). In literature, PDD matrices sometimes are called {\it diagonally quasidominant} (see, for example, \cite{MOY}, \cite{BHK3}).
\item[\rm 5.] {\it totally positive (TP)} if all of its minors up to the order $n$ are nonnegative and {\it strictly totally positive (STP)} if all of its minors up to the order $n$ are positive (see \cite{GANT}).
\item[\rm 6.] {\it oscillatory} if it is totally positive and there is a positive integer $p$ such that ${\mathbf A}^p$ is strictly totally positive (see \cite{GANT}).
\item[\rm 7.] a {\it $Z$-matrix} if all the off-diagonal entries $a_{ij}$, $i \neq j$, are nonpositive.
\item[\rm 8.] a {\it Metzler matrix} if all the off-diagonal entries $a_{ij}$, $i \neq j$, are nonnegative (i.e. if $-{\mathbf A}$ is a $Z$-matrix).
\item[\rm 9.] a {\it $P$-matrix ($P_0$-matrix)} if all its principal minors are positive (respectively, nonnegative), i.e the inequality $A \left(\begin{array}{ccc}i_1 & \ldots & i_k \\ i_1 & \ldots & i_k \end{array}\right) > 0$ (respectively, $\geq 0$)
holds for all $(i_1, \ \ldots, \ i_k), \ 1 \leq i_1 < \ldots < i_k \leq n$, and all $k, \ 1 \leq k \leq n$.
\item[\rm 10.] a {\it $Q$-matrix ($Q_0$-matrix)} if the inequality
$$\sum_{(i_1, \ldots, i_k)}A \begin{pmatrix}i_1 & \ldots & i_k \\ i_1 & \ldots & i_k \end{pmatrix} > 0 \quad (\mbox{respectively,} \ \geq 0)$$
holds for all $k, \ 1 \leq k \leq n$.
\item[\rm 11.] a {\it $P_0^+$-matrix)} if it is a $P_0$-matrix and, in addition, the sums of all principal minors of every fixed order $i$ are positive $(i = 1, \ \ldots, \ n)$ (i.e. if it is a $P_0$-matrix and a $Q$-matrix at the same time).
\item[\rm 12.] a {\it Hicksian matrix} if $-{\mathbf A}$ is a $P$-matrix.
\item[\rm 13.] an {\it almost Hicksian matrix} if $-{\mathbf A}$ is a $P_0^+$-matrix.
\item[\rm 14.] an {\it $M$-matrix} if its off-diagonal entries are nonpositive and the principal minors are all positive (i.e. if it is a $Z$-matrix and $P$-matrix at the same time).
\item[\rm 15.]  {\it sign-symmetric} if the inequality
\begin{equation}\label{SYM} A\begin{pmatrix}
  i_1 & \ldots & i_k \\
  j_1 & \ldots & j_k
\end{pmatrix}A\begin{pmatrix}
  j_1 & \ldots & j_k \\
  i_1 & \ldots & i_k
\end{pmatrix} \geq 0 \end{equation}
holds for all sets of indices $(i_1, \ \ldots, \ i_k), \ (j_1, \ \ldots, \ j_k)$, where $1 \leq i_1 < \ldots < i_k \leq n $, $1 \leq j_1 < \ldots < j_k \leq n $ (see \cite{CARL2}).
A matrix $\mathbf A$ is called {\it anti-sign-symmetric} if Inequalities \eqref{SYM} all has the opposite sign $\leq 0$ (see \cite{HERK2}).
\item[\rm 16.]  {\it weakly sign-symmetric} if the inequality
\begin{equation}\label{GKK} A\begin{pmatrix}
  \alpha \\
  \beta
\end{pmatrix}A\begin{pmatrix}
  \beta \\
  \alpha
\end{pmatrix} \geq 0 \end{equation}
holds for all sets of indices $\alpha =(i_1, \ \ldots, \ i_k), \ \beta=(j_1, \ \ldots, \ j_k)$, where $1 \leq i_1 < \ldots < i_k \leq n $, $1 \leq j_1 < \ldots < j_k \leq n $ and $|\alpha|=|\beta|= |\alpha\cap\beta| + 1$), i.e. the products of symmetrically located with respect to the principal diagonal almost principal minors are all nonnegative (see \cite{HOLS}).
\item[\rm 17.]  a {\it GKK matrix} (after Gantmacher--Krein--Kotelyansky) if it is a weakly sign-symmetric $P$-matrix. $\mathbf A$ is called a {\it strictly GKK matrix} if, in addition, Inequalities \eqref{GKK} are all strict.
\item[\rm 18.] {\it tridiagonal} if $a_{ij} = 0$ whenever $|i-j| > 1$.
\item[\rm 19.] {\it normal} if it commutes with its transpose: ${\mathbf A}{\mathbf A}^T = {\mathbf A}^T{\mathbf A}$.
\item[\rm 20.] {\it right (left) stochastic}, if it is (entry-wise) nonnegative and $\sum_i a_{ij} =1$ (respectively, $\sum_j a_{ij} =1$). A matrix ${\mathbf A}$ is called {\it doubly stochastic matrix} if it is right and left stochastic.
\item[\rm 21.]  {\it strictly row square diagonally dominant for every order of minors} if the following inequalities hold:
$$\left(A\left(\begin{array}{c} \alpha \\
\alpha \end{array}\right)\right)^2 > \sum_{\alpha,\beta \in [n], \alpha \neq \beta}\left( A\left(\begin{array}{c} \alpha \\
\beta \end{array}\right)\right)^2$$
for any $\alpha = (i_1, \ \ldots, \ i_k)$, $\beta = (j_1, \ \ldots, \ j_k)$ and all $k = 1, \ \ldots, \ n$.
A matrix $\mathbf A$ is called {\it strictly column square diagonally dominant} if ${\mathbf A}^T$ is strictly row square diagonally dominant.
\item[\rm 22.] A matrix ${\mathbf A}$ is called {\it a Kotelyansky matrix (K-matrix)} if all its principal minors are positive and all its almost principal minors are nonnegative. A matrix ${\mathbf A}$ is called a {\it strictly Kotelyansky matrix (SK-matrix)} if all its principal and almost principal minors are positive (see \cite{BAO}).
\item[\rm 23.] {\it reducible} if there exists a permutation of the indices that puts $\mathbf A$ into a block-triangular form, i.e.
 $${\mathbf P}{\mathbf A}{\mathbf P}^T = \widetilde{{\mathbf A}}, \qquad \widetilde{{\mathbf A}} = \begin{pmatrix}{\mathbf A}_{11} & {\mathbf O} \\ {\mathbf A}_{21} & {\mathbf A}_{22} \\ \end{pmatrix}$$
 where $\mathbf P$ is a permutation matrix, ${\mathbf A}_{11}$ and ${\mathbf A}_{22}$ are square matrices (see \cite{GANT}).
 \item[\rm 24.]{\it irreducible}, if it is not reducible.
\item[\rm 25.] {\it an $H$-matrix} if its comparison matrix is an $M$-matrix. A matrix $\mathbf A$ is called {\it an $H^+$-matrix} if, in addition, $a_{ii} \geq 0$ $(i = 1, \ \ldots, \ n)$. An $H$-matrix ${\mathbf A} = \{a_{ij}\}_{i,j = 1}^n$ is called an $H_+$-matrix if all its principal diagonal entries are nonnegative ($a_{ii} \geq 0, \ i = 1, \ \ldots, \ n$).
\end{enumerate}


\begin{thebibliography}{100}

\bibitem{AB1}
E.H. Abed, {\it Strong $D$-stability,} Systems Control Lett. {\bf 7} (1986), 207-212.

\bibitem{AB3}
E.H. Abed, {\it Multiparameter singular perturbation problems: Iterative expansions and asymptotic stability,} Systems Control Lett. {\bf 5} (1985), 279-282.

\bibitem{AB4}
E.H. Abed, {\it Decomposition and stability of multiparameter singular perturbation problems,} IEEE Transactions on Automatic Control {\bf AC-31} (1986), 925-934.

\bibitem{AB2}
E.H. Abed, {\it Singularly perturbed Hopf bifurcation,} IEEE Trans. Circuits and Systems {\bf 32} (1985), 1270-1280.

\bibitem{AB3}
E.H. Abed, S.P. Boyd, {\it Perturbation bounds for structured robust stability,} Proc. 27th IEEE Conf. Dec. Contr. (1988), 1029-1031.

\bibitem{ABT}
E.H. Abed, A.L. Tits, {\it On the stability of multiple time-scale systems,} International Journal of Control {\bf 44} (1986), 211-218.

\bibitem{ASCM}
M. Allouche, M. Souissi, M. Chaabane, D. Mehdi, {\it Robust $D$-stability analysis of an induction motor,} in 16th Mediterranean Conference on Control and Automation, (2008), 255-260.

\bibitem{ALT}
C. Altafini, {\it Stability analysis of diagonally equipotent matrices,} Automatica {\bf 49} (2013), 2780-2785.

\bibitem{ALAL}
Yu.A. Al'pin, V.S. Al'pina, {\it Combinatorial structure of $k$-semiprimitive matrix families,} Sbornik: Mathematics {\bf 207} (2016), pp. 639--651.

\bibitem{ANB}
B. Anderson, N. Bose, E. Jury, {\it A simple test for zeros of a complex polynomial in a sector,} IEEE Transactions on Automatic Control {\bf 19} (1974), 437--438.

\bibitem{AR}
M. Araki, {\it Application of $M$-matrices to the stability problems of composite dynamical systems,} Journal of Math. Analysis Appl. {\bf 52} (1975), 309-321.

\bibitem{ARC}
M. Arcak, {\it Linear matrix inequality tests for synchrony of diffusively coupled nonlinear systems,} Proceedings of the 58th Allerton Conference (2010), 1651--1656.

\bibitem{ARC2}
M. Arcak, {\it Diagonal stability on cactus graphs and application to network stability analysis,} IEEE Transactions on Automatic Control {\bf 56} (2011), 2766--2777.

\bibitem{ARCS1}
M. Arcak, E.D. Sontag, {\it Diagonal stability of a class of cyclic systems and its connection with the secant criterion,} Automatica {\bf 42} (2006), 1531--1537.

\bibitem{ARCS2}
M. Arcak, E.D. Sontag, {\it Connections between diagonal stability and the secant condition for cyclic systems,} Proceedings of the 2006 American Control Conference (2006), 1493--1498.

\bibitem{ART}
 A. Arhangel'skii, M. Tkachenko, {\it Topological groups and related structures,} Atlantis
Press, 2008.

\bibitem{AR1}
K.J. Arrow, M.D. Intriligator (Eds.) {\it Handbook on mathematical economics,} Vol. I, II. North-Holland, Amsterdam, 1987

\bibitem{AR2}
K.J. Arrow, S. Karlin, P. Suppes (Eds.), {\it Mathematical methods in the social sciences,} Stanford. Univ. Press, Stanford, 1960.

\bibitem{AM1}
K.J. Arrow, M. McManus, {\it A theorem on expectations and the stability of equilibrium,} Econometrica {\bf 24} (1956), 288-293.

\bibitem{AM}
K.J. Arrow, M. McManus, {\it A note on dynamical stability,} Econometrica {\bf 26} (1958), 448-454.

\bibitem{BACR}
A. Bacciotti, L. Rosier, {\it Liapunov functions and stability in control theory,} 2nd edition, Springer-Verlag, 2005.

\bibitem{BAP}
O. Bachelier, B. Pradin, {\it $\partial\mathcal{D}$-regularity for robust matrix root clustering}, IFAC Proceedings Volumes, {\bf 36} (2003), pp. 237-242.

\bibitem{BBM}
O. Bachelier, J. Bosche, D. Mehdi, {\it On matrix root-clustering in a combination of first order regions}, IFAC Proceedings Volumes, {\bf 39} (2006), pp. 405-410.

\bibitem{BAM}
O. Bachelier, D. Mehdi, {\it Robust matrix ${\mathcal D}_U$-stability analysis}, International Journal of Robust and Nonlinear Control, {\bf 13} (2003), pp. 533-558.

\bibitem{BHPM}
O. Bachelier, D. Henrion, B. Pradin, D. Mehdi, {\it Robust root-clustering of a matrix in intersections or unions of regions}, SIAM J. Control Optim., {\bf 43} (2004), pp. 1078--1093.

\bibitem{BC}
C. Bahl, B. Cain, {\it The inertia of diagonal multiples of $3\times 3$ real matrices}, Linear Algebra Appl., {\bf 18} (1997), pp. 267--280.

\bibitem{BAL}
C.S. Ballantine, {\it Stabilization by a diagonal matrix}, Proc. Amer. Math. Soc., {\bf 25}
(1970), pp. 728--734.

\bibitem{BAR}
G.P. Barker, {\it Common solutions to the Lyapunov equations}, Linear Algebra Appl., {\bf 16} (1977), pp. 233--235.

\bibitem{BARBERPL}
G.P. Barker, A. Berman, R.J. Plemmons, {\it Positive diagonal solutions to the Lyapunov equations}, Linear and Multilinear Algebra \textbf{5} (1978), 249--256.

\bibitem{BARK}
Y.S. Barkovsky, {\it Rank-one perturbations method and differential operators of oscillatory type} (in Russian),
PhD Thesis, Rostov-on-Don, 1980.

\bibitem{BAO}
Y.S. Barkovsky, T.V. Ogorodnikova, {\it On matrices with positive and simple spectra}, Izvestiya SKNC VSH Natural sciences \textbf{4} (1987), 65--70.

\bibitem{BAY}
Y.S. Barkovsky, V.I. Yudovich, {\it The momenta problem and spectral theory of the operators}, Izvestiya SKNC VSH Natural sciences \textbf{4} (1975), 49--53.

\bibitem{BAR}
B.R. Barmish, {\it New tools for robustness of linear systems}, Macmillan, New York, 1994.

\bibitem{BED}
B. Beavis, I.M. Dobbs, {\it Optimization and stability theory for economic analysis}, Cambridge University Press, 1990.

\bibitem{BELL}
R. Bellman, {\it Introduction to matrix analysis}, McGraw Hill, New York, 2nd edition, 1970.

\bibitem{BERH}
A. Berman, D. Hershkowitz, {\it Matrix diagonal stability and its implications,} SIAM J. Alg. Disc. Meth., {\bf 4} (1983), pp. 377-382.

\bibitem{BERH1}
A. Berman, D. Hershkowitz, {\it Graph theoretical methods in studying stability}, Contemporary Math., {\bf 47} (1985), pp. 1--6.

\bibitem{BERH2}
A. Berman, D. Hershkowitz, {\it Characterization of acyclic $D$-stable matrices,} Linear Algebra Appl., {\bf 58} (1984), pp. 17-31.

\bibitem{BERPL}
A. Berman, R.J. Plemmons, {\it Nonnegative Matrices in
the Mathematical Sciences}, Academic Press, New York, 1979.

\bibitem{BERW}
A. Berman, R.C. Ward, {\it Classes of stable and semipositive matrices,} Linear Algebra Appl., {\bf 21} (1978), pp. 163-174.

\bibitem{BERVW}
A. Berman, R.S. Varga, R.C. Ward, {\it Matrices with nonpositive off-diagonal entries,} Linear Algebra Appl., {\bf 21} (1978), pp. 233-244.

\bibitem{BEF}
B. Besselink, H.R. Feyzmahdavian, H. Sandberg, M. Johansson, {\it $D$-stability and delay-independent stability of monotone nonlinear systems with max-separable Lyapunov functions}, IEEE Conference on Decision and Control (CDC),
(2016), pp. 3172--3177.

\bibitem{BHK}
A. Bhaya, E. Kaszkurewicz, {\it On discrete-time diagonal and $D$-stability}, Linear Algebra Appl., {\bf 187} (1993), pp. 87--104.

\bibitem{BHK2}
A. Bhaya, E. Kaszkurewicz, {\it Control perspectives on numerical algorithms and matrix problems}, SIAM, 2006.

\bibitem{BHK3}
A. Bhaya, E. Kaszkurewicz, R. Santos {\it Characterizations of classes of stable matrices}, Linear Algebra Appl., {\bf 374} (2003), pp. 159--174.

\bibitem{BHAT2}
R. Bhatia, {\it Positive definite matrices}, Princeton University Press, 2007.

\bibitem{BH}
S. Bhattacharya, H. Chapellat, L. Keel, {\it Robust control: the parametric approach}, Prentice-Hall, New Jersey, 1995.

\bibitem{BICJ1}
T.A. Bickart, E.I. Jury, {\it Regions of polynomial root clustering,} Journal of the Franklin Institute, {\bf 304} (1977), pp. 149-160.

\bibitem{BICJ2}
T.A. Bickart, E.I. Jury, {\it The Schwarz--Christoffel transformation and polynomial root clustering,} IFAC Proceedings Volumes {\bf 11} (1978), pp. 1171-1176.

\bibitem{BICJ3}
T.A. Bickart, E.I. Jury, {\it Polynomial root clustering,} Journal of the Franklin Institute, {\bf 308} (1979), pp. 487-496.

\bibitem{BIR}
J. Bierkens, A. Ran, {\it A singular $M$-matrix perturbed by a nonnegative rank-one matrix has positive principal minors; is it $D$-stable?} Linear Algebra Appl., {\bf 457} (2014), pp. 191--208.

\bibitem{BFG}
F. Blanchini, E. Franco, G. Giordano, {\it Determining the structural properties of a class of biological models}, 51th IEEE Conference on Decision and Control (2012), pp. 5505--5510.

\bibitem{BGF}
S. Boyd, L. El Ghaoui, E. Feron, V. Balakrishnan, {\it Linear matrix inequalities in system and control theory}, SIAM, 1994.

\bibitem{BOR}
P. Borwein, T. Erdelyi, {\it Polynomials and polynomial inequalities}, Springer, 1995.

\bibitem{BURL1}
L.A. Burlakova, {\it $D$-stable 4th-order matrices,} J. Sovremennie technologii. Systemniy analiz. Modelirovanie {\bf 1 (21)} (2009), 109-116.

\bibitem{BURL2}
L.A. Burlakova, {\it Conditions of $D$-stability of the fifth-order matrices,} in: Gerdt V.P., Mayr E.W., Vorozhtsov E.V. (eds) Computer Algebra in Scientific Computing. Lecture Notes in Computer Science, Springer, Berlin, Heidelberg {\bf 5743} (2009), pp. 54-65.

\bibitem{BUT}
P. Butkovi\u{c}, {\it Max-linear systems: theory and algorithms,} Springer-Verlag London Limited 2010.

\bibitem{CA}
B. Cain, {\it Real, $3 \times 3$, $D$-stable matrices,} J. Res. Nat. Bur. Standards Sect. B, {\bf 80B} (1976), 75–77.

\bibitem{CA1}
B. Cain, {\it Inside the $D$-stable matrices,} Linear Algebra Appl., {\bf 56} (1984), 237–243.

\bibitem{CA2}
B. Cain, {\it Convergent multiples of convergent operators,} Linear Algebra Appl., {\bf 299} (1999), 171–173.

\bibitem{CADHJ}
B. Cain, L.M. DeAlba, L. Hogben, C.R. Johnson, {\it Multiplicative perturbations of stable and convergent operators,} Linear Algebra Appl., {\bf 268} (1998), pp. 151--169.

\bibitem{CALN}
B. Cain, T.D. Lenker, S.K. Narayan, P. Vermeire {\it Classes of stable complex matrices defined via the theorems of Ger\u{s}gorin and Lyapunov,} Linear and Multilinear Algebra, {\bf 56} (2008), pp. 713--724.

\bibitem{CAMM}
P.J. Campo, M. Morari, {\it Achievable closed-loop properties of systems under decentralized control: conditions involving the steady-state gain}, IEEE Transactions on Automatic Control, {\bf 39} (1994), pp. 932--943.

\bibitem{CARL1}
D. Carlson, {\it Controllability, inertia and stability for tridiagonal matrices}, Linear Algebra Appl., {\bf 56} (1984), pp. 207--220.

\bibitem{CARL2}
D. Carlson, {\it A class of positive stable matrices}, J. Res. Nat. Bur. Standards Sect. B, {\bf 78B} (1974), pp. 1--2.

\bibitem{CARL3}
D. Carlson, {\it A new criterion for $H$-stability of complex matrices}, Linear Algebra Appl., {\bf 1} (1968), pp. 59--64.

\bibitem{CARD}
D. Carlson, B. Datta, {\it The Lyapunov matrix equation ${\mathbf S}{\mathbf A} + {\mathbf A}^*{\mathbf S}={\mathbf S}^*{\mathbf B}^*{\mathbf B}{\mathbf S}$}, Linear Algebra Appl., {\bf 28} (1979), pp. 43--52.

\bibitem{CARDJ}
D. Carlson, B. Datta, C. Johnson, {\it A semi-definite Lyapunov theorem and the characterization of tridiagonal $D$-stable matrices}, SIAM J. Alg. Disc. Meth., {\bf 3} (1982), pp. 293--304.

\bibitem{CAHS}
D. Carlson, D. Hershkowitz, D. Shasha, {\it Block diagonal semistability factors and Lyapunov semistability of block triangular matrices}, Linear Algebra Appl., {\bf 172} (1992), pp. 1--25.

\bibitem{CAS}
D. Carlson, H. Schneider, {\it Inertia theorems for matrices: the semidefinite case,} Journal of Math. Analysis and Appl. {\bf 6} (1963), 430--446.

\bibitem{CAU}
A.L. Cauchy, {\it Calcul des indices des fonctions,} J. \'{E}cole Polytech. {\bf 15}, 176–229 (1837) ({\OE}uvres {\bf 1}(2),
416–466).

\bibitem{CHEN}
C.-T. Chen, {\it A generalization of the inertia theorem}, SIAM J. Appl. Math., {\bf 25} (1973), pp. 158--161.

\bibitem{CHEN2}
C.-T. Chen, {\it Linear system theory and design}, 3rd Ed. Oxford University Press, 1999.

\bibitem{CFY}
J. Chen, M. Fan, Ch.-Ch. Yu, {\it On $D$-stability and structured singular values}, Systems and Control letters, {\bf 24} (1995), pp. 19--24.

\bibitem{CHIGA}
M. Chilali, P. Gahinet, {\it $H_{\infty}$ design with pole placement constraints: an LMI approach},
IEEE Transactions on Automatic Control, {\bf 41} (1996), pp. 358--367.

\bibitem{CGA}
M. Chilali, P. Gahinet, P. Apkarian, {\it Robust pole placement in LMI regions},
Proceedings of the 36th Conference on Decision and Control
San Diego, USA, 1997, pp. 1291--1296.

\bibitem{CHO}
D. Choi, {\it Inequalities related to partial trace and block Hadamard product}, Linear and Multilinear Algebra, {\bf 66} (2018), pp. 280-284.

\bibitem{CHU}
T. Chu, {\it An equivalent condition for stability properties of Lotka--Volterra systems}, Physics Letters {\bf A 368} (2007), pp. 235-237.

\bibitem{COH}
A. Cohn, {\it \"{U}ber die Anzahl der Wurzeln einer algebraischen Gleichung in einem Kreise}, Mathematische Zeitschrift, {\bf 14} (1922), pp. 110-148.

\bibitem{CROSS}
G.W. Cross, {\it Three types of matrix stability}, Linear Algebra Appl., {\bf 20} (1978), pp. 253--263.

\bibitem{CURT}
M.L. Curtis, {\it Matrix groups}, Springer-Verlag, 1984.

\bibitem{DARW}
S. Dashkovskiy, B.S. R\"{u}ffer, F.R. Wirth, {\it An ISS small gain theorem for general networks}, Math. Control Signals Syst., {\bf 19} (2007), pp. 93--122.

\bibitem{DATTA}
B.N. Datta, {\it Stability and $D$-stability}, Linear Algebra Appl., {\bf 21} (1978), pp. 135--141.

\bibitem{DATTA1}
B.N. Datta, {\it Stability and inertia}, Linear Algebra Appl., {\bf 302-303} (1999), pp. 563--600.

\bibitem{DEC}
R. Descartes, {\it La G\'{e}om\'{e}trie}, Leyden France: Maire, 1637. (Translation: {\it The Geometry of Ren\'{e} Descartes}. La Salle, France: Open Court, 1925.)

\bibitem{DOP}
M. Domijan, E. P\'{e}cou, {\it The interaction graph structure of mass-action reaction networks,} Mathematical Biology, {\bf 65} (2012), pp. 375--402.

\bibitem{DJD}
J. Drew, C. Johnson, P. van den Driessche, {\it Strong forms of nonsingularity}, Lin. Algebra Appl., {\bf 162-164} (1992), pp. 187--204.

\bibitem{DUP}
G.-R. Duan, R.J. Patton, {\it A note on Hurwitz stability of matrices}, Automatica, {\bf 34} (1998), pp. 509--511.

\bibitem{ENT}
A.C. Enthoven, K.J. Arrow, {\it A theorem on expectations and the stability of equilibrium}, Econometrica, {\bf 24} (1956), pp. 288--293.

\bibitem{FAJ}
S.M. Fallat, C.R. Johnson, {\it Sub-direct sums and positivity classes of matrices}, Linear Algebra Appl., {\bf 288} (1999), pp. 149--173.

\bibitem{FER}
A. S. R. Ferreira, {\it Exploiting Structure and Input-Output Properties in Networked Dynamical Systems}, Technical Report No. UCB/EECS-2015-245 (2015)
http://www.eecs.berkeley.edu/Pubs/TechRpts/2015/EECS-2015-245.html

\bibitem{FID}
 M. Fiedler, {\it Additive compound matrices and an inequality for eigenvalues of symmetric stochastic matrices}, Czech. Math. J., 24
(1974), pp. 392--402.

\bibitem{FIF}
M.E. Fisher, A.T. Fuller, {\it On the stabilization of matrices and the convergence of linear iterative processes}, Proc. Cambridge Philos. Soc., {\bf 54} (1958), pp. 417--425.

\bibitem{FISH}
F.M. Fisher, {\it A simpler proof of the Fisher--Fuller theorem}, Proc. Cambridge Philos. Soc., {\bf 71} (1972), pp. 523--525.

\bibitem{FL1}
 R. Fleming, G. Grossman, T. Lenker, S. Narayan, S.-C. Ong, {\it Classes of Schur $D$-stable matrices,} Linear Algebra Appl., {\bf 306} (2000), pp. 15--24.

 \bibitem{FL2}
 R. Fleming, G. Grossman, T. Lenker, S. Narayan, S.-C. Ong, {\it On Schur $D$-stable matrices,} Linear Algebra Appl., {\bf 279} (1998), pp. 39--50.

\bibitem{FGR}
J.H. Fourie, G.J. Groenewald, D.B. Janse van Rensburg, A.C.M. Ran, {\it Rank-one perturbations of $H$-positive real matrices,} Linear Algebra Appl., {\bf 439} (2013), pp. 653--674.

\bibitem{GAQ}
Z. Gajic, M. Qureshi, {\it The Lyapunov matrix equation in system stability and control,} Academic Press, Inc, 1995.

\bibitem{GANT}
 F. Gantmacher, {\it The Theory of Matrices,} Volume 1,
Volume 2. Chelsea. Publ. New York, 1990.

\bibitem{GANT2}
F. Gantmacher, {\it Applications of the Theory of Matrices,} Dover Publications, 2005.

\bibitem{GEA}
 X. Ge, M. Arcak, {\it A sufficient condition for additive $D$-stability and application to reaction-diffusion models,} Systems and Control Letters, {\bf 58} (2009), pp. 736--741.

\bibitem{GER}
 J. C. Geromel, {\it On the determination of a diagonal solution of the Lyapunov equation,} IEEE Transactions on Automatic Control, {\bf AC-30} (1985), pp. 404--406.

\bibitem{GEOH}
 J. C. Geromel, M.C. de Oliveira, L. Hsu, {\it LMI characterization of structural and robust stability,} Linear Algebra Appl., {\bf 285} (1998), pp. 69--80.

\bibitem{GIZ}
 G. Giorgi, C. Zuccotti, {\it An overview on $D$-stable matrices,} Universit\`{a} di Pavia, Department of Economics and Management, DEM Working Paper Series {\bf 97} (2015), pp. 1--28.

 \bibitem{GOH1}
 B.H. Goh, {\it Global stability in two species interactions,} Journal of Mathematical Biology, {\bf 3} (1976), pp. 313--318.

\bibitem{GUH}
M. Gumus, J. Xu, {\it A new characterization of simultaneous Lyapunov diagonal stability via Hadamard products,} Linear Algebra Appl., {\bf 531} (2017), pp. 220--233.

\bibitem{GUH1}
M. Gumus, J. Xu, {\it Some new results related to $\alpha$-stability,} Linear and Multilinear Algebra, {\bf 65} (2017), pp. 325--340.

\bibitem{GUT}
S. Gutman, {\it Matrix root clustering in algebraic regions,} International Journal of Control, {\bf 39} (1984), pp. 773--778.

\bibitem{GUT2}
S. Gutman, {\it Root clustering in parameter space,} Springer-Verlag Berlin, Heidelberg, 1990.

\bibitem{GUJU}
S. Gutman, E. Jury, {\it A general theory for matrix root-clustering in subregions of the complex plane,} IEEE Transactions on Automatic control, {\bf AC-26} (1981), pp. 853--863.

\bibitem{HAD}
 K.P. Hadeler, {\it Nonlinear diffusion equations in biology,} in Proceedings of the Conference on Differential Equations, Dundee 1976, Springer Lecture Notes.

\bibitem{HART}
 D.J. Hartfiel, {\it Concerning the interior of the $D$-stable matrices,} Linear Algebra Appl., {\bf 30} (1980), pp. 201--207.

\bibitem{HENBS}
 D. Henrion, O. Bachelier, M. \v{S}ebek, {\it $D$-stability of polynomial matrices,} International Journal of Control, {\bf 74} (2001), pp. 845--856.

\bibitem{HENG}
 D. Henrion, A. Garulli (eds), {\it Positive polynomials in control,} Springer, 2005.

\bibitem{HERMI}
C. Hermite, {\it “On the number of roots of an algebraic equation between two limits,” Extract of a letter from Mr. C. Hermite of Paris to Mr. Borchardt of Berlin}, J. Reine angew. Math., {\bf 52} (1856), pp. 39–51. Translation by P.C. Parks, Int. J. Cont. {\bf 26} (1977), pp. 183–196.

\bibitem{HER1}
 D. Hershkowitz, {\it Recent directions in matrix stability,} Linear Algebra Appl., {\bf 171} (1992), pp. 161--186.

\bibitem{HERB}
 D. Hershkowitz, A. Berman, {\it Localization of the spectra of $P$- and $P_0$-matrices,} Linear Algebra Appl., {\bf 52/53} (1983), pp. 383--397.

\bibitem{HERJ1}
 D. Hershkowitz, C.R. Johnson, {\it Spectra of matrices with $P$-matrix powers,} Linear Algebra Appl., {\bf 80} (1986), pp. 159--171.

\bibitem{HERK2}
D. Hershkowitz and N. Keller, {\it Positivity of principal minors, sign symmetry and stability}, Linear Algebra Appl., {\bf 364} (2003), pp. 105--124.

\bibitem{HERM}
D. Hershkowitz, N. Mashal, {\it $P^\alpha$-matrices and Lyapunov scalar stability}, ELA. {\bf 4} (1998), 39-47.

\bibitem{HERSS1}
D. Hershkowitz, H. Shneider, {\it Lyapunov diagonal semistability of real $H$-matrices}, Linear Algebra Appl., {\bf 71} (1985), pp. 119--145.

\bibitem{HERSS3}
D. Hershkowitz, H. Shneider, {\it Scalings of vector spaces and the uniqueness of Lyapunov scaling factors}, Linear and Multilinear Algebra, {\bf 17} (1985), pp. 203--226.

\bibitem{HERSS2}
D. Hershkowitz, H. Schneider, {\it On Lyapunov scaling factors of real symmetric matrices}, Linear and Multilinear Algebra, {\bf 22} (1988), pp. 373--384.

\bibitem{HIL}
R.D. Hill, {\it Inertia theory for simultaneously triangulable complex matrices,} Linear Algebra Appl. {\bf 2} (1969), 131-142.

\bibitem{HOLS}
O. Holtz, H. Schneider, {\it Open problems on GKK $\tau$-matrices,} Linear Algebra Appl. {\bf 345} (2002), 263-267.

\bibitem{HOJ}
 R. Horn, C.R. Johnson, {\it Topics in matrix analysis,} Cambridge University
Press, 1991.

\bibitem{HMN}
 R. Horn, R. Mathias, Y. Nakamura, {\it Inequalities for unitarily invariant norms and bilinear matrix products,} Linear and Multilinear Algebra, {\bf 30} (1991), pp. 303--314.

\bibitem{HHP}
F.-H. Hsiao, J.-D. Hwang, S.-P. Pan, {\it $D$-stability analysis for discrete uncertain time-delay systems,} Appl. Math. Lett., {\bf 11} (1998), pp. 109--114.

\bibitem{HU}
H. Hu, {\it An algorithm for rescalling a matrix positive definite,} Linear Algebra Appl., {\bf 96} (1987), pp. 131--147.

\bibitem{HUR}
A. Hurwitz, {\it \"{U}ber die Begingungen, unter welchen eine Gleichung nur Wurzeln mit negativoen reelen Teilen besitzt,} Math. Ann., {\bf 46} (1895), pp. 273--284 (Werke, {\bf 2}, pp. 533--545).

\bibitem{IJP}
S.T. Impram, R. Johnson, R. Pavani, {\it The $D$-stability problem for $4 \times 4$ real matrices,} Archivum Mathematicum (BRNO) {\bf 41} (2005), 439--450.

\bibitem{JEM}
R. Jeltsch, M. Mansour (eds), {\it Stability theory,} Hurwitz Centenary Conference, Ascona, 1995. Birkh\"{a}user Verlag, 1996.

\bibitem{JOHN6}
C.R. Johnson, {\it Positive definite matrices,} The American Mathematical Monthly {\bf 77} (1970), 259-264.

\bibitem{JOHNN}
C.R. Johnson, {\it $D$-stability and real and complex quadratic forms,} Linear Algebra Appl. {\bf 9} (1974), 89-94.

\bibitem{JOHN1}
C.R. Johnson, {\it Sufficient conditions for $D$-stability,} Journal of Economic Theory {\bf 9} (1974), 53-62.

\bibitem{JOHN3}
C.R. Johnson, {\it Hadamard products of matrices,} Linear and Multilinear Algebra {\bf 1:4} (1974), 295-307.

\bibitem{JOHN4}
C.R. Johnson, {\it Second, third and fourth order $D$-stability,} J. Research Nat. Bureau Standards USA {\bf B78(1)} (1974), 11-13.

\bibitem{JOHN5}
C.R. Johnson, {\it A characterization of the nonlinearity of $D$-stability,} Journal of Mathematical Economics {\bf 2} (1975), 87-91.

\bibitem{JOHD}
C.R. Johnson, P. van den Driessche, {\it Interpolation of $D$-stability and sign stability,} Linear and Multilinear Algebra, {\bf 23} (1988), 363-368.

\bibitem{JOHN2}
C.R. Johnson, S. Narayan, {\it When the positivity of the leading principal minors implies the positivity of all principal minors of a matrix,} Linear Algebra Appl., {\bf 439} (2013), 2934-2947.

\bibitem{JOHNT}
R. Johnson, A. Tesi, {\it On the $D$-stability problem for real matrices,} Bollettino dell'Unione Matematica Italiana, {\bf 2-B} (1999), 299--314.

\bibitem{JONCK}
E.A. Jonckheere, {\it Algebraic and differential topology of robust stability,} Oxford University Press, 1997.

\bibitem{JU2}
 E.I. Jury, {\it Inners and stability of dynamic systems,} 2nd edition, Florida: R.E. Krieger, 1982.

\bibitem{JU3}
 E.I. Jury, {\it Stability, root clustering and inners,} IFAC Proceedings, {\bf 5} (1972), 153--159.

\bibitem{JUA}
E.I. Jury, S.M. Ahn, {\it Symmetric and innerwise matrices for the root-clustering and root-distribution of a polynomial,} Journal of the Franklin Institute {\bf 293} (1972), pp. 433-450.


\bibitem{KAF}
W. Kafri, {\it Robust $D$-stability,} Applied Math. Letters {\bf 15} (2002), 7-10.

\bibitem{KAL}
R.E. Kalman, {\it On the Hermite--Fujiwara theorem in stability theory,} Quart. Appl. Math. (1965).

\bibitem{KANO}
G.V. Kanovei, {\it On one necessary condition of the $D$-stability of matrices having no less than two zero elements on the main diagonal}, Automat. Remote Control, {\bf 62} (2001), pp. 704--708.

\bibitem{KL}
G.V. Kanovei and D.O. Logofet, {\it $D$-stability of 4-by-4 matrices}, Comput. Math. Math. Phys., {\bf 38} (1998), pp. 1369--1374.

\bibitem{KL2}
G.V. Kanovei and D.O. Logofet, {\it Relations, properties and invariant transformations of $D$- and $aD$-stable matrices}, Vestnik Moskov. Univ. Ser I Mat. Mekh., {\bf 6} (2001), pp. 40--43.

\bibitem{KAB}
E. Kaszkurewicz, A. Bhaya, {\it Matrix diagonal stability in systems and computation}, Springer, 2000.

\bibitem{KAB1}
E. Kaszkurewicz, A. Bhaya, {\it Qualitative stability of discrete-time systems}, Linear Algebra Appl. {\bf 117} (1989), pp. 65-71.

\bibitem{KAB2}
E. Kaszkurewicz, A. Bhaya, {\it Robust stability and diagonal Liapunov functions}, SIAM J. Matrix Anal. Appl. {\bf 14} (1993), pp. 508-520.

\bibitem{KEL}
R.B. Kellogg, {\it On complex eigenvalues of $M$- and $P$-matrices,} Numer. Math. {\bf 19} (1972), 170-175.

\bibitem{KEK}
M.C. Kemp, Y. Kimura, {\it Introduction to Mathematical Economics,} Springer Verlag, New York, 1978.

\bibitem{KHAL2}
H.K. Khalil, {\it Asymptotic stability of nonlinear multiparameter singularly perturbed systems,} Automatica {\bf 17} (1981), pp. 797-804.

\bibitem{KHAL3}
H.K. Khalil, {\it A new test for $D$-stability,} Journal of Economic Theory {\bf 23} (1980), pp. 120-122.

\bibitem{KHAL4}
H.K. Khalil, {\it On the existence of positive diagonal $P$ such that $PA + A^TP <0 $,} IEEE Transactions on Automatic Control {\bf AC-27} (1982), pp. 181-184.

\bibitem{KHAK1}
H.K. Khalil, P.V. Kokotovic, {\it Control of linear systems with multiparameter singular perturbations,} Automatica {\bf 15} (1979), pp. 197-207.

\bibitem{KHAK2}
H.K. Khalil, P.V. Kokotovic, {\it $D$-stability and multi-parameter singular perturbation,} SIAM J. Control Optim. {\bf 17} (1979), pp. 56-65.

\bibitem{KHAR1}
V.L. Kharitonov, {\it The asymptotic stability of the equilibrium state of a family of systems of linear differential equations}, Differ. Uravn., {\bf 14} (1978), p. 2086--2088. (Russian)

\bibitem{KHAR2}
V.L. Kharitonov, {\it Stability of imbedded families of polynomials}, Automat. i Telemekh., {\bf 11} (1995), p. 169--178. (Russian)

\bibitem{KIMB}
K.-K. Kim, R.D. Braatz, {\it Continuous- and discrete-time $D$-stability, Joint $D$-stability and their applications: $\mu$ theory and diagonal stability approaches,}
 Proceedings of the 51st IEEE Conference on Decision and Control (2012), pp. 2896-2901.

\bibitem{KIM}
Y. Kimura, {\it A note on sufficient conditions for D-stability,}
Journal of Mathematical Economics, {\bf 8} (1981), pp. 113-120.

\bibitem{KINN}
C. King, M. Nathanson, {\it On the existence of a common quadratic Lyapunov function for a rank one difference,}
Linear Algebra Appl., {\bf 419} (2006), pp. 400-416.

\bibitem{KINS}
C. King, R. Shorten, {\it Singularity conditions for the non-existence of a common quadratic Lyapunov function for pairs of third order linear time-invariant dynamic systems,}
Linear Algebra Appl., {\bf 413} (2006), pp. 24-35.

\bibitem{KOG}
J. Kogan, {\it Robust stability and convexity,}
Springer--Verlag, 1995.

\bibitem{KOS1}
A. Kosov, {\it About a class of systems preserving the stability property at negative feedbacks,}
Automation and Remote Control, {\bf 69} (2008), pp. 764-773.

\bibitem{KOS}
A. Kosov, {\it On the D-stability and additive D-stability of matrices and Svicobians,}
Journal of Applied and Industrial Mathematics, {\bf 4} (2010), pp. 200-212.

\bibitem{KR}
J. Kraaijevanger, {\it A characterization of Lyapunov diagonal stability using Hadamard products,} Linear Algebra Appl., {\bf 151} (1991), pp. 245--254.

\bibitem{KU1}
O.Y. Kushel, {\it On a criterion of $D$-stabiity for $P$-matrices}, Special Matrices, {\bf 4} (2016), pp. 181-188.

\bibitem{KU2}
O.Y. Kushel, {\it Generalized diagonal and $D$-stability in LMI regions}, in preparation.

\bibitem{KUSHN}
H.J. Kushner, {\it Stochastic stability and control}, Academic Press, New York, 1967.

\bibitem{LLK}
C.-H. Lee, T.-H. Li, F.-C. Kung, {\it $D$-stability analysis for discrete systems with a time delay,} Systems and Control Letters {\bf 19} (1992), 213-219.

\bibitem{LEE}
J. Lee, T. Edgar, {\it Real structured singular value conditions for the strong $D$-stability,} Systems and Control Letters {\bf 44} (2001), 273-277.

\bibitem{LIW}
M.Y. Li, L. Wang, {\it A criterion for stability of matrices,} Journal of Mathematical Analysis and Applications {\bf 225} (1998), 249-264.

\bibitem{LIW}
X.-B. Liang, J. Wang, {\it An additive diagonal-stability condition for absolute exponential stability of a general class of neural networks,} IEEE Transactions on circuits and systems I: Fundamental Theory and Applications {\bf 48} (2001), 1308-1317.

\bibitem{LOC}
A. Locatelli, N. Schiavoni, {\it A necessary and sufficient condition for the stabilization of a matrix and its principal submatrices,} Linear Algebra Appl. {\bf 436} (2012), 2311-2314.

\bibitem{LOG2}
D.O. Logofet, {\it Svicobians of the compartment models and $DaD$-stability of the Svicobians: aggregating "$0$-dimensional" models of global biogeochemical cycles,} Ecological Modelling {\bf 104} (1997), 39-49.

\bibitem{LOG1}
D.O. Logofet, {\it Matrices and graphs stability problems in mathematical ecology,} 2nd edition, CRC Press, 2018.

\bibitem{LOG}
D.O. Logofet, {\it Stronger-than-Lyapunov notions of matrix stability, or how "flowers" help solve problems in mathematical ecology,} Linear Algebra Appl. {\bf 398} (2005), 75-100.

\bibitem{LYA}
M.A. Lyapunov, {\it General problem of the stability of motion}, Math. Soc. Kharkov (1892) (in Russian). Translation: Probleme generale de la stabilite du mouvement, Princeton Univ. Press, N.J. (1949).

\bibitem{MAO}
W.-J. Mao, {\it An LMI approach to $D$-stability and $D$-stabilization of linear discrete singular systems with state delay,} Appl. Math. and Computations {\bf 218} (2011), 1694--1704.

\bibitem{MAM}
M. Marcus, H. Minc, {\it A survey of matrix theory and matrix inequalities}, Allyn and Bacon. Inc., Boston, 1964.

\bibitem{MAR}
M. Marden, {\it Geometry of polynomials}, AMS, Providence, 1966.

\bibitem{MART}
A.A. Martynyuk, {\it Stability by Liapunov's matrix function methods with applications}, Marcel Dekker, Inc., NY, 1998.

\bibitem{MASS}
O. Mason, R. Shorten, {\it On the simultaneous diagonal stability of a pair of positive linear systems,}
Linear Algebra Appl., {\bf 413} (2006), pp. 13-23.

\bibitem{MAT}
D. Matignon, {\it Stability results for fractional differential equations with applications to control processing,}
Computational Engeneering in Systems Appl., {\bf 2} (1996), pp. 963--968.

\bibitem{MAT1}
D. Matignon, {\it Stability properties for generalized fractional differential systems,}
ESAIM Proceedings Fractional Differential Systems: Models, Methods and Applications, {\bf 5} (1998), pp. 145--158.

\bibitem{MELN}
V.G. Melnikov, {\it A sweeping method for matrix root clustering,} IFAC Proceedings Volumes {\bf 44} (2011), pp. 168-171.

\bibitem{MEME}
F. Mesquine, D. Mehdi, {\it Pole assignment in LMI Regions for linear constrained control systems,} in: Proceedings of the 15th Mediterranean conference on control and automation, Athens-Greece; July 27–29; 2007.

\bibitem{ME}
L. Metzler, {\it Stability of multiple markets: the Hicks conditions,} Econometrica {\bf 13} (1945), 277--292.

\bibitem{MOMOK}
T. Mori, Y. Mori, H. Kokame, {\it Common Lyapunov function approach to matrix root clustering,} System and Control Letters, {\bf 44} (2001), 73--78.

\bibitem{MOY}
P.J. Moylan, {\it Matrices with positive principal minors,} Linear Algebra Appl., {\bf 17} (1977), 53--58.


\bibitem{MOYH}
P. Moylan, D. Hill, {\it Stability criteria for large-scale systems,} IEEE Transactions on Automatic Control, {\bf AC-23} (1978), 143--149.

\bibitem{NAS}
K.S. Narendra, R. Shorten, {\it A characterization of the Hurwitz stability of Metzler matrices,} American Control Conference (2009), 1833--1837.

\bibitem{NIS}
H.J. Nieuwenhuis, L. Schoonbeek, {\it Stability of matrices with negative diagonal submatrices,} Linear Algebra Appl., {\bf 353} (2002), 183--196.

\bibitem{OLN}
N. Oleng, K. Narendra, {\it On the existence of diagonal solutions to the Lyapunov equations for third order systems,} Proceedings of American Control Conference, 2003.

\bibitem{OBG}
M.C. de Oliveira, J. Bernussou, J.C. Geromel, {\it A new discrete-time robust stability condition,} Systems and Control Letters {\bf 37} (1999), 261--265.

\bibitem{OGH}
M.C. de Oliveira, J.C. Geromel, L. Hsu {\it LMI characterization of structural and robust stability: the discrete-time case,} Linear Algebra Appl., {\bf 296} (1999), 27--38.

\bibitem{OLP}
R. Oliveira, P. Peres {\it A simple and less conservative test for $D$-stability,} SIAM J. Matrix Anal. Appl., {\bf 26} (2005), 415--425.

\bibitem{OSS}
A. Ostrowski, H. Schneider, {\it Some theorems on the inertia of general matrices,} Journal of Math. Analysis and Appl. {\bf 4} (1962), 72--84.

\bibitem{PAR1}
 P.C. Parks, {\it A new proof of the Routh--Hurwitz stability criterion using the second method of Lyapunov,} Proc. Comb. Philos.Soc., {\bf 58} (1962), pp. 694--702.

\bibitem{PAR2}
 P.C. Parks, {\it Lyapunov and Schur--Cohn stability criterion,} IEEE Trans. Automat. Contr., {\bf AC-9} (1964), p. 121.

\bibitem{PAR3}
P.C. Parks, {\it A new proof of Hermite's stability criterion and a generalization of Orlando's formula,} Int. J. Control, {\bf 26} (1977), pp. 197--206.

\bibitem{PAV1}
R. Pavani, {\it About characterization of $D$-stability by a computer algebra approach,} AIP Conf. Proc., {\bf 1558} (2013), pp. 309--312.

\bibitem{PAV2}
R. Pavani, {\it A new efficient approach to the characterization of $D$-stable matrices,} Mathematical Methods in the Applied Sciences, {\bf 41} (2018), pp. 1--10.

\bibitem{PABB}
D. Peaucelle, D. Arzelier, O. Bachelier, J.Bernussou, {\it A new robust $\mathcal D$-stability condition for real convex polytopic uncertainty,} System \& Control letters, {\bf 40} (2000), pp. 21--30.

\bibitem{PER}
S.K. Persidskii, {\it Problem of absolute stability,} Automat. Remote Control, {\bf 12} (1969), pp. 1889--1895.

\bibitem{PONT}
 L. Pontrjagin, {\it Topological groups,} Princeton University Press, 1946.

\bibitem{PRY1}
O. Pryporova, {\it Types of convergence of matrices,} PhD Thesis, Iowa State University, Ames, Iowa, 2009.

\bibitem{PRY2}
O. Pryporova, {\it Qualitative convergence of matrices,} Linear Algebra Appl. {\bf 431} (2009), p. 28-38.

\bibitem{QR}
J.P. Quirk, R. Ruppert, {\it Qualitative economics and the stability of equilibrium,} Rev. Econom. Studies {\bf 32} (1965), 311-326.

\bibitem{QS}
J.P. Quirk, R. Saposnik, {\it Introduction to general equilibrium theory and welfare economics,} McGraw-Hill, New-York, 1968.

\bibitem{RAS}
Q.I. Rahman, G. Schmeisser, {\it Analytic theory of polynomials,} Clarendon Press, Oxford, 2002.

\bibitem{REDH1}
R. Redheffer, {\it Volterra multipliers I,} SIAM J. Alg. Disc. Meth., {\bf 6} (1985), 592--611.

\bibitem{REDH2}
R. Redheffer, {\it Volterra multipliers II,} SIAM J. Alg. Disc. Meth., {\bf 6} (1985), 612--623.

\bibitem{RS1}
I.M. Romanishin, L.A. Sinitskii, {\it On the additive $D$-stability of matrices on the basis of the Kharitonov criterion,} Mathematical Notes, {\bf 72} (2002), pp. 237-240.

\bibitem{ROU1}
E.J. Routh, {\it The advanced part of a treatise on the dynamics of a system of rigid bodies}, Macmillan and Co., London, 1884, 4th ed., pp. 168--176.

\bibitem{ROU2}
E.J. Routh, {\it A treatise on the stability of a given state of motion}, Macmillan and Co., London, 1877.

\bibitem{ROU3}
E.J. Routh, {\it The Advanced Part of a Treatise on the Dynamics of a System of Rigid Bodies. Being Part II of a Treatise on the Whole Subject. With Numerous Examples}, Dover Publications, 1955.

\bibitem{ROU4}
E.J. Routh, {\it Dynamics of a system of rigid bodies}, Macmillan, 1892.

\bibitem{RUGH}
W.J. Rugh, {\it Linear system theory}, Prentice-Hall, 1996.

\bibitem{SADR}
R.A. Satnoianu, P. van den Driessche, {\it Some remarks on matrix stability with application to Turing instability}, Linear Algebra Appl. {\bf 398} (2005), 69--74.

\bibitem{SCHNE}
H. Schneider, {\it Positive operators and an inertia theorem,} Numerische Matematik, {\bf 7} (1965), pp. 11--17.

\bibitem{SHU}
I. Schur, {\it Bemerkungen zur Theorie der beschr\"{a}nkten Bilinearformen mit unendlich vielen Ver\"{a}nderlichen}, J. Reine Angew. Math., {\bf 140} (1911), pp. 1--28.

\bibitem{SHU2}
I. Schur, {\it \"{U}ber Potenzreihen, die im Innern des Einheitskreises beschr\"{a}nkt sind}, J. Reine Angew. Math., {\bf 147} (1917), pp. 205-232.

\bibitem{SJK}
 R. Sepulchre, M. Jankovic, P.V. Kokotovic, {\it Constructive nonlinear control}, Springer (2011).

\bibitem{SMK}
 R. Shorten, O. Mason, C. King, {\it An alternative proof of the Barker, Berman, Plemmons (BBP) result on diagonal stability and extensions}, Linear Algebra Appl., {\bf 430} (2009), pp. 34--40.

\bibitem{SHN}
 R. Shorten, K.S. Narendra, {\it Diagonal stability and strict positive realness}, Proceedings of IEEE Conference on Desicion and Control, 2006.

\bibitem{SHN2}
 R. Shorten, K.S. Narendra, {\it On a theorem of Redheffer conserning diagonal stability}, Linear Algebra Appl., {\bf 431} (2009), pp. 2317--2329.

\bibitem{SHN3}
 R. Shorten, K.S. Narendra, {\it On common quadratic Lyapunov functions for pairs of stable LTI systems whose system matrices are in companion form}, IEEE Transactions on Automatic Control, {\bf 48} (2003), pp. 618--621.

\bibitem{SHN4}
 R. Shorten, K.S. Narendra, {\it On the diagonal stability of a class of almost positive switched systems}, Proceedings of 2010 American Control Conference, Baltimore, USA (2010), pp. 6250--6255.

\bibitem{SIN}
 R. Sinkhorn, {\it Matrices which commute with doubly stochastic matrices}, Linear and Multilinear Algebra, {\bf 4} (1976), pp. 201--203.

\bibitem{SOF}
 Y.C. Soh, Y.K. Foo, {\it Kharitonov regions: it suffices to check a subset of vertex polynomials}, IEEE Transactions on Automatic Control, {\bf 36} (1991), pp. 1102--1105.

\bibitem{STE}
 P. Stein, {\it Some general theorems on iterants}, J. Res. Nat. Bur. Standards, {\bf 48} (1952), pp. 82--83.

\bibitem{STU}
 C. Sturm, {\it Analyse d'un m\'{e}moire sur la r\'{e}solution des \'{e}quations num\'{e}riques}, Bull. Sci. Math. Ferussac., {\bf II} (1829), pp. 419--422.

\bibitem{SGH}
 Y.-J. Sun, R.-S. Gau, J.-G. Hsieh, {\it Simple criteria for sector root clustering of uncertain systems with multiple time delays}, Chaos, Solutions and Fractals, {\bf 39} (2009), pp. 65--71.

\bibitem{SUN}
 W. Sun, {\it The interior and closure of strongly stable matrices}, Linear Algebra Appl., {\bf 165} (1992), pp. 53--58.

\bibitem{TAU}
O. Taussky, {\it A generalization of a theorem of Lyapunov},
J. Soc. Indust. Appl. Math., {\bf 9} (1961), p. 640--643.

\bibitem{TAU1}
O. Taussky, {\it A remark on a theorem of Lyapunov},
Journal of Math. Analysis Appl., {\bf 2} (1961), p. 105--107.

\bibitem{TAU2}
O. Taussky, {\it Matrices $C$ with $C^n \rightarrow 0$},
Journal of Algebra, {\bf 1} (1964), p. 5--10.

\bibitem{THR1}
C.D. Thron, {\it The secant condition for instability in biochemical feedback control --- I. The role of cooperativity and saturability},
Bulletin of Mathematical Biology, {\bf 53} (1991), p. 383--401.

\bibitem{THR2}
C.D. Thron, {\it The secant condition for instability in biochemical feedback control --- II. Models with upper Hessenberg Jacobian matrices},
Bulletin of Mathematical Biology, {\bf 53} (1991), p. 403--424.

\bibitem{TIM}
D. Timotin, {\it Redheffer products and characteristic functions},
Journal of Mathematical Analysis Appl., {\bf 196} (1995), p. 823--840.

\bibitem{TOG}
Y. Togawa, {\it A geometric study of the $D$-stability problem},
Linear Algebra Appl., {\bf 33} (1980), p. 133--151.

\bibitem{TK}
 W. Tru\"{o}l, F.J. Kraus, {\it Robust $D$-stability in frequency domain with Kharitonov-like properties},
in: IFAC Design Methods of Control Systems, Zurich (1991), pp. 149-154.

\bibitem{TSA}
 M. Tsatsomeros, {\it Generating and detecting matrices with positive principal minors},
Asian Information-Science-Life, {\bf 1} (2002), p. 115--132.

\bibitem{TYO}
J.J. Tyson, H.G. Othmer, {\it The dynamics of feedback control circuits in biochemical pathways},
In R. Rosen and F. Snell, {\it Progress in theoretical byology,} Academic Press, New York {\bf 5} (1978), p. 1--62.

\bibitem{TSOY}
T.-W. Ma, {\it Classical analysis on normed spaces}, World Scientific Publishing, 1995.

\bibitem{TZO}
 M. Tzoumas, {\it On sign-symmetric circulant matrices},
Applied Math. and Computations, {\bf 195} (2008), p. 604--617.

\bibitem{WAN}
 M. Wanat, {\it The $\alpha$-scalar diagonal stability of block matrices},
Linear Algebra Appl., {\bf 414} (2006), p. 304--309.

\bibitem{WANGL}
L. Wang, M.Y. Li, {\it Diffusion-driven instability in reaction-diffusion systems},
Journal of Math. Analysis and Appl., {\bf 254} (2001), p. 138--153.

\bibitem{WAR}
E. Waring, {\it Problems}, Philos, Trans. Roy. Soc. London, {\bf 53} (1763), p. 294--299.

\bibitem{WIM}
 H.K. Wimmer, {\it Generalizations of theorems of Lyapunov and Stein},
Linear Algebra Appl., {\bf 10} (1975), p. 139--146.

\bibitem{WIM3}
 H.K. Wimmer, {\it Diagonal stability of matrices with cyclic structure and the secant condition},
System and Control Letters, {\bf 58} (2009), p. 309--313.

\bibitem{WOOD}
J.E. Woods, {\it Mathematical economics. Topics in multi-sectoral economics}, Longman Group Ltd, London, 1978.

\bibitem{YAR}
 R. Yarlagadda, {\it Stabilization of matrices},
Linear Algebra Appl., {\bf 21} (1978), p. 271--288.

\bibitem{ZAL}
N. Zavalishin, D. Logofet, {\it Modelling ecological systems according to a given "storage-flow" diagram},
Mathematical Models and Computer Simulations {\bf 9} (1997), pp. 3-17 (in Russian).
\end{thebibliography}
\end{document}